\documentclass[12pt,reqno]{article} 

\usepackage[left=1cm,right=2cm,top=2cm,bottom=2cm]{geometry}

\usepackage{amsmath,amssymb,amsthm,amsfonts}

\usepackage{graphicx}   
\usepackage{float}      
\usepackage{booktabs}   

\usepackage{algorithm}
\usepackage{algpseudocode}

\usepackage{titlesec}
\usepackage{tocloft}

\usepackage[utf8]{inputenc}
\usepackage{xcolor}
\usepackage{hyperref}
\usepackage{cleveref}

\usepackage{titlesec} 

\titleclass{\subsubsubsection}{straight}[\subsubsection]

\newcounter{subsubsubsection}[subsubsection]
\renewcommand\thesubsubsubsection{\thesubsubsection.\arabic{subsubsubsection}}

\titleformat{\subsubsubsection}
  {\normalfont\normalsize\bfseries}{\thesubsubsubsection}{1em}{}

\titlespacing*{\subsubsubsection}
  {0pt}{3.25ex plus 1ex minus .2ex}{1.5ex plus .2ex}

\makeatletter
\newcommand\l@subsubsubsection{\@dottedtocline{4}{7em}{4em}}
\makeatother

\usepackage[font={small}]{caption}
\newcommand{\mathsym}[1]{{}}

\usepackage{color}
\usepackage[utf8]{inputenc}
\usepackage{amsmath,amsthm,amsfonts}  
\usepackage{sectsty}  
\usepackage{graphicx} 
\usepackage{hyperref} 
\usepackage{cleveref}

 \author{Noura Al Helwani \thanks{Department of Physics, American University of Beirut, Lebanon (nmh82@mail.aub.edu)} \and Sophie Moufawad \thanks{Department of Mathematics, American University of Beirut, Lebanon (sm101@aub.edu.lb)} \and Georges Sakr \thanks{Department of Mathematics, American University of Beirut, Lebanon (gcs01@mail.aub.edu)} \thanks{The authors would like to acknowledge: \\
$\hspace{2mm}\qquad \qquad \bullet$ The contribution of  Mr Christophe El Chabab, and Mr. Ali Ibrahim during the Summer Research Camp 2024\\
$\qquad$ $\normalsize \bullet$ AUB Center for Advanced Mathematical Sciences (CAMS) and Math department's support 
} }
\title{Solving Inverse PDE Problems using Minimization Methods and AI } 
\begin{document}
\maketitle
\begin{abstract}
\noindent Many physical and engineering systems require solving direct problems to predict behavior and inverse problems to determine unknown parameters from measurement. In this work, we study both aspects for systems governed by differential equations, contrasting well-established numerical methods with new AI-based techniques, specifically Physics-Informed Neural Networks (PINNs). We first analyze the logistic differential equation, using its closed-form solution to verify numerical schemes and validate PINN performance. We then address the Porous Medium Equation (PME), a nonlinear partial differential equation with no general closed-form solution,  building strong solvers of the direct problem and testing techniques for parameter estimation in the inverse problem. Our results suggest that PINNs can closely estimate solutions at competitive computational cost, and thus propose an effective tool for solving both direct and inverse problems for complex systems.
\end{abstract}
\hrule

\tableofcontents

\vspace{0.75cm}

\hrule

\section{Introduction}
Compared to direct problems, which start with a mathematical model and search for its solution, inverse problems start in the opposite direction, beginning from measurements or data and attempting to infer unknown parameters or inputs of the underlying model. Such problems appear everywhere in applied mathematics because many real-world challenges are essentially inverse problems. For example, in oil and gas exploration, seismic waves are pumped into the earth, and the reflected waves are analyzed to make inferences about rock types and fault locations; this is a classic inverse problem.

\noindent For a few more mathematically formal examples, paper \cite{Moufawad} addresses an inverse problem of gas diffusion in polar firn, modeled by a parabolic partial differential equation. The coefficients are all assumed to be constant, and the diffusion coefficient is assumed to be unknown and reconstructed using generated data. Study \cite{Kaltenbacher} concerns the inversion of an inverse problem in a reaction-diffusion system to recover source terms, which are functions of the state variables, using extra information such as values of the state at some time or time measurements at a fixed boundary point. Lastly, research \cite{Caro} addresses an inverse quantum mechanics problem, in which the task is to reconstruct an unknown Hamiltonian operator controlling particle dynamics from initial and final quantum states of particles that evolve under the same Hamiltonian.

\noindent Inverse problems are generally stated in the form of ordinary differential equations (ODEs) and partial differential equations (PDEs) and are usually cast as an optimization or minimization problem. While there are exact solutions-either in closed form, integral form, or series form-for some linear first-order and second-order ODEs, most ODEs do not have exact solutions. More importantly, PDEs very rarely have closed-form solutions.

\noindent Numerical techniques such as finite difference, finite element, and spectral methods are often used for solving PDEs. These techniques approximate the continuous problem by discretizing the domain and mapping the equations to systems of algebraic equations that can be solved computationally. There exist, however, several challenges in numerical methods for solving PDEs such as discretization errors, high mesh-resolution requirements, and instability in iterative schemes. To explain, as discretization approximates continuous spatial or temporal spaces by finite points, it adds errors that can vitiate accuracy unless the mesh is refined enough, but refinement involves increased computational expense. In situations involving complicated geometries or 3D simulations, creating a high-resolution mesh can have an element size in the billions, which significantly necessitates the use of high-performance computing. Also, the majority of the numerical solvers are based on iterative techniques to refine an approximation further; these iterations  may be unstable and slowly converge. Along comes Physics Informed Neural Networks (PINNs) to the rescue. PINNs offer a mesh-free alternative to traditional numerical methods by applying physics directly through neural network loss functions. The primary objective of this research is to discover the PINNs' methodology and demonstrate both their efficacy for finding the solutions for the forward and inverse problems of PDEs and their restrictions.
\noindent We will consider 2 different Inverse problems:
\begin{enumerate}
\item The logistic equation, 
\begin{equation}\label{eq:1}
\begin{cases}
p'(t) = (1 - \dfrac{p(t)}{K})\;r\;p(t)& t\geq t_0\\
p(t_0) = p_0
\end{cases}
\end{equation}
where the goal is to find the growth rate $r$ of some population, given the population size at some time $t_i$, $p(t_i)$, the initial population size $p_0$, and the environmental carrying capacity or
saturation level $K$ (available food, space,..). It is assumed that $r$ is a constant. 
  
\item The one-dimensional porous medium equation (PME),
\begin{equation}
\begin{cases}
u_t(t,x) - \Delta(|u|^{\beta-1}u) = 0\quad & t> t_0, \; x\in (a,b)\\
u(t,a) = 0 & t> t_0\\
u(t,b) = 0 & t> t_0\\
u(t_0,x) = u_0(x)&x\in (a,b) \\
\end{cases}
\end{equation}
 where the goal is to find the  polytropic exponent $\beta$, given $u(t_i) \geq 0$, the scaled density of the considered gas  at different times $t_i$, and $u_0(x)$, the initial scaled density.

The porous medium equation (PME) is used in different fields to model different phenomena such as the diffusion of an ideal gas in porous media, nonlinear heat transfer, groundwater flow, population dynamics,..

\end{enumerate}
\noindent As a benchmark, we begin with the logistic differential equation, using its known analytical solution to validate our computational approach. However, the main focus is the porous medium equation (PME), a non-linear partial differential equation with no general closed-form solution. We develop and test several methods for parameter recovery and evaluate their accuracy and computational efficiency using real and generated data (with and without added noise).
In Chapter 2, we start by applying some classical methods to solve the direct and inverse problems of Equations 1 and 2. Then, in Chapter 3, we move to Physics-Informed Neural Networks to tackle the same equations. Throughout the process, we compare methods, analyze results, and study their strengths and weaknesses.

\section{Classical Methods}

Section~2 recapitulates the classical numerical toolkit that will be our reference point for the subsequent AI-based comparisons. In Section~2.1 we revisit the logistic equation: its forward (direct) dynamics are solved by high-order Runge-Kutta schemes and MATLAB built-in solvers (\texttt{ode45}), whereafter the respective inverse problem is cast as an optimization problem solved by Newton's method, the secant method, steepest descent with Armijo rule, and MATLAB's \textit{fminunc}/\textit{fmincon} functions. {In Section~2.2} we proceed with the one-dimensional porous-medium equation, where finite-difference discretizations in space and time combined with minimization techniques enable parameter recovery in a nonlinear PDE setting. Accuracy and computational cost are evaluated along the way to establish a clear performance benchmark for the modern methods introduced in the next section. Test results for the inverse problem for the logistic equation problem are in the Appendix \ref{sec:app}.

\subsection{Logistic Equation}\label{sec:Log}
Ideally, population growth would be exponential, with the population increasing by a fixed percentage. This would require unlimited resources to meet everyone's needs. This situation can be described by the following simple ODE:

\[ \frac{dp(t)}{dt} = rp(t) \]
 where p(t) $>$ 0 is the population and r $>$ 0 is the rate of growth. 
However, in reality, resources are limited, and many factors affect population growth. As the population approaches the carrying capacity K $>$ 0, resources become scarcer, and the growth rate slows down until it eventually stops. 

\noindent On the other hand, when resources are unlimited (i.e K can go to infinity), the situation is again ideal and growth should be exponential.
To satisfy these two cases, one can multiply the exponential growth equation by \((1 - \frac{p}{K})\) and get

\begin{equation} \label{eq: logistic}
\begin{cases}
\dfrac{dp(t)}{dt} = (1 - \dfrac{p(t)}{K})\;r\;p(t),& t\geq t_0\\
p(t_0) = p_0
\end{cases}
\end{equation}

\noindent This equation is the logistic equation and it expresses one model of the real-world conditions.
In fact, one can study the population's variations before solving the ODE.
Given that $r$ and $p(t)$ are positive, $ \dfrac{dp(t)}{dt}$ takes the sign of $(1 - \dfrac{p(t)}{K})$ for $t\geq t_0$:
\begin{enumerate}
 \item Underpopulation: if $p(t) < K$, then $ \dfrac{dp(t)}{dt}$ is positive and $p(t)$ increases and tends to K as t goes to infinity.
 \item Equilibrium point: if $p(t) = K$, then $ \dfrac{dp(t)}{dt}$ is zero and $p(t)$ stays constant.
 \item Overpopulation: if $p(t) > K$, then $ \dfrac{dp(t)}{dt}$ is negative and $p(t)$ decreases and tends to K as t goes to infinity.
 \end{enumerate}
 Turning our attention to the factors that can influence population growth and the carrying capacity K. Several models in the literature depict different factors such as land availability \cite{Land}, economic growth and resource consumption \cite{AT1}, climate change \cite{AT2}, fertility \cite{BTT}, and war and armed conflicts \cite{Brezis} among others.

\subsubsection{Direct Problem}
In this section, we address the direct problem of the logistic equation, aiming to find the solution of Equation~\ref{eq: logistic}. We begin by solving it analytically, as it is a separable first-order ODE; this is presented in Section~\ref{sec:logistic analytical sol}. 

\noindent Next, we turn to numerical methods for solving the ODE. We briefly introduce the Runge-Kutta method of order 4 (RK4) and MATLAB's \texttt{ode45} in Sections~\ref{sec:RK4} and~\ref{sec:ODE45}, respectively. Finally, in Section~\ref{sec:test}, we test these methods and compare their results to the exact analytical solution.

\subsubsubsection{Analytical Solution} \label{sec:logistic analytical sol}
Given values of $r$ and $k$, we can determine the future population at a specific time $t$ by solving the logistic ODE using the method of separation of variables:
$$p'(t) = r p(t) \left(1 - \frac{p(t)}{K}\right), \qquad \quad \implies  r = \frac{p'(t)}{\left(1 - \dfrac{p(t)}{K}\right) p(t)}, \qquad \quad \implies  C e^{r t} = \frac{p(t)}{K - p(t)}\vspace{-6mm} $$
Thus, $ p(t) = \dfrac{C K e^{r t}}{1 + C e^{r t}}$ is the solution to the ODE.

\noindent To determine the constant \( C \), we apply the initial condition \( p(t_0) = p_0 \):
$$  \frac{C K e^{r t_0}}{1 + C e^{r t_0}} = p_0, \qquad \quad \implies  C = \frac{p_0}{(K - p_0) e^{r t_0}}$$
We obtain the following exact analytical solution to the logistic initial value problem:
\begin{eqnarray}
p(t) &=& \dfrac{Kp_0e^{r(t-t_0)}}{K-p_0+p_0e^{r(t-t_0)}} \label{eq:17}
\end{eqnarray}
\subsubsubsection{RK4}\label{sec:RK4}
The Runge-Kutta 4th order method (RK4) is a widely used numerical technique for solving first order ODEs  with initial conditions of the form $$\begin{cases}
   y'(t) = f(t,y(t)), & t>t_0\\
    y(t_0) = y_0&
\end{cases}$$ due to its balance between accuracy and computational efficiency. It estimates the solution by computing four intermediate slopes ($k_1, k_2, k_3, k_4$) at each time step and then combining them using a weighted average to approximate the next value as detailed in Algorithm \ref{alg:RK4}. The method achieves fourth-order accuracy, meaning the global error scales with the fourth power of the time step size $h$ ($\mathcal{O}(h^4)$, with $h \in (0,1]$), making it more precise than Euler (RK-1; $\mathcal{O}(h)$) or midpoint methods (RK2-M; $\mathcal{O}(h^2)$) for the same step size. The RK4 method can be seen as a discrete approximation to the Taylor series expansion of the exact solution without requiring all higher-order derivatives. \vspace{-2mm}
\begin{algorithm}[H]
\centering
\caption{Runge-Kutta 4th order (RK4)}\label{alg:RK4}
{\renewcommand{\arraystretch}{1.3}
\begin{algorithmic}[1]
\Statex \textbf{Input:} \quad $f$: function handle (population growth); $t$: time vector; $y_0$: initial value; $h$: time step; $r$: function handle for $r(t)$ (can be constant or time-dependent).
\Statex \textbf{Output:} $y$: vector of computed values. \vspace{2mm}

\Statex \texttt{function y = RK4Logistic(f, t, y0, h, r)}
\State $n = \texttt{length}(t); \quad y = \texttt{zeros}(n, 1)$
\State $y(1) = y_0$
\For {$i = 1 : n-1$}
    \State $t_i = t(i)$; \quad $t_{ih} = t_i + h$; \quad $t_{ih2} = t_i + \frac{h}{2}$
    \State $k_1 = f(t_i, r(t_i), y(i))$
    \State $k_2 = f(t_{ih2}, r(t_{ih2}), y(i) + 0.5*h*k_1)$
    \State $k_3 = f(t_{ih2}, r(t_{ih2}), y(i) + 0.5*h*k_2)$
    \State $k_4 = f(t_{ih}, r(t_{ih}), y(i) + h*k_3)$
    \State $y(i+1) = y(i) + \frac{h}{6} (k_1 + 2k_2 + 2k_3 + k_4)$
\EndFor
\State \textbf{end}
\end{algorithmic}}
\end{algorithm}
\subsubsubsection{ODE45}\label{sec:ODE45}
\texttt{ODE45} is a MATLAB function that solves ODEs using an adaptive Runge-Kutta method of orders 4 and 5 (combining 4th and 5th-order formulas to estimate both the solution and the local error). The local error estimate is used to adjust the step size $h$: if the estimated error is too large, the step is reduced; if it is sufficiently small, the step is increased. This dynamic adjustment allows \texttt{ODE45} to maintain accuracy while minimizing the number of function evaluations needed to compute the solution efficiently.

\subsubsubsection{Testings}\label{sec:test}
We now compare the performance of RK4 (Algorithm~\ref{alg:comprk4})  and ODE45 (Algorithm~\ref{alg:compode45}) for solving the logistic ODE under various $r$, $K$, $p_0$, and time intervals. The goal is to assess both accuracy and computational behavior of fixed-step versus adaptive solvers.

\begin{algorithm}[H]
\centering
\caption{Comparison of RK4 Numerical Solution with Analytical Solution for the Logistic ODE}\label{alg:comprk4}
{\renewcommand{\arraystretch}{1.3}
\begin{algorithmic}[1]
\Statex{\textbf{Input:} \;\;\;$f$: function handle; $t_0$: initial time; $T$: end time; $y_0$: initial value; $n$: number of steps; $r$: function handle for $r(t)$; exactSol: analytical solution}
\Statex{\textbf{Output:} $p$: vector of computed values; $relerr$: average relative error between numerical and analytical solution}\vspace{2mm}
\Statex {function [p, relerr] = CompareLogisticRK4 ($f, t_0, T, p_0, n, r, exactSol$)}
\State $h = \frac{T - t_0}{n}$;
\State $t = t_0 : h : T$;
\State $tresp = exactSol(t)$;
\State $p = \texttt{RK4Logistic}(f, t, p_0, h, r)$;
\State $relerr = \dfrac{\texttt{norm}(p - tresp)}{\texttt{norm}(tresp)}$;
\State $relerr = \dfrac{relerr}{n+1} \qquad \triangleright\text{\textit{average relative error}}$;
\State \textbf{end}
\end{algorithmic}}
\label{alg:Runge-Kutta-Euler}
\end{algorithm}

\begin{algorithm}[H]
\centering
\caption{Comparison of \texttt{ode45} Numerical Solution with Analytical Solution for the Logistic ODE}\label{alg:compode45}
{\renewcommand{\arraystretch}{1.3}
\begin{algorithmic}[1]
\Statex \textbf{Input:} \;\;\;$f$: function handle; $t_0$: initial time; $T$: end time; $y_0$: initial value; $n$: number of steps; $r$: function handle for $r(t)$; $sol$: analytical solution
\Statex \textbf{Output:} $p$: vector of computed values; $relerr$: average relative error between numerical and analytical solution{\vspace{2mm}}
\Statex {function [p, relerr] = CompareLogisticOde45 ($f, t_0, T_{}, y_0, r, n, exactSol$)}
\State $(t_{ODE}, p) = \texttt{ode45}(f, [t_0, T_{}], y_0)$;
\State $t = \texttt{linspace}(t_0, T_{}, n)$;
\State $analytical\_sol = sol(t)$;
\State $relerr = \dfrac{\texttt{norm}(p - analytical\_sol)}{\texttt{norm}(p) \times \texttt{length}(t)} \qquad \triangleright  \text{\textit{average relative error}}$;
\end{algorithmic}}
\label{alg:ODE-Solver}
\end{algorithm}
\noindent The logistic equation under consideration is given by:
\[
f = @(t, r, y)\, r \cdot y \cdot \left(1 - \frac{y}{K} \right)
\]
with the exact analytical solution:
\[
\texttt{exactSol} = @(t, r, y, p_0)\, \frac{K \cdot p_0 \cdot \exp(r \cdot (t - t_0))}{K - p_0 + p_0 \cdot \exp(r \cdot (t - t_0))}
\]

We apply both methods across different scenarios by varying the time interval, number of steps, and the values of $r$, $K$, and $p_0$. The computed solutions are then compared to the analytical solution using the average relative error described in the algorithms.

\noindent The table below summarizes the results obtained using RK4 and \texttt{ode45} under each parameter setting.

\begin{table}[h] 
    \centering

    \begin{tabular}{|c|c|c|c|c|c|c|}
        \hline
        Time Interval & Intervals & $K$ & $p_0$ & $r$ &  RK4 (RelErr) & ODE45 (RelErr) \\
        \hline
        2011 - 2022 & 100 & 10 & 20 & 0.079 &$4.1640\times 10^{-3}$&$9.1448 \times 10^{-8}$ \\
        450 - 500 & 200 & 90 & 10 & 0.05 & $9.8764\times 10^{-3}$&$5.6998 \times 10^{-8}$\\
        1 - 100 & 100000 & $10^3$ & $10^2$ & 0.9 & $5.4090 \times 10^{-5}$ &$3.5415 \times 10^{-4}$\\ 
        \hline
    \end{tabular}
    \caption{Relative errors (RelErr) for RK4 and ODE45 with different parameters.}
    \label{table: rlt errors}
\end{table}

\noindent Both RK4 and ODE45 accurately solve the logistic ODE, with ode45 consistently achieving higher precision due to its adaptive step size control. However, this increased accuracy comes at a higher computational cost per step. In Algorithm ~\ref{alg:RK4}, the classical RK4 method requires four evaluations of the right-hand side per time step. In contrast, ode45 is based on the seven-stage Dormand-Prince embedded Runge-Kutta formula, which computes additional stages to enable local error estimation and adaptivity \cite{ShampineReichelt1997}. Consequently, RK4 can be more efficient for problems where fixed step sizes are sufficient. 

In practice, RK4 is preferred for its simplicity and speed when the problem is well-behaved, while ode45 is better suited for scenarios requiring high accuracy or involving rapidly changing dynamics. The choice ultimately depends on the desired trade-off between computational efficiency and solution precision.

\subsubsection{Inverse Problem}
Realistically, parameters such as r, k, and $p_0$ can not be easily calculated. One could try to find them based on observations of the current and previous populations. This consists of solving the inverse problem, whose objective is to find parameters that best minimize the cost or loss function defined using the sum of squared errors:
$$f(u)= \sum_{i=1}^m(p(t_i) - p_i)^2$$ where $p_i$ is the observed population  at a given times $t_i$ taken from existing data and $p(t_i)$ is the computed solution for the parameter (or set of parameters) $u$.
In our case, we consider the following optimization scenarios:
\begin{enumerate}
    \item Minimize $f$ with respect to $r$, where $u = r$
    \item Minimize $f$ with respect to $r$ and $k$, where $u = (r, k)$
    \item Minimize $f$ with respect to $r$, $k$, and $p_0$, where $u = (r, k, p_0)$
\end{enumerate}

\noindent In this section, we focus on the first two cases. 
In Section~\ref{sec: inv prob analyt logistic}, we address the problem analytically, while in Section~\ref{sec: mini algo}, we approach it using numerical minimization algorithms.

\subsubsubsection{Analytical Solution }\label{sec: inv prob analyt logistic}
Let us derive an analytical expression for $r$. We assume that $r$ is a function of time, taking the value $r_i$ at $t = t_i$. Let $p_r(t_i)$ denote the approximation of the population at time $t_i$ when using $r = r_i$. We proceed as follows.

\begin{equation}\label{eq:28}
\frac{df(r)}{dr} = 0 \qquad \iff \qquad
\frac{d(p_r(t_i) - p_i)^2}{dr} = 0 \qquad \iff \qquad
2(p_r(t_i) - p_i) \cdot \frac{\partial p_r(t_i)}{\partial r} = 0 
\end{equation}

Taking the partial derivative with respect to $r$ of the analytical solution, we get:
\begin{equation}
\frac{\partial p_r(t_i)}{\partial r} = \frac{K p_0 (t_i - t_0)(K - p_0) e^{r(t_i - t_0)}}{(K - p_0 + p_0 e^{r(t_i - t_0)})^2}
\end{equation}

Substituting this into Equation~\eqref{eq:28} and writing $p_r(t_i)$ explicitly, we get:
\begin{eqnarray}
2\left( \frac{K p_0 e^{r(t_i - t_0)}}{K - p_0 + p_0 e^{r(t_i - t_0)}} - p_i \right) 
\cdot 
\left( \frac{K p_0 (t_i - t_0)(K - p_0) e^{r(t_i - t_0)}}{(K - p_0 + p_0 e^{r(t_i - t_0)})^2} \right)
&=& 0
\end{eqnarray}
Solving for $r$, we get:\vspace{-2mm}
\begin{equation}
    r_i = \dfrac{ln({\dfrac{P_i(K-p_0)}{p_0(K-p_i)}})}{t_i-t_0} 
\end{equation}

\noindent However, this results in an interpolation error of zero (see Appendix~\ref{app:analytical-logistic}). Nonetheless, our primary interest lies in prediction, and interpolation does not guarantee low extrapolation error. Indeed, we are at risk of overfitting on our training data, and at danger of low generalization to unseen data.\\[1mm]

\noindent This observation motivates the use of numerical minimization algorithms to estimate $r$. In what follows, we assume $r$ is constant and test its estimation using simulated data. Our objective in this section is to solve the first two inverse problems, namely:
finding $r$, or $(r, K)$, such that the loss function is minimized. We begin by presenting several minimization algorithms. Then, we numerically approximate $r$ using two datasets: first, data generated using the exact analytical solution, and second, real-world population data.

\subsubsubsection{Minimization Algorithms} \label{sec: mini algo}
Numerical methods for function minimization are computational techniques used to iteratively determine the point at which a given function attains its lowest value. These methods are particularly valuable when analytical solutions are complex or unavailable. In the following sections, we present several algorithms used in such one-variable minimization problems.\\[0.5em] The first step involves approximating the derivative of a function symbolically in MATLAB.\\

\begin{algorithm}[H]
\centering
\caption{Symbolic Derivative Calculation}
{\renewcommand{\arraystretch}{1.3}
\begin{algorithmic}[1]
\Statex{\textbf{Input:} $f$: function handle of one variable}
\Statex{\textbf{Output:} $df$: symbolic approximation of the derivative of $f$ \vspace{3mm}}
\Statex   function $df$ =symbolicDerivative($f$)
\State $fsym = \text{sym}(f)$ \Comment{Converts a function handle into a symbolic function}
\State $dfsym = \text{diff}(fsym)$ \Comment{Computes the derivative symbolically}
\State $df = \text{matlabFunction}(dfsym)$ \Comment{Converts a symbolic function back to a function handle}
\end{algorithmic}}
\label{alg:SymbDeriv}
\end{algorithm}

\noindent We aim to minimize a function , which is done by finding the roots of its derivative. To do so, we first construct a symbolic derivative using the above procedure. The resulting $df$ function can then be used as an input for root-finding algorithms, such as Newton's method and the Secant method. The Newton's and Secant methods are used on the derivative to find its roots, whereas \texttt{fminunc}, \texttt{fmincon} and the steepest descent are applied to the function itself. 
\subsubsubsection*{Newton's Method}\label{section: newtons}
Newton’s method is a classical iterative technique used to approximate the roots of a real-valued function. It works by linearizing the function near a current estimate and using the tangent line to approximate where the function crosses the x-axis. More precisely, given a current approximation \( x_n \), the method updates it using the formula:
\[
x_{n+1} = x_n - \frac{f(x_n)}{f'(x_n)}
\]
This formula is derived from the \textit{first-order Taylor expansion} of the function:
\[
f(x) \approx f(x_n) + f'(x_n)(x - x_n)
\]
As we previously mentioned, we approximate the function where it crosses the x-axis, so setting this approximation to zero and solving for \( x \) yields the update rule. Newton’s method works well when the function is smooth and the initial guess is sufficiently close to the actual root. In such cases, it exhibits quadratic convergence, meaning the error roughly squares at each iteration, leading to very rapid convergence. However, the method is local in nature; in fact, poor initial guesses can lead to divergence, oscillations, or convergence to the wrong root. Additionally, it requires knowledge of the function’s derivative, which may be challenging to compute. Despite these challenges, Newton’s method remains a fundamental and powerful tool in numerical root-finding. The algorithm below implements Newton’s method using a symbolic derivative. Note that we can replace the symbollic derivative by the exact analytical derivative if it is known to us.
\begin{algorithm}[H]
\centering
\caption{Newton's Method for Root Finding}
{\renewcommand{\arraystretch}{1.3}
\begin{algorithmic}[1]
\Statex{\textbf{Input:} \;\;\;$f$: function handle (derivative function from previous symbolic derivative algorithm); $i$: initial condition; $n$: max number of Newton iterations; $tol$: error tolerance}
\Statex{\textbf{Output:} $sol$: approximation of the root \vspace{3mm}}
\Statex function $sol$ = Newton($f, i, n, tol$)
\State $error = 1 + tol$;  
\quad $r(1) = i$ 
\State $df = \text{Symb derv}(f)$
\State $count = 1$ 
\While{$count \leq n$ and $error > tol$}
    \State $r(count + 1) = r(count) - {f(r(count))}/{df(r(count))}$
    \State $error = {\text{norm}(r(count + 1) - r(count))}/{\text{norm}(r(count))}$
    \State $count = count + 1$
\EndWhile
\If{$error < tol$}
    \State $sol = r(\text{length}(r) - 1)$
\Else
    \State $sol = \text{``No convergence"}$
\EndIf
\end{algorithmic}}
\label{alg:NewtonRootFinding}
\end{algorithm}

\subsubsubsection*{Secant Method}
Similarly, we will input the derivative function in the Secant algorithm:\\
\noindent The Secant method is a root-finding algorithm closely related to Newton’s method, but it eliminates the need for computing the derivative (analytically or even symbolically). Instead, it approximates the derivative using the slope of a secant line through two recent approximations of the root. The update formula is given by
\[
x_{n+1} = x_n - f(x_n) \cdot \frac{x_n - x_{n-1}}{f(x_n) - f(x_{n-1})}
\]
This approach makes the Secant method particularly useful when the derivative of the function is difficult or expensive to compute. Although it typically exhibits superlinear convergence (faster than linear but slower than Newton’s quadratic rate), it remains efficient in practice due to its lower per-iteration computational cost. Like Newton’s method, it is also local in nature and sensitive to the choice of initial guesses. Nevertheless, the Secant method strikes a good balance between accuracy and ease of implementation when derivative information is unavailable or unreliable.

\begin{algorithm}[H]
\centering
\caption{Secant Method for Root Finding}
{\renewcommand{\arraystretch}{1.3}
\begin{algorithmic}[1]
\Statex{\textbf{Input:} \;\;\;$f$: function handle; $i_1$: first initial condition; $i_2$: second initial condition; $n$: max number of iterations; $tol$: error tolerance}
\Statex{\textbf{Output:} $sol$: approximation of the root \vspace{3mm}}
\Statex  function $sol$ = Secant(f, $i_{1}$, $i_{2}$, $n$, $tol$)
\State $error = 1 + tol$ \Comment{To enter the while loop}
\State $r(1) = i_1$ \Comment{Vector to keep track of the root approximations}
\State $r(2) = i_2$
\State $count = 2$ \Comment{Keep count of the iterations}
\While{$count \leq n+1$ and $error > tol$}
    \State $r(count + 1) = r(count) - {f(r(count)) \cdot (r(count) - r(count-1))}/{(f(r(count)) - f(r(count-1)))}$
    \State $error = {\text{norm}(r(count + 1) - r(count))}/{\text{norm}(r(count))}$
    \State $count = count + 1$
\EndWhile
\If{$error < tol$}
    \State $sol = r(\text{length}(r) - 1)$
\Else
    \State $sol = \text{``No Convergence"}$
\EndIf
\end{algorithmic}}
\label{alg:SecantMethod}
\end{algorithm}

\subsubsubsection*{fminunc}\label{subsec:fminunc}

\texttt{fminunc} is MATLAB’s high-level routine for finding a local minimum of a smooth and unconstrained scalar objective.  
Starting from a user supplied initial guess $x_0$, the routine iteratively updates the decision variables until a first order optimality criterion is met and then returns the final estimate $x^\ast$, with $x^\ast$ = \texttt{fminunc}(f,\;$x_0$).

\noindent In the call above, \verb|f| is a function handle that evaluates the objective and (optionally) its gradient, while $x_0$ is the initial point.  The decision vector can be one- or multi-dimensional.  Internally the solver relies on gradient-based quasi-Newton updates (BFGS) by default, but these details are transparent to the user. One need only set the desired tolerances through an \verb|optimset| structure.

\begin{algorithm}[H]
\caption{MATLAB's \texttt{fminunc}}\label{alg:fminuncwrapper}
\begin{algorithmic}[1]
\Statex \textbf{Input:} $f$: function handle; $x_0$: initial guess; $n_{\max}$: max number of iterations; $tol$: error tolerance
\Statex \textbf{Output:} sol: approximation of the minimizer
\Statex
\Statex function $sol$ = \texttt{myfminunc}($x_0$, $n_{max}, tol$)
\State {opts} = {optimset}(`{Display',`off',`MaxIterations'}, $n_{\max}$,`TolX', $tol$)
\State {[sol, fval, exitflag, output]} = {fminunc}($f$, $x_0$, {opts})
\If{{exitflag} $\leq 0$}
    \State sol = NaN
\EndIf
\end{algorithmic}
\end{algorithm}

\subsubsubsection*{fmincon}\label{subsec:fmincon}
Whereas \texttt{fminunc} minimises smooth objectives without restrictions, \texttt{fmincon} targets the more common situation in which the design variables must satisfy explicit constraints. It searches for a local minimiser of a continuously differentiable scalar objective
\[
    \min_{x} f(x)
\]
subject to bound limits $\ell b \leq x \leq u b$, linear inequalities $A x \leq b$, linear equalities $A_{\mathrm{eq}} x = b_{\mathrm{eq}}$, and user-defined nonlinear conditions $g(x) \leq 0,\; h(x) = 0$. A typical engineering example is optimizing the cooling schedule of a solidification process: the objective might be the total solidification time, while constraints enforce energy balance, phase-fraction limits, and temperature bounds. All such restrictions can be encoded directly in \texttt{fmincon}’s call, leaving the numerical details to the solver. 

Internally, \texttt{fmincon} proceeds by solving a sequence of approximate subproblems using methods such as interior-point or SQP. It handles constraints through penalty terms and updates both the decision variables and associated multipliers until some first-order optimality conditions are met.

\begin{algorithm}[H]
\caption{MATLAB's \texttt{fmincon}}\label{alg:fminconwrapper}
\begin{algorithmic}[1]
\Statex \textbf{Input:} $f$: objective handle; $x_0$: initial guess; $A,\,b,\,A_{\mathrm{eq}},\,b_{\mathrm{eq}},\,\ell b,\,u b$: constraints; \textit{nonlcon}: nonlinear constraint handle; $n_{\max}$: max iterations; $tol$: error tolerance
\Statex \textbf{Output:} \texttt{sol}: approximation of the constrained minimiser
\Statex
\Statex function $sol$= \texttt{myfmincon} ($f, x_0, A, b, A_{eq}, b_{eq}, lb,ub, n_{max}, tol$) 
\State {opts} = {optimoptions('fmincon', 'Algorithm', 'interior-point', 'Display', 'off', 'MaxIterations', n\_max, 'TolFun', $tol$)}
\State [sol, fval, exitflag, output] = {fmincon}($f$, $x_0$, $A$, $b$, $A_{\mathrm{eq}}$, $b_{\mathrm{eq}}$, $\ell b$, $u b$, {nonlcon}, {opts})
\If{{exitflag} $\le 0$}
    \State {sol} = \texttt{NaN}  \Comment{no convergence}
\EndIf
\end{algorithmic}
\label{alg fmincon}
\end{algorithm}

\subsubsubsection*{Steepest Descent with Armijo Rule}\label{subsec:steepest-descent}
The Steepest–Descent method is a first-order iterative technique that moves in the direction of the negative gradient of a smooth scalar objective $f\colon\mathbb R^{n}\!\to\!\mathbb R$.  Given a current iterate $x_k$, the update takes the form
\[
    x_{k+1}=x_k-\alpha_k\nabla f(x_k),
\]
where the step size $\alpha_k>0$ must be chosen carefully.  A first-order Taylor expansion gives
\[
    f(x_k-\alpha\nabla f(x_k)) \;\approx\; f(x_k)-\alpha\,\|\nabla f(x_k)\|^{2}+O(\alpha^{2}),
\]
so one expects a decrease proportional to $\alpha\|\nabla f(x_k)\|^{2}$ when $\alpha$ is small.  The Armijo (sufficient-decrease) condition
\[
    f\!\bigl(x_k-\alpha_k\nabla f(x_k)\bigr)\;\le\;f(x_k)-c\,\alpha_k\|\nabla f(x_k)\|^{2}, \quad 0<c<1,
\]
formalises this expectation by demanding that the actual reduction in $f$ is at least a fraction $c$ of the linear model’s prediction.  Choosing $\alpha_k$ via back-tracking until this inequality holds prevents overshooting, maintains descent, and stabilises convergence (especially on poorly scaled or ill-conditioned problems).

\begin{algorithm}[H]
\caption{Steepest Descent with Armijo Line Search}\label{alg:SDArmijo}
\begin{algorithmic}[1]
\Statex \textbf{Input:} $f$: objective function handle; $x_0$: initial guess; $n_{\max}$: max iterations; ${tol}$: error tolerance
\Statex \textbf{Output:} {sol}: approximation of the minimiser
\Statex
\State $x(1)= x_0$, \quad ${iter}= 1$, \quad ${relerr}= 1+{tol}$ \hfill $\triangleright$ initialise
\State $\nabla f = \text{gradient}(f)$ \hfill $\triangleright$ symbolic or automatic differentiation
\While{${iter}<n_{\max}$ \textbf{and} ${relerr}\ge{tol}$}
      \State $g =  \nabla f\bigl(x({iter})\bigr)$,\quad $\|g\| = \|g\|_2$
      \State $\alpha = \texttt{Armijo\_LS}\bigl(x(\text{iter}), g, f\bigr)$ \hfill $\triangleright$ Alg.\,\ref{alg:ArmijoLS}
      \State $x({iter}+1) =  x({iter})-\alpha\,g$
      \State ${relerr} = \dfrac{\|x({iter}+1)-x({iter})\|_2}{\|x({iter})\|_2}$
      \State ${iter} = {iter}+1$
\EndWhile
\If{$\textit{relerr}<\textit{tol}$}
      \State {sol}$ = x({iter})$
\Else
      \State {sol}$ = $ ``No convergence''
\EndIf
\end{algorithmic}
\end{algorithm}

\begin{algorithm}[H]
\caption{Armijo Line Search}\label{alg:ArmijoLS}
\begin{algorithmic}[1]
\Statex \textbf{Input:} $x$: current iterate;\; $g=\nabla f(x)$: gradient;\; $f$: objective handle
\Statex \textbf{Parameters:} initial step $\alpha=1$;\; reduction factor $\beta=0.5$;\; Armijo constant $c=0.1$
\Statex \textbf{Output:} final step size $\alpha$
\Statex
\Statex\texttt{function } $\alpha$ \texttt{ = Armijo\_LS(x, g, f)}
\State $f_\text{current} = f(x), \quad \texttt{norm\_g} = (\texttt{norm}(g,2))^{2}$
\While{$f\bigl(x-\alpha g\bigr) \;>\; f_\text{current}-c\,\alpha\,\|g\|^{2}$}
      \State $\alpha = \beta\,\alpha$ \hfill $\triangleright$ reduce step
\EndWhile
\State \textbf{return} $\alpha$
\end{algorithmic}
\end{algorithm}

\subsubsubsection{Testing Sample}
We test classical methods on two inverse problems: (i) logistic equation (ii) the porous medium equation.
We use a normalized MSE loss and, train on the first 50\% of the data points, and use the others for testing; noise (when used) is additive white Gaussian at  3\% amplitude. Methods compared: steepest descent, \texttt{fminunc}, \texttt{fmincon}, Secant, and Newton (with analytic/symbollic derivatives).
All procedures and results can be found in ~Appendix \S\ref{sec:app}.

\subsection*{Key Findings:}
The experiments on the Logistic equation show that initialization matters, as local methods like Newton and Secant require starting close enough to the true parameters; from distant guesses they either diverge or converge to different points. In contrast, the non-linear solvers \texttt{fminunc}/\texttt{fmincon} tolerate poorer initial guesses better than Newton/secant, although they can still be trapped by flat valleys; among them, \texttt{fmincon} was the most robust overall, likely due to the use of constraints. Beyond initialization, the topology of the space also matters: the minimization procedure behaves well when considering $r$ only, however things get complicated when introducing $K$. After plotting the loss function with respect to $K$ only, we observe that it is negative exponential (no convex, and with no minimizer), which explains the observed difficulties. Re-parameterizing the space by doing a change of variable $ \tilde{K} = log(K)$ as described in Appendix \ref{sec:rk} improves a bit the convergence of some solvers. 
\subsection{The Porous Medium Equation (PME)}

The Porous Medium Equation (PME) is a nonlinear partial differential equation that describes diffusion-like processes in porous media, with applications in gas flow, heat transfer, and biological systems. It is given by 
\begin{equation}
\begin{cases}
u_t(t,x) - \Delta(|u|^{\beta-1}u) = 0\quad & t> t_0, \; x\in (a,b)\\
u(t,a) = 0 & t> t_0\\
u(t,b) = 0 & t> t_0\\
u(t_0,x) = u_0(x)&x\in (a,b) \\
\end{cases}
\end{equation}
where $u(t_i) \geq 0$ is the scaled density of the considered gas at different times $t_i$, $u_0(x)$, the initial scaled density, and \( \beta \) is the  polytropic exponent. 

\noindent  The general one-dimensional Porous Medium Equation (PME) can be written in the following equivalent form as given in Chapter 21 in \cite{Rostamian}:
\[
\begin{cases}
\partial_t u(t,x) = \partial_x\left( \beta u^{\beta - 1} \partial_x u(t,x) \right), & x \in [-1,1],\; t \in [0,1], \\
u(0, x) = u_0(x), & x \in [-1,1],
\end{cases}
\]
where we assume
\[
\beta \in \mathbb{R}^+, \quad u(x,0) > 0, \quad \text{and that } u \text{ is sufficiently regular to allow differentiation.}\]

\noindent As opposed to the classical heat equation, which allows instantaneous spreading of information, the PME allows a finite speed of propagation and hence the disturbances travel at a finite speed. In fact, the PME is fully treated analytically in Vázquez's monograph~\cite{vazquez2007porous}, which can be considered the principal theoretical book of the PME. For the case of $\beta >1$, we will consider the Barenblatt solution, discussed in the book, as a benchmark for validating our numerical methods. 
For the rest of this work, we will only consider the case of $\beta = 3$ for which Barenblatt’s solution takes on a particularly simple form \cite{Rostamian}:
\begin{equation}
\label{eq:exact_solution}
u(t, x) = \frac{1}{(t + \delta)^{1/4}} \sqrt{ \max( 0, \ 1 - \frac{x^2}{12 \sqrt{t + \delta}}  )},
\end{equation}
In the following sections, different methods of discretization are employed. The purpose of discretization is to bring the continuous partial differential equation to a system of algebraic equations that can be solved step by step. Starting with the well-known and well-studied heat equation ($\beta=1$) in Section~\ref {sec: heat}, we will consider some discretization schemes. 
We will next consider the direct and inverse classical problem of the general PME with $\beta > 1$ in Section ~\ref{sec: general pme} and ~\ref{sec: general pme inv}, respectively.

\subsubsection{The Heat Equation} \label{sec: heat}
In the special case where \( \beta = 1 \), the PME reduces to the heat equation.
The heat equation has been one of the most extensively studied partial differential equations since the early work of Fourier and is a model for parabolic problems. For simple domains and boundary conditions, analytical or series solutions can be obtained through separation of variables and Fourier transforms. However, in more general situations, e.g., irregular geometries or non-constant coefficients, numerical methods are required.
\noindent Both analytical and numerical methods for the heat equation are defined and analyzed in a vast amount of literature \cite{firmansah2019study,butt2009numerical,makhtoumi2018numerical, tveito1998finite,foster1972diffusion,vaidya2021heat}.  We do not repeat the full derivations of these methods in the present work. Instead, we only refer the reader to Appendix~\ref{appendix:schemes}, where full discretizations, matrix form and step-by-step calculations of some classical schemes are provided.

\subsubsection{The Direct Problem of the General PME} \label{sec: general pme}

In this section, we solve the direct problem of the general PME. 

\paragraph*{Numerical Scheme and Residual Formulation}
\label{sec: classical}

To approximate the solution of
\[
\partial_t u(t,x) = \partial_x\!\left( \beta u^{\beta-1}(t,x)\,\partial_x u(t,x) \right),
\]
we use a finite difference method with a backward Euler scheme in time and central differences in space.

The spatial and temporal grids are defined by
\begin{align*}
x_i &= -1 + i \Delta x, \quad i = 0,\dots,N, \quad \Delta x = \frac{2}{N}, \\
t^n &= n \Delta t, \quad n = 0,\dots,T.
\end{align*}

We denote the numerical approximation of \( u(t^{n+1},x_i) \) by
\[
u_i^{n+1} := u(t^{n+1},x_i).
\]

Evaluating the PDE at the grid point \( (t^{n+1},x_i) \), we obtain
\[
\partial_t u(t^{n+1},x_i)
=
\partial_x\!\left( \beta u^{\beta-1}\partial_x u \right)(t^{n+1},x_i).
\]

The time derivative is discretized using a backward difference:
\[
\partial_t u(t^{n+1},x_i) \approx
\frac{u_i^{n+1}-u_i^n}{\Delta t}.
\]

The nonlinear diffusion term is discretized in divergence form using central differences:
\[
\partial_x\!\left( \beta u^{\beta-1}\partial_x u \right)(t^{n+1},x_i)
\approx
\frac{1}{\Delta x}
\left[
\left( \beta u^{\beta-1}\partial_x u \right)_{i+\frac12}^{n+1}
-
\left( \beta u^{\beta-1}\partial_x u \right)_{i-\frac12}^{n+1}
\right],
\]
where the numerical fluxes at half-grid points are approximated by
\[
\left( \beta u^{\beta-1}\partial_x u \right)_{i+\frac12}^{n+1}
\approx
\beta\left(\frac{u_{i+1}^{n+1}+u_i^{n+1}}{2}\right)^{\beta-1}
\frac{u_{i+1}^{n+1}-u_i^{n+1}}{\Delta x},
\]
\[
\left( \beta u^{\beta-1}\partial_x u \right)_{i-\frac12}^{n+1}
\approx
\beta\left(\frac{u_i^{n+1}+u_{i-1}^{n+1}}{2}\right)^{\beta-1}
\frac{u_i^{n+1}-u_{i-1}^{n+1}}{\Delta x}.
\]

\noindent Combining the above approximations yields, for each interior grid point \( i=1,\dots,N-1 \), the nonlinear residual equation
\[
F_i(u^{n+1}) =
u_i^{n+1} - u_i^n
-
\frac{\beta \Delta t}{\Delta x^2}
\left[
\left(\frac{u_{i+1}^{n+1}+u_i^{n+1}}{2}\right)^{\beta-1}
( u_{i+1}^{n+1}-u_i^{n+1} )
-
\left(\frac{u_i^{n+1}+u_{i-1}^{n+1}}{2}\right)^{\beta-1}
( u_i^{n+1}-u_{i-1}^{n+1} )
\right].
\]

We define the vector of unknowns at time step \( n+1 \) as
\[
u^{n+1} = (u_1^{n+1},\dots,u_{N-1}^{n+1})^\top \in \mathbb{R}^{m},
\quad m = N-1,
\]
and the residual vector
\[
F(u^{n+1}) = (F_1,\dots,F_{N-1})^\top.
\]
The Dirichlet boundary conditions \( u_0^{n+1}=u_N^{n+1}=0 \) are imposed explicitly.

\paragraph*{Newton's Solver}

The scheme was implemented in Python. 

\noindent At each time step, the nonlinear system is solved using Newton’s method discussed in~\ref{section: newtons}, but in our multidimensional case, the update rule becomes:
\[
J (u^{(k)}) \, \delta u^{(k)} = -F(u^{(k)}), \quad 
u^{(k+1)} = u^{(k)} + \delta u^{(k)},
\]
where \( u^{(k)} \in \mathbb{R}^m \) denotes the vector of discrete solution values at the \(k\)-th Newton iteration and is initialized as \( u^{(0)} = u^n \), the converged solution from the previous time step; for the first time step, \( u^n \) is taken as the initial condition of the PDE. The Jacobian matrix \(J(u^{(k)})\) is defined by
\[
\bigl[J\bigl(u^{(k)}\bigr)\bigr]_{ij}
= \frac{\partial F_i}{\partial u_j}\Bigg|_{u=u^{(k)}},
\qquad i,j=1,\dots,m.
\]
It is approximated numerically using finite differences. Specifically, each column of \( J \) is computed by perturbing a single component of \( u \) with a small value \( h = 10^{-6} \), and evaluating the resulting change in \( F \). The iterations stop when the maximum absolute value of the residual components, \( \max_i |F_i| \), falls below a tolerance of \( 10^{-6} \), or when the maximum number of iterations is reached. After computing the solution over all time steps, we visualize the results using a 2D heatmap of \( u(x,t) \).

\noindent Below, one can find the pseudocode of the PME direct problem for \( \beta = 3 \), using the Barenblatt solution~\ref{eq:exact_solution} for the initial and boundary conditions.

\begin{algorithm}[H]
\caption{Solve Direct PME using Newton's Method with Barenblatt IC/BC}
\begin{algorithmic}[1]
\Statex Set parameters: $m$, $N$, $T$, $dt$, tolerance, max iterations
\State Define spatial grid $x \in [-1, 1]$ and time grid $t$
\State Set initial condition $u^0(x)$ using Barenblatt formula
\State Set Dirichlet boundary values 
\For{$k = 1$ to total time steps}
    \State Set $u_{\text{old}} \gets u^{k-1}$
    \State Initialize $u^{(0)} \gets u_{\text{old}}$
    \For{Newton iteration = 1 to max\_iter}
        \State Compute residual vector $F(u^{(i)})$
        \If{$\|F\|_\infty < \text{tolerance}$} \State \textbf{break} \EndIf
        \State \%Approximate Jacobian $J$ using finite differences:
        \For{each component $j$ of $u$}
            \State Perturb $u_j$ by $h$
            \State Compute $F(u + h e_j)$
            \State $J[:,j] \gets \frac{F(u + h e_j) - F(u)}{h}$
        \EndFor
        \State Solve $J \, \delta u = -F$
        \State Update $u^{(i+1)} \gets u^{(i)} + \delta u$
    \EndFor
    \State Save $u^k \gets u^{(i)}$ to history
\EndFor
\State Plot heatmap of $u(x, t)$ over space-time grid
\end{algorithmic}
\label{alg 11}
\end{algorithm}

\subsubsubsection{Testing and results}
All parameters of the implementation are listed in Table~\ref{tab:parameters}:

\begin{table}[H]
\centering
\begin{tabular}{|l|l|}
\hline
\textbf{Parameter} & \textbf{Value / Description} \\ \hline
Nonlinearity parameter \(\beta \) & \(3.0\)  \\ \hline
Spatial domain & \(x\in[-1,1]\) with \(N=100\) intervals (\(N+1\) grid points) \\ \hline
Time domain & \(t\in[0,1]\) with \(\Delta t=0.01\) (\(100\) time steps) \\ \hline

Boundary conditions & Dirichlet values from Barenblatt at \(x=-1\) and \(x=1\) \\ \hline
Tolerance (Newton) & \(10^{-6}\) \\ \hline
Maximum Newton iterations & \(20\) per time step \\ \hline
Jacobian approximation & Finite differences with perturbation \(h=10^{-6}\) \\ \hline
\end{tabular}
\caption{Parameters used in the Newton-based PME solver with Barenblatt initial and boundary conditions.}
\label{tab:parameters}
\end{table}

\noindent With this configuration, the relative \(L^{2}\) error over the full space–time grid was \(1.56\times10^{-2}\).

\subsubsection{The Inverse Problem of the General PME} \label{sec: general pme inv}

For the inverse problem of the general PME, we use the MATLAB's fmincon  function discussed in Section ~\ref{sec: mini algo} 

\noindent The inverse problem here consists in recovering the $\beta$ parameter given some measurements of $u(x,t)$. It is worth mentioning that solving the inverse problem involves the direct problem as well as will be shown in the following enumeration. We start by synthesizing the data using the exact Barenblatt solution ~\ref{eq:exact_solution} on the same spatial and temporal grid initialized in the direct problem ~\ref{tab:parameters}. The Barenblatt solution is then used as a reference for optimization.

\noindent We estimate \( \beta \) by an objective function that computes the error between the numerical and the exact one. The optimization uses MATLAB's fmincon and is limited to a realistic range of \( \beta \), typically between 1.1 and 10.\\
The key steps are as follows:

\begin{enumerate}
    \item Generate synthetic data from the Barenblatt solution with a known \( \beta \).
    \item Define an objective function that computes the error between the numerical and exact solutions.
        \[
J(\beta) = \sum_{i,k} \big[ u_{\text{num}}(x_i,t_k;\beta) - u_{\text{exact}}(x_i,t_k) \big]^2,
\]
    \item Use Newton’s method ~\ref{alg 11}  at each time step to solve the nonlinear system for a given \( m \).
    \item Perform the optimization using MATLAB’s \texttt{fmincon} ~\ref{alg fmincon}, which iteratively updates \( \beta \) to minimize the error.
\end{enumerate}
\subsubsubsection*{Testing and results} \label{sec: test3}
\noindent The algorithm was tested using the same parameters in Table \ref{tab:parameters} for the Newton's method.  \\  
The synthetic “true” data were generated from the Barenblatt self–similar solution using the exponent
\(\beta=3.0\).  
The initial condition and the Dirichlet boundary values for all times were taken directly from this Barenblatt profile as well.  
Starting from an initial guess \(\beta_{0}=2.0\) for the exponent, the model estimated \( \beta_{\text{estimated}} = 3.0980 \), resulting in an L2 relative error of \( 3.267 \times 10^{-2} \).
\noindent  
The MATLAB code implementing this procedure can be applied to any exponent \(\beta\); we present here the case \(\beta=3.0\) specifically because it will later be compared against the PINN that was trained for \(\beta=3.0\). Other tests using different methods, initializations,  and for $\beta = 2$ can be found in ~\ref{sec:appPME}

\section{Physics-Informed Neural Networks}
Classical numerical methods of differential equations, such as finite difference, finite element, or spectral methods, are usually limited by mesh resolution requirements, discretization error, stability constraints, and inefficient computations in large or complex domains. This is in comparison to novel advancements based on machine learning that have presented new possibilities. Among them, Physics-Informed Neural Networks (PINNs) present a promising, mesh-free alternative for solving forward and inverse problems of differential equations. It directly integrates physical governing equations into the neural network training mechanism and allows us to approximate solutions without the need to employ traditional discretization techniques.

\subsection{Preliminaries: Neural Networks}
At the core of any PINN is a regular feedforward neural network (NN). Conceptually, a neural network is a universal function approximator \cite{hornik1991universal}: it learns a mapping \( f: \mathbb{R}^n \rightarrow \mathbb{R}^m \) from input to output by adjusting internal parameters, namely weights and biases. The network is composed of:
\begin{itemize}
    \item An input layer that receives data (e.g., spatial or temporal coordinates),
    \item One or more hidden layers where transformations and non-linear activations occur,
    \item An output layer that returns the predicted quantity of interest.
\end{itemize}

\begin{figure}[H]
    \centering
    \includegraphics[width=0.8\textwidth]{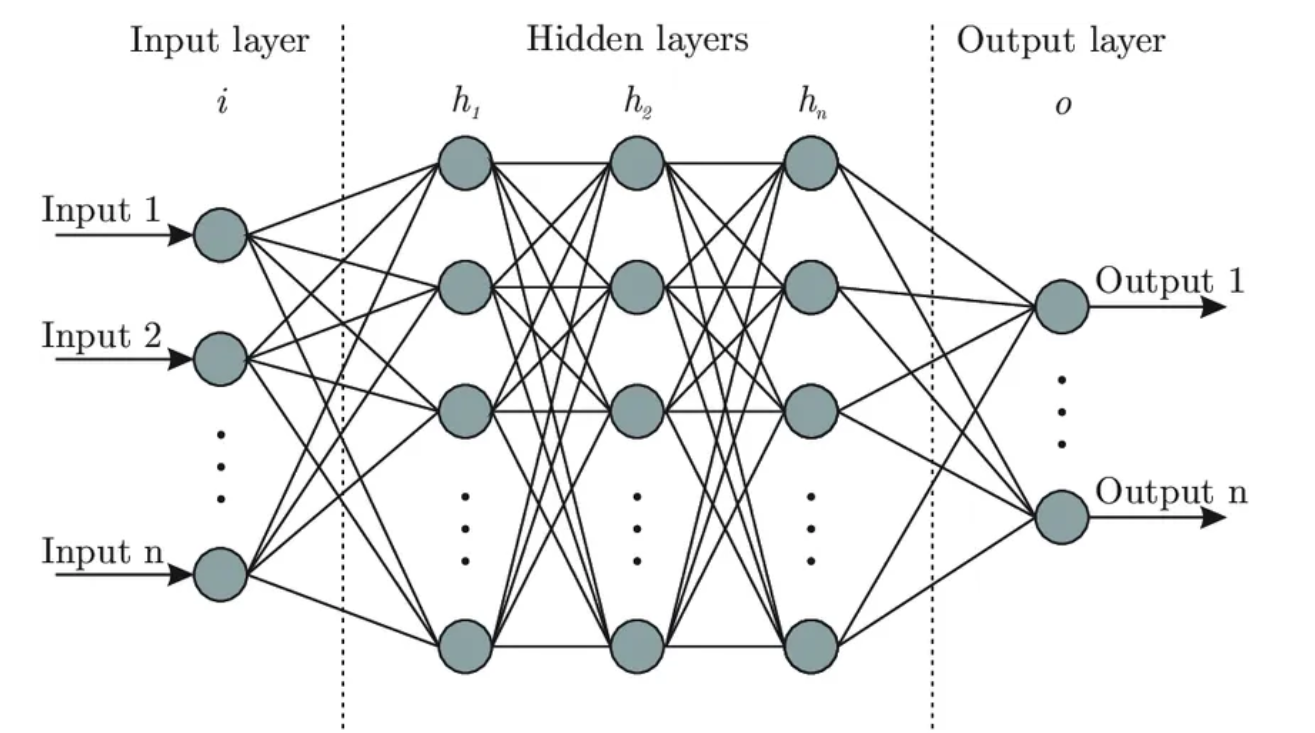}
    \caption{Fully Connected NN example}
\end{figure}

Mathematically, a single layer transforms input \( \mathbf{x} \in \mathbb{R}^n \) to output \( \mathbf{y} \in \mathbb{R}^m \) via:
\[
\mathbf{y} = \sigma(W\mathbf{x} + \mathbf{b}),
\]
where \( W \) is the weight matrix, \( \mathbf{b} \) is the bias vector, and \( \sigma \) is a nonlinear activation function (e.g., tanh, ReLU, or sigmoid). Additionally, training a neural network involves minimizing a loss function, therefore, there exist several Optimization algorithms, such as Adam (a stochastic gradient-based optimizer) and L-BFGS (a quasi-Newton second-order method). These are both commonly used to update weights and biases during training. Adam is efficient for large, noisy datasets, while L-BFGS often converges faster and more accurately in small, smooth problems like PDEs. According to \cite{raissi2019physics} \cite{liu1989limited} \cite{kingma2014adam}, we summarize the difference between them in the following table

\begin{table}[H]
\centering
\caption{Comparison between Adam and L-BFGS optimizers}
\begin{tabular}{|l|c|c|}
\hline
\textbf{Feature} & \textbf{Adam} & \textbf{L-BFGS} \\
\hline
Type & First-order optimizer & Second-order (quasi-Newton) optimizer \\
\hline
Speed per iteration & Fast & Slower (due to line search and history) \\
\hline
Convergence behavior & Good at global search & Good at local refinement \\
\hline
Memory usage & Low & Higher (stores gradient history) \\
\hline
Common use in PINNs & Early training phase & Fine-tuning after Adam \\
\hline
\end{tabular}
\label{tab:optimizers}
\end{table}

\subsection{From Neural Networks to Physics-Informed Neural Networks}
What distinguishes PINNs from traditional neural networks is the fact that they include physical laws, expressed as differential equations, directly in the loss function. Consequently, instead of training only on labeled data, PINNs use the residuals of the governing PDE as part of the loss. This means the network is penalized when it mispredicts known data and also when its predictions violate the underlying physics.

\noindent For instance, suppose we have a PDE of the form:
\[
\mathcal{N}[u](\mathbf{x}, t) = 0,
\]
where \( \mathcal{N} \) is a differential operator acting on the unknown function \( u(\mathbf{x}, t) \). A PINN seeks to approximate \( u \) with a neural network \( u_\theta(\mathbf{x}, t) \), parameterized by \( \theta \). During training, the loss function is composed of multiple terms:
\begin{itemize}
    \item A data loss, ensuring the network fits any available initial, boundary, or measurement data,
    \item A physics loss, ensuring the PDE residual \( \mathcal{N}[u_\theta] \) remains small across the domain.
\end{itemize}

The total loss will take the form:
\[
\mathcal{L}(\theta) = \mathcal{L}_{\text{data}}(\theta) + \mathcal{L}_{\text{physics}}(\theta).
\]

\noindent Thanks to this physics-based guidance, PINNs can perform surprisingly well even when data is scarce. This makes them especially valuable in inverse problems, where direct observations are often limited or noisy.

\subsection{Logistic Equation: The Direct Problem}
To illustrate the implementation and effectiveness of PINNs, we begin with a well-understood problem: the logistic differential equation, which has an analytical solution, making it ideal for validating our PINN methodology.

\subsubsection{Neural Network Architecture}
A fully connected feedforward neural network (FNN) is used to approximate the solution \( p_\theta(t) \), where \( \theta \) denotes the trainable parameters (weights and biases). The network consists of:
\begin{itemize}
    \item An input layer accepting the scalar time variable \( t \),
    \item Two hidden layers with 32 neurons each and hyperbolic tangent (\texttt{Tanh}) activations,
    \item An output layer returning a scalar prediction for \( p(t) \).
\end{itemize}

\subsubsection{Loss Function and Optimization}
The PINN is trained by minimizing a composite loss function:
\[
\mathcal{L} = \mathcal{L}_{\text{ODE}} + \mathcal{L}_{\text{IC}},
\]
where:
\begin{itemize}
    \item \( \mathcal{L}_{\text{ODE}} = \mathbb{E}_t \left[ \left( \frac{dp_\theta}{dt} - r p_\theta(t) \left(1 - \frac{p_\theta(t)}{K} \right) \right)^2 \right] \) for the PDE residual,
    \item \( \mathcal{L}_{\text{IC}} = \left( p_\theta(t_0) - p_0 \right)^2 \) for the initial condition.
\end{itemize}

\noindent To compute derivatives, we use PyTorch's automatic differentiation (\texttt{autograd}).

\subsubsection{Training Procedure}
The model is trained using the Adam optimizer with a learning rate of \(10^{-3}\), for 5000 epochs. The training set consists of 100 collocation points uniformly sampled over the interval \( [t_0, 5] \). The network learns to minimize the residual of the logistic ODE at these collocation points, along with the error on the initial condition.

\subsubsection{Evaluation and Accuracy}
After training, the PINN is evaluated on a denser set of 200 points in the same domain. The relative \( L^2 \) error between the PINN prediction and the analytical solution is computed as:
\[
\text{Relative } L^2 \text{ Error} = \frac{ \| p_\theta - p_{\text{exact}} \|_2 }{ \| p_{\text{exact}} \|_2 }.
\]

\noindent Lastly, the results are visualized by plotting the predicted and exact solutions.

\subsubsection{Testing and Results}
To evaluate our model, we tested it on three different parameter settings (following the same setup described in Table ~\ref{table: rlt errors}).

\subsubsubsection*{Case 1: \( K = 10.0, \quad r = 0.079, \quad p_0 = 20.0 \)}
\begin{figure}[H]
    \centering \includegraphics[width=0.7\textwidth]{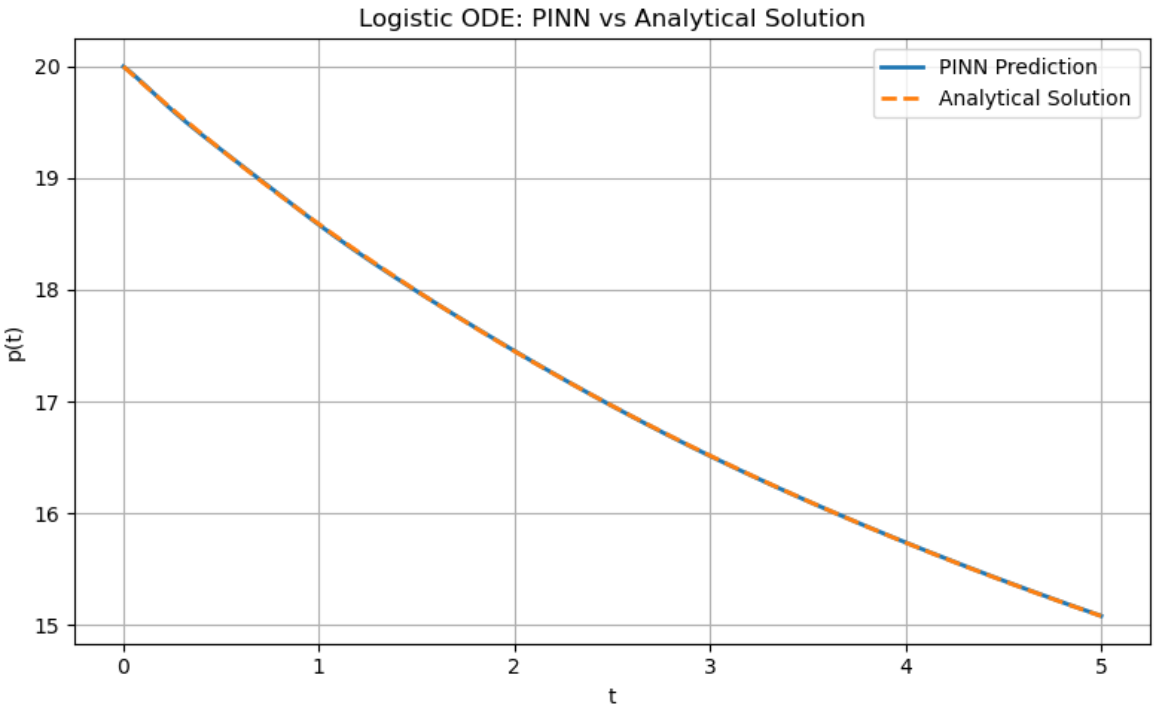} 
    \caption{Direct Problem: PINN Prediction vs Analytical Solution for Case 1.}
    \label{fig:caseX}
\end{figure}

\noindent
Relative \( L^2 \) Error = \text{2.067135 \( \times 10^{-4} \)} 

\subsubsubsection*{Case 2: 
\( K = 90.0, \quad r = 0.05, \quad p_0 = 10.0 \)}

\begin{figure}[H]
    \centering \includegraphics[width=0.75\textwidth]{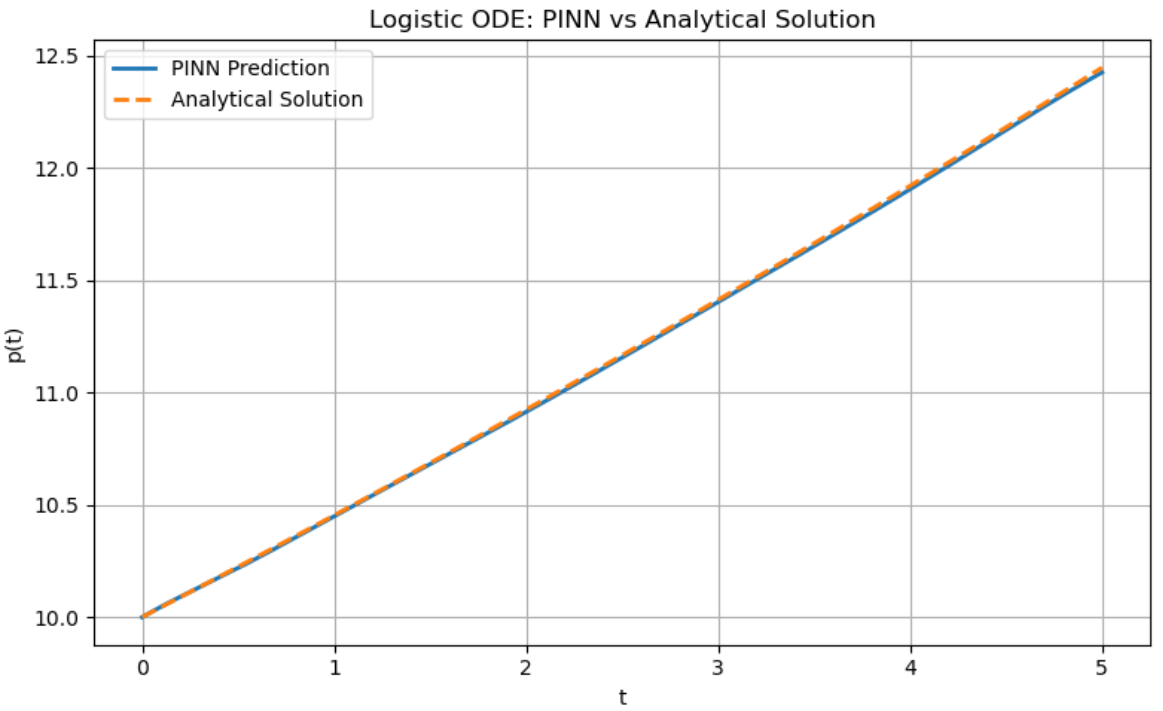} 
    \caption{Direct Problem : PINN Prediction vs Analytical Solution for Case 2.}
    \label{fig:caseX2}
\end{figure}

\noindent
Relative \( L^2 \) Error = \text{9.286355 \( \times 10^{-4}\)}
\subsubsubsection*{Case 3: \( K = 1000.0, \quad r = 0.9, \quad p_0 = 100.0 \)}
\begin{figure}[H]
    \centering \includegraphics[width=0.78\textwidth]{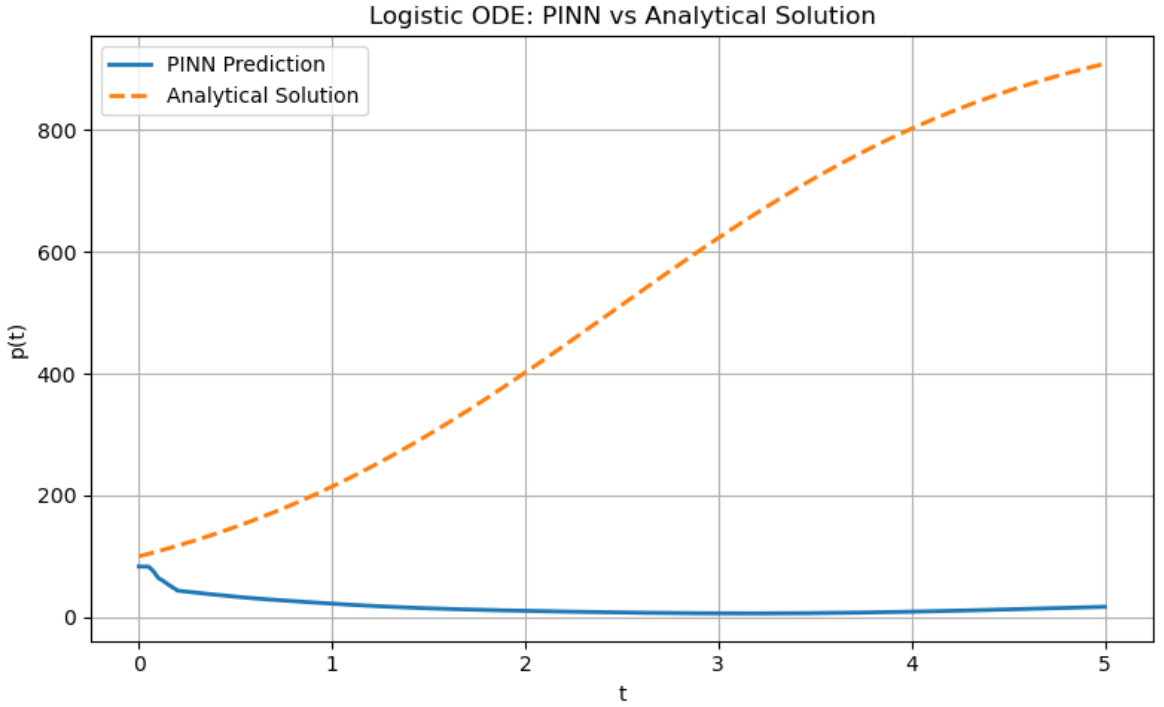} 
    \caption{Direct problem: PINN Prediction vs Analytical Solution for Case 3.}
    \label{fig:caseX3}
\end{figure}

\noindent
Relative \( L^2 \) Error = \text{9.822342 \( \times 10^{-1}\)} 

\subsubsection{Analysis and Interpretation}

\noindent Our main objective is to approximate the direct solution \( p(t) \) over a given time interval using a neural network. We implemented a standard feedforward neural network to approximate \( p(t) \), and we trained it on uniformly sampled time points by minimizing a loss consisting of: the mean squared residual of the (1) differential equation, and (2) the initial condition.

\noindent This initial approach performed well for moderate parameter values. For example, in test cases 1 and 2, the network converged successfully and produced predictions that closely matched the analytical solution, with relative \( L^2 \) errors on the order of a few percent or less.

\noindent However, when testing with significantly larger parameters, specifically case 3, the PINN failed to learn the solution. The training loss plateaued early, and the relative error exceeded 98\% . 

\noindent Upon analysis, two main sources of difficulty were identified: 
\begin{enumerate}
    \item The large value of \( K \) and high growth rate \( r \) produced a steep initial solution profile. This results in a stiffness-like behavior that is challenging for standard neural networks to approximate. 
    \item The scale mismatch between the network output and the true solution range introduced instability in the optimization process. Neural networks with $\tanh$  activations are not naturally equipped to model functions that vary sharply over large output ranges.
\end{enumerate} 

\noindent To address this, we reformulated the problem using a normalized variable $
    u(t) = \dfrac{p(t)}{K}$.\\
\noindent  The normalized equation becomes:

\begin{equation}
    u'(t) = r u(t) (1 - u(t)), \quad u(t_0) = \frac{p_0}{K}
\end{equation}

\noindent Then, we trained the network to predict \( u(t) \) instead of \( p(t) \), and then recovered the original variable using the relation \( p(t) = K \cdot u(t) \). To further support this, we used a sigmoid activation function in the output layer, which ensures that the predictions stayed within the physically meaningful interval \( (0, 1) \).

\noindent This reformulation resulted in a significant improvement: the training became stable and the final relative error \( L^2 \) decreased to \text{2.464935 \( \times 10^{-4} \)}

\begin{figure}[H]
    \centering \includegraphics[width=0.7\textwidth]{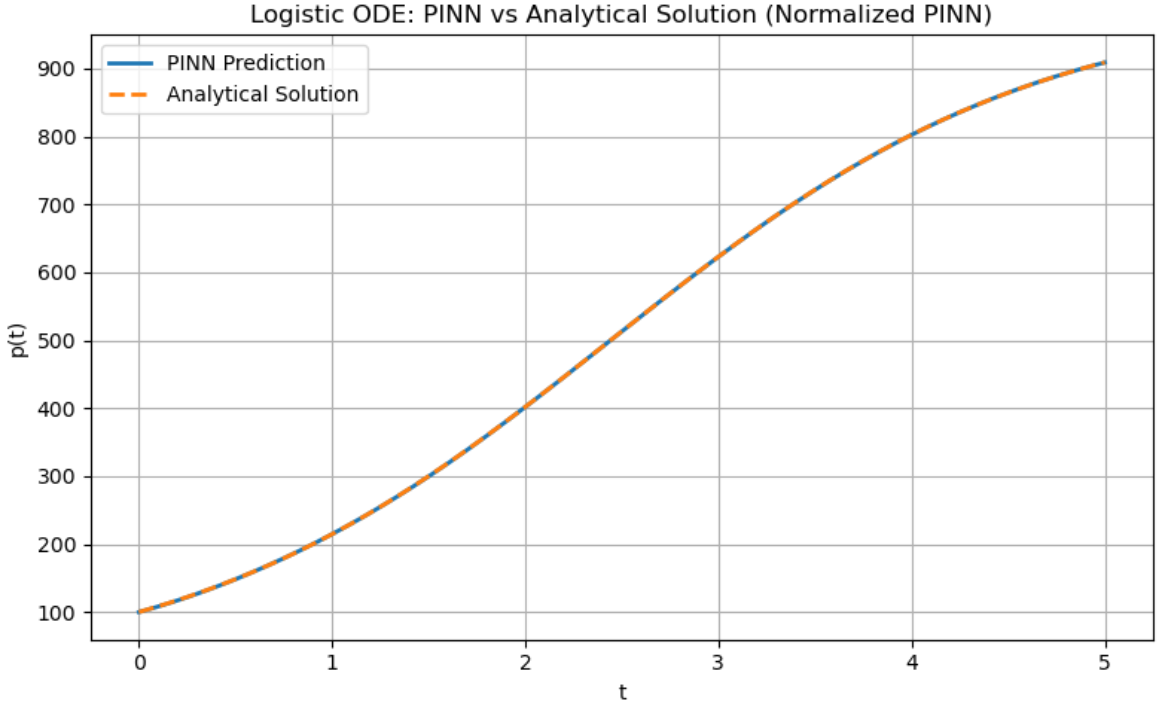} 
    \caption{PINN: direct problem for Case 3 after normalization.}
    \label{fig:caseX4}
\end{figure}

\subsubsection{Classical VS PINN}

\begin{table}[h!]
    \centering
    \begin{tabular}{|c|c|c|c|}
        \hline
        \textbf{Case} & \textbf{PINN Error} & \textbf{ODE45 Error} & \textbf{RK4 Error} \\
        \hline
        Case 1 & \(2.07 \times 10^{-4}\) & \(9.14 \times 10^{-8}\) & \(4.16 \times 10^{-3}\) \\
        Case 2 & \(9.29 \times 10^{-4}\) & \(5.70 \times 10^{-8}\) & \(9.88 \times 10^{-3}\) \\
        Case 3 & \(9.82 \times 10^{-4}\) & \(3.54 \times 10^{-4}\) & \(5.41 \times 10^{-5}\) \\
        \hline
    \end{tabular}
    \caption{Comparison of errors between PINN, ODE45, and RK4 methods}
    \label{tab:error_comparison}
\end{table}

\noindent While all methods give accurate results, PINN performs better than RK4 in most cases but worse than ODE45. ODE45 consistently provides the highest accuracy across most cases.
\subsection{Logistic Equation: The Inverse Problem Recovering One Parameter}

\noindent In the inverse problem, we aim to uncover unknown quantities based on limited observations of the solution. In our case, this means estimating the growth rate \( r \), and later both \( r \) and the carrying capacity \( K \), using only a small set of population measurements over time. One of the most appealing aspects of the PINN approach is how naturally it extends to this setting: the same network architecture used for the direct problem can be adapted for parameter discovery with only minimal modifications.

\subsubsection{Methodology and Differences from the Direct Problem}

The overall construction of the neural network and training procedure is almost identical to the direct problem: we use a fully connected feedforward neural network with two hidden layers and $\tanh$ activations, train on uniformly sampled collocation points, and compute gradients using PyTorch's automatic differentiation.

\noindent However, the inverse formulation introduces two key modifications:
\begin{enumerate}
    \item The unknown parameter \( r \) is treated as a learnable variable. we initialize it with a guess, and it updates during training using gradient descent.
    \item we add an additional data loss term in the loss function to penalize deviations between the network prediction \( p_\theta(t) \) and the observed data points.
\end{enumerate}

\noindent The total loss function becomes:\vspace{-7.5mm}
\begin{eqnarray}
\mathcal{L} = \mathcal{L}_{\text{ODE}} + \mathcal{L}_{\text{IC}} + \lambda_{\text{data}} \cdot \mathcal{L}_{\text{data}},
\end{eqnarray}
where:
\begin{itemize}
    \item \( \mathcal{L}_{\text{ODE}} = \mathbb{E}_t \left[ \left( \dfrac{dp_\theta}{dt} - r \cdot p_\theta(t) \left(1 - \dfrac{p_\theta(t)}{K} \right) \right)^2 \right] \),\vspace{-2mm}
    \item \( \mathcal{L}_{\text{IC}} = \left( p_\theta(t_0) - p_0 \right)^2 \),
    \item \( \mathcal{L}_{\text{data}} = \mathbb{E}_{t_i} \left[ \left( p_\theta(t_i) - p_{\text{data}}(t_i) \right)^2 \right] \),
    \item \( \lambda_{\text{data}} \) is a weighting factor that balances the two parts of the loss function. 
\end{itemize}
\subsubsection{Training and Evaluation}

We generate synthetic data using the analytical solution with known parameters \( K  \), \( r  \), and \( p_0  \), and sample 30 points uniformly over the interval \( [0, 10] \). The PINN is trained using the Adam optimizer on a combined set of physics and data constraints for 10,000 epochs. 

\subsubsection{Testing and results}
We tested the PINN inverse problem using the same cases as in the previous section
\subsubsubsection*{Case 1: 
\( K = 10.0,  \quad p_0 = 20.0 \)}\vspace{-2mm}
\begin{figure}[H]
    \centering \includegraphics[width=0.7\textwidth]{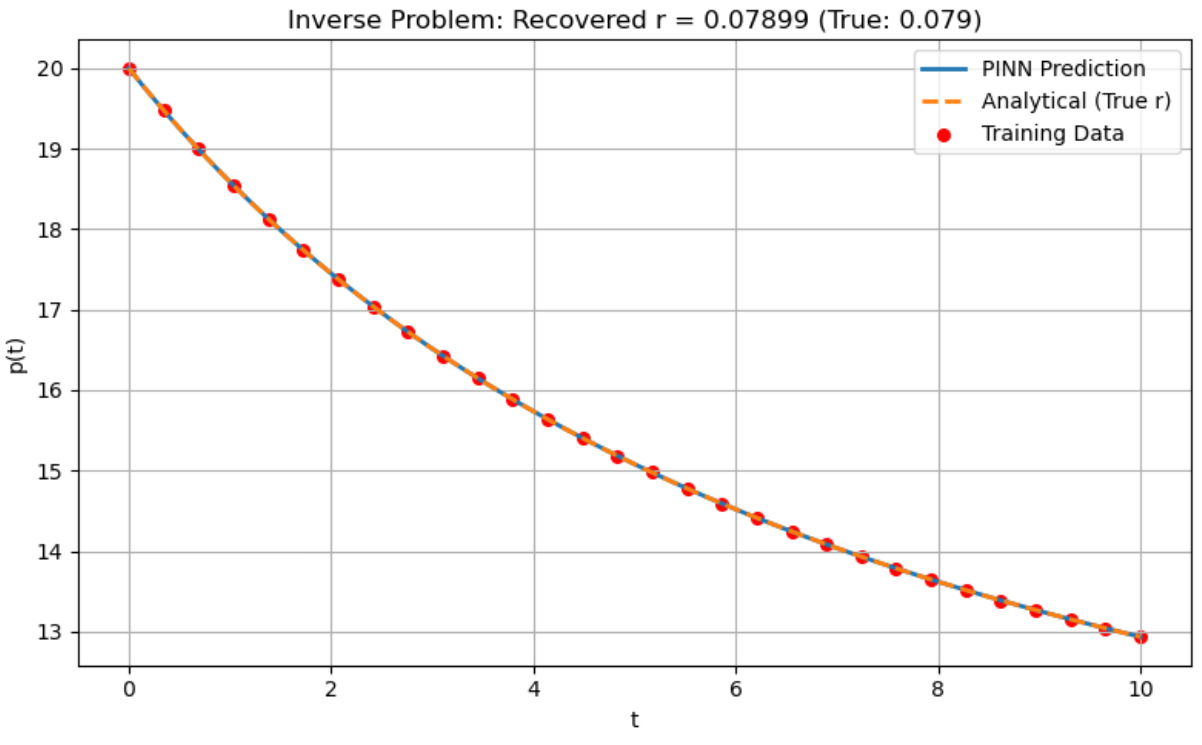} \vspace{-1mm}
    \caption{Inverse problem over r for Case 1.}
    \label{fig:caseX5}
\end{figure}
\noindent
{Recovered r = 0.07899 \quad
True r = 0.079 \quad 
Relative \( L^2 \) Error = \text{5.370742 \( \times 10^{-6 }\)} }
\subsubsubsection*{Case 2: \( K = 90.0, \quad p_0 = 10.0 \)}\vspace{-2mm}
\begin{figure}[H]
    \centering \includegraphics[width=0.7\textwidth]{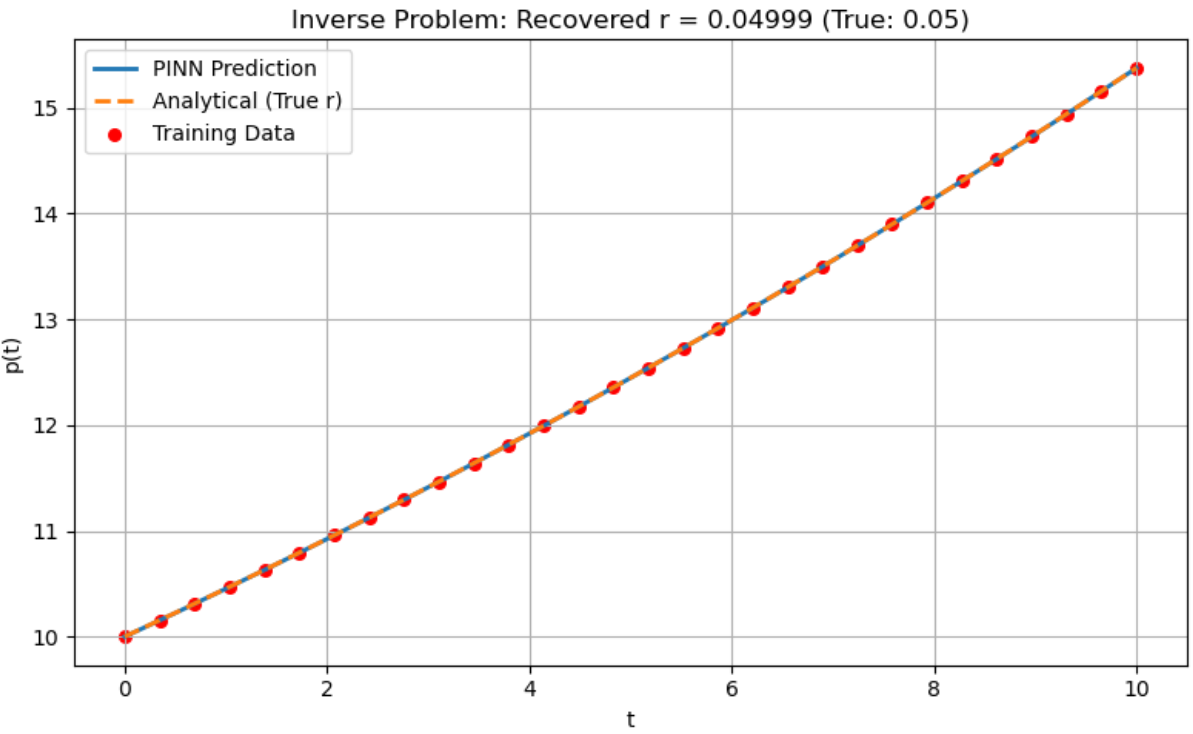} \vspace{-1mm}
    \caption{Inverse Problem over r for Case 2.}
    \label{fig:caseX6}
\end{figure}

\noindent
Recovered r = 0.04999 \quad
True r = 0.05 \quad
Relative \( L^2 \) Error = \text{4.948125e \( \times 10^{-6} \)}

\subsubsubsection*{Case 3: \( K = 1000.0, \quad p_0 = 100.0 \)}\vspace{-2mm}
\begin{figure}[H]
    \centering \includegraphics[width=0.7\textwidth]{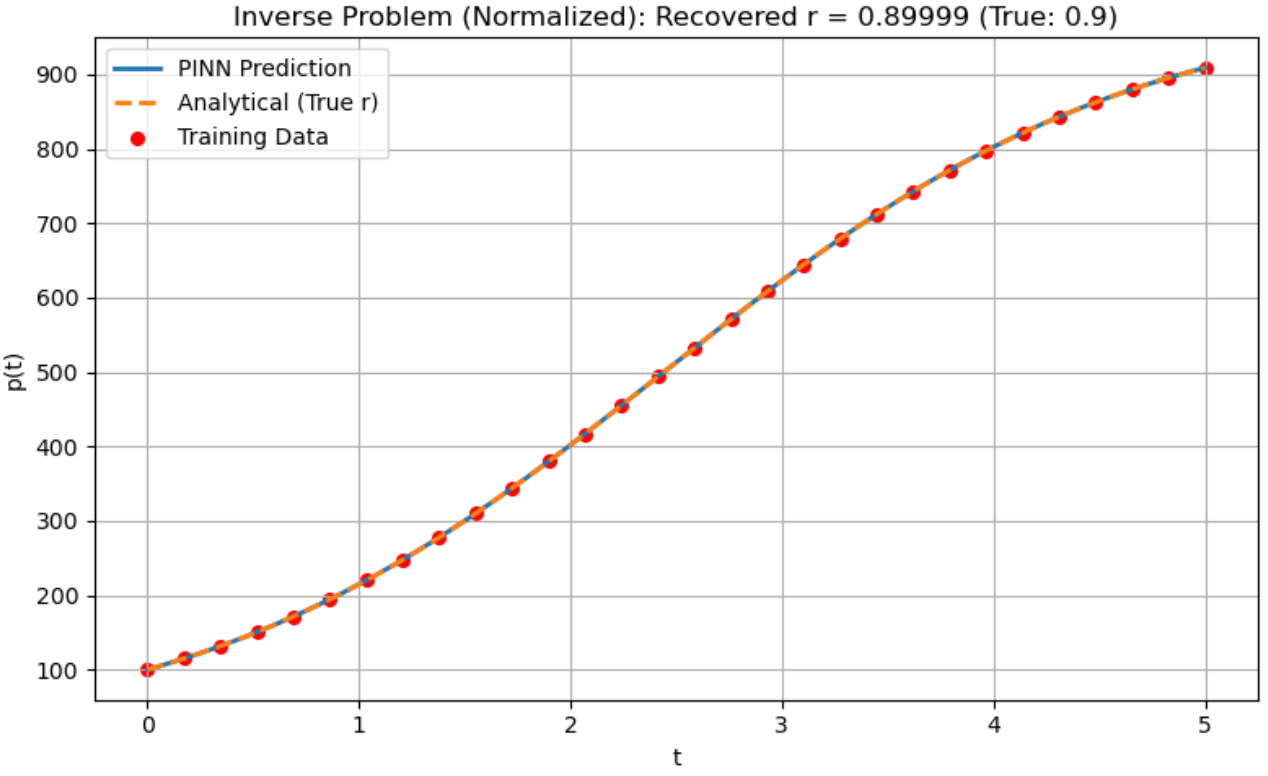} 
    \caption{Inverse Problem over r for Case 3.}
    \label{fig:caseX7}
\end{figure}

\noindent
Recovered r = 0.89999 \quad
True r = 0.9 \quad
Relative \( L^2 \) Error = \text{ 4.256486\( \times 10^{-9} \)}

\subsubsection{Analysis and Interpretation}
The results show that the PINN method works extremely well in recovering the growth rate parameter \( r \) accurately in all three scenarios. The very low relative \( L^2 \) errors show us precisely how well the model works. This tells us that PINNs can effectively combine the underlying differential equation with just observations to estimate the correct parameters, even if the observations available are limited.

\subsection{Logistic Equation: The Inverse Problem Recovering Two Parameters}
We aim to simultaneously infer both the growth rate $r$ and the carrying capacity $K$ in the logistic growth model using PINNs. The approach required minimizing the residual of the ODE together with data and initial condition losses, while treating $r$ and $K$ as trainable parameters.

\subsubsection{First Attempt: Traditional Procedure}
In our initial setup, we designed the PINN to directly approximate the population $p(t)$. The network output was unconstrained, and $r$, $K$ were initialized as raw variables transformed via a softplus function to ensure positivity. The loss function was the same as eq(16)

\noindent Training was performed using the Adam optimizer for 8000 epochs.

\noindent However, this implementation did not succeed in any of the 3 cases we have.

\subsubsection{Second Attempt: Optimizer Improvement and Normalization}
To improve performance, we experimented with several hyperparameters: increasing depth and width of the neural network, changing learning rates, and adjusting loss weights. None of these changes led to meaningful improvements in parameter recovery.

\noindent We then adopted a two-stage training approach, using:

\begin{enumerate}
    \item Adam for initial coarse optimization, followed by
    \item L-BFGS, a second-order optimizer known for fine-tuning convergence in PINNs.
\end{enumerate}
\noindent This change was successful. With this two-stage optimization strategy and good initializations:

\begin{itemize}
    \item Case 1 was successfully solved, accurately recovering $r$ and $K$.
\begin{figure}[H]
    \centering \includegraphics[width=0.7\textwidth]{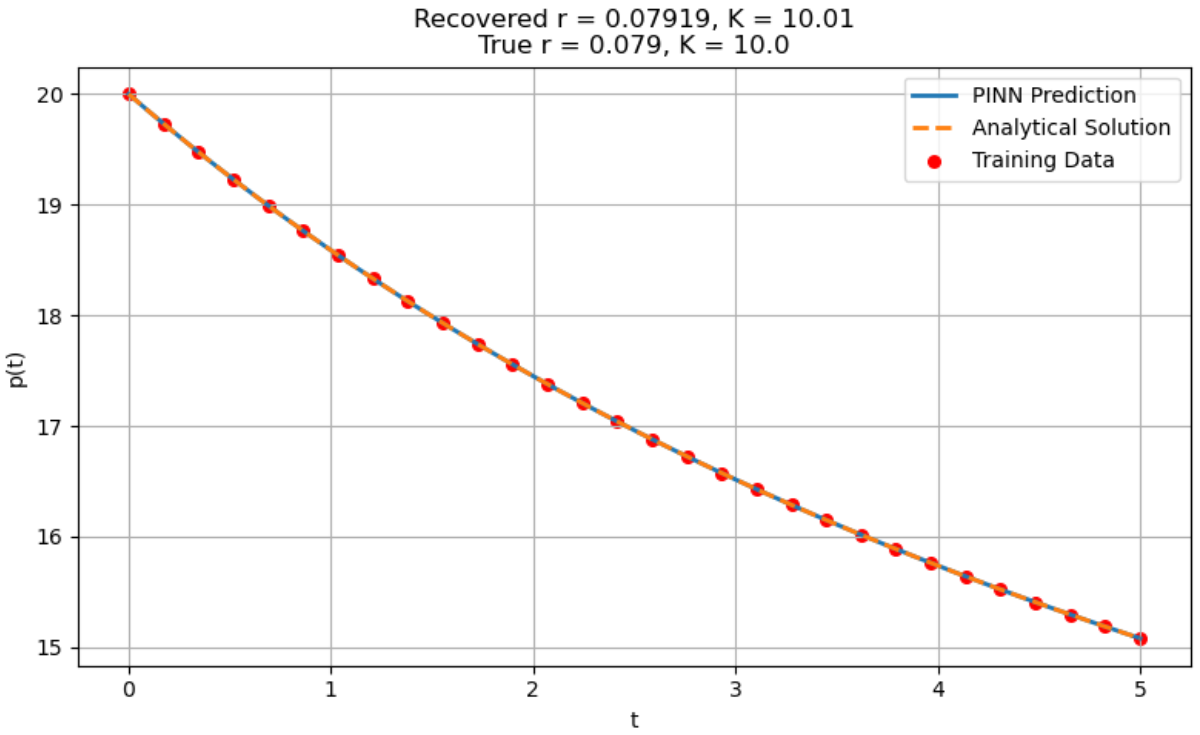} 
    \caption{Inverse problem over r and k for Case 1.}
    \label{fig:caseX8}
\end{figure}
    
    \item Case 2 still showed poor performance
\begin{figure}[H]
    \centering \includegraphics[width=0.7\textwidth]{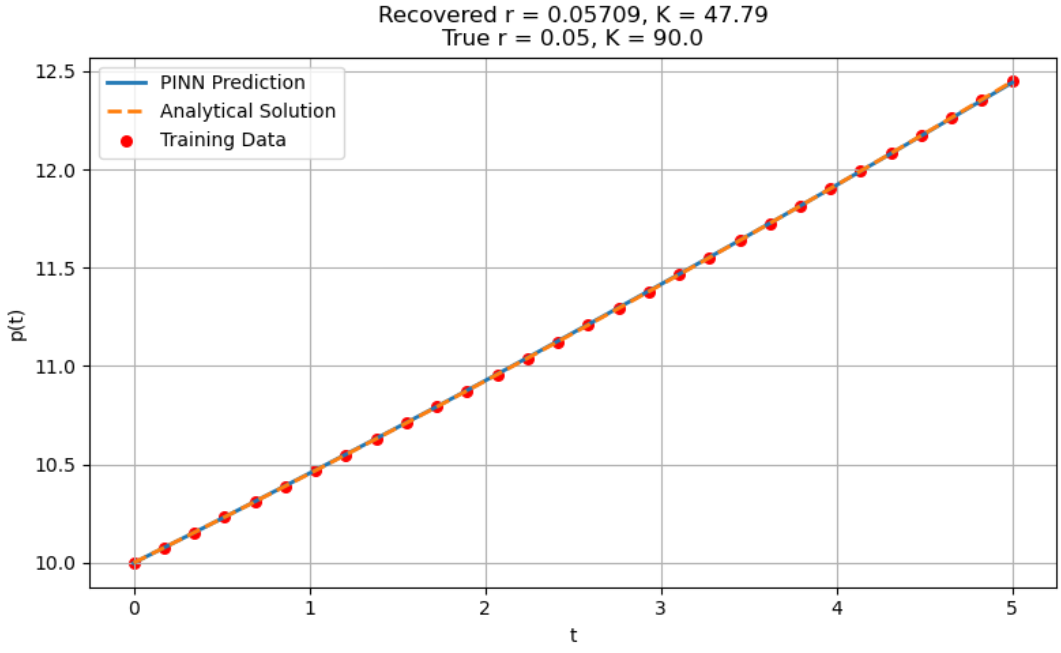} 
    \caption{Inverse problem over r and k for Case 2}
    \label{fig:caseX9}
\end{figure}
    
    \item Case 3, with normalization, showed clear improvements, but still not quite accurate

\begin{figure}[H]
    \centering \includegraphics[width=0.7\textwidth]{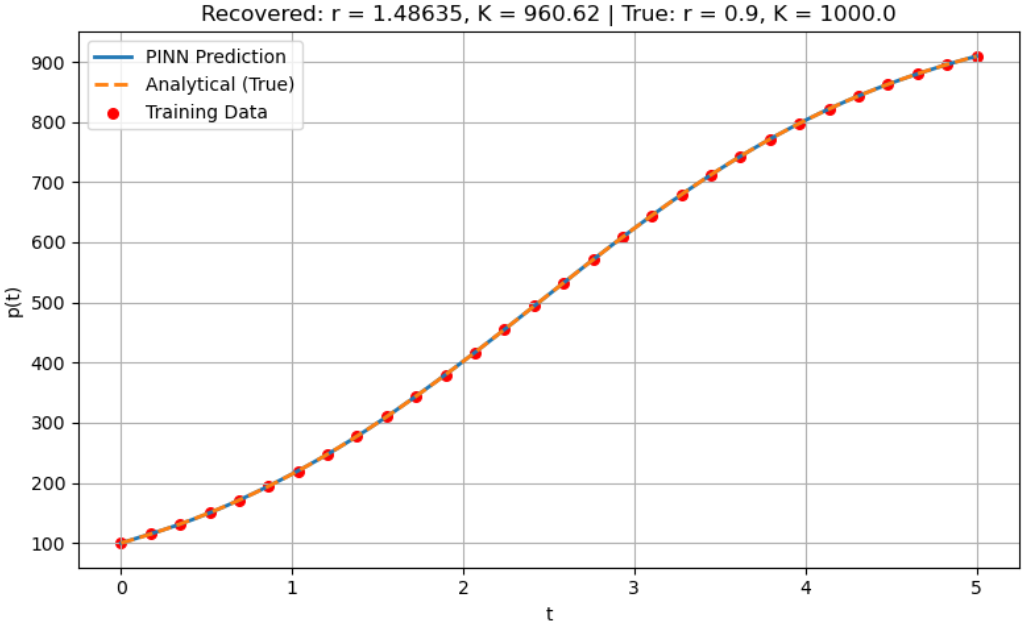} 
    \caption{Inverse problem over r and k for Case 3 after normalization.}
    \label{fig:caseX10}
\end{figure}
    
\end{itemize}

\subsection{PME: Direct Problem for $\beta$ = 3}
We now consider the numerical approximation of solutions to the PME using PINNs. The specific form studied in this section corresponds to $\beta$ = 3 whose form, exact solution, and analysis can be found in Chapter 21 in \cite{Rostamian}:

\begin{equation}
\label{eq:pme}
\frac{\partial u}{\partial t} = \frac{\partial}{\partial x} \left( 3u^2 \frac{\partial u}{\partial x} \right), \quad x \in [-1, 1],\ t \in [0, 1],
\end{equation}

\noindent we consider an exact solution shifted by a small positive parameter \( \delta > 0 \):

\begin{equation}
\label{eq:exact_solution}
u(t, x) = \frac{1}{(t + \delta)^{1/4}} \sqrt{ \max( 0, \ 1 - \frac{x^2}{12 \sqrt{t + \delta}}  )},
\end{equation}

\subsubsection{Neural Network Architecture}

We construct a fully-connected feedforward neural network \( \mathcal{N}_\theta: \mathbb{R}^2 \to \mathbb{R} \), with input \( (t, x) \), output \( u_\theta(t,x) \), and trainable parameters \( \theta \). The architecture is defined as follows:

\begin{itemize}
    \item Input layer: size 2 (for \( t \) and \( x \)).
    \item Hidden layers: 4 layers, each with 20 neurons.
    \item Activation function: hyperbolic tangent \( \tanh \)
    \item Output layer: single neuron returning \( u_\theta(t,x) \).
\end{itemize}

\noindent All weights of the neural network are initialized with the \\textit{Xavier initialization} procedure that is especially tailored for the \( \tanh \) activation function. This process initializes the weights by sampling random values from a uniform distribution over the range:
\[
W \sim \mathcal{U}\left( -\sqrt{\frac{6}{n_{\text{in}} + n_{\text{out}}}},\ \sqrt{\frac{6}{n_{\text{in}} + n_{\text{out}}}} \right),
\]
\noindent where \( n_{\text{in}} \) and \( n_{\text{out}} \) are the number of input and output units of a layer, respectively. Such initialization helps maintain stable gradients during training, especially when activation functions like \(\tanh \) are employed since it maintains the signal variance across layers constant.

\subsubsection{Loss Function and Optimization}
The training objective is to minimize a composite loss \( \mathcal{L} \) that enforces agreement with the initial and boundary conditions as well as the underlying PDE. Let \( \mathcal{D}_t \), \( \mathcal{D}_b \), and \( \mathcal{D}_\Omega \) denote sampled points on the initial boundary \( t=0 \), spatial boundaries \( x = \pm 1 \), and interior domain \( \Omega \), respectively.

\noindent The total loss is:

\begin{equation}
\mathcal{L}_\theta = \log_{10} \left(\underbrace{ \lambda_u \cdot ( \mathcal{L}_b + \mathcal{L}_t ) + \mathcal{L}_{\text{PDE}}}_{\text{Physics Loss}} \right),
\end{equation}

where:

\begin{align}
\mathcal{L}_b &= \mathbb{E}_{(t,x)\in\mathcal{D}_b} \left[ \left( u_\theta(t,x) - u_{\text{BC}}(t,x) \right)^2 \right], \\
\mathcal{L}_t &= \mathbb{E}_{(t,x)\in\mathcal{D}_t} \left[ \left( u_\theta(t,x) - u_0(x) \right)^2 \right], \\
\mathcal{L}_{\text{PDE}} &= \mathbb{E}_{(t,x)\in\mathcal{D}_\Omega} \left[ \left( \partial_t u_\theta - \partial_x \left( 3u_\theta^2 \partial_x u_\theta \right) \right)^2 \right].
\end{align}

Gradients are computed using automatic differentiation.

\subsubsection{Training Procedure}
The neural network is trained using the following optimization approach:

\begin{enumerate}
    \item ADAM Optimizer
    \item followed by L-BFGS in another trial
    \item Early Stopping: The training halts if the validation loss does not improve over a pre-defined patience window.
\end{enumerate}

\noindent The optimizer minimizes the logarithmic loss, and training points are generated using a Sobol sequence, which ensures a quasi-random coverage of the space-time domain.

\subsubsection{Evaluation and Accuracy}

The solution \( u_\theta(t,x) \) is evaluated against the analytical form . We compute the relative \( L^2 \)-error over a dense set of 50{,}000 test points generated by Sobol:

\begin{equation}
\epsilon_{rel} = \frac{ \left\| u_\theta - u \right\|_{L^2} }{ \left\| u \right\|_{L^2} } = \sqrt{ \frac{ \sum_{i=1}^N \left( u_\theta(t_i, x_i) - u(t_i, x_i) \right)^2 }{ \sum_{i=1}^N u(t_i, x_i)^2 } }.
\end{equation}

\subsubsection{Testing and Results}
The network was trained with \( n_{\text{int}} = 256 \) , \( n_{\text{sb}} = 64 \) (for each side), \( n_{\text{tb}} = 64 \), \( L = 4 \) hidden layers, and \( N = 20 \) neurons per layer.

\textbf{Case 1 : Adam Optimizer alone}
\begin{itemize}
    \item 10000 epochs
    \item Final log loss: 0.266907
    \item PDE L2 Relative Error Norm: 4.435748 $\times 10^{-2}$
\end{itemize}

\begin{figure}[H]
\centering
\includegraphics[width=0.85\textwidth]{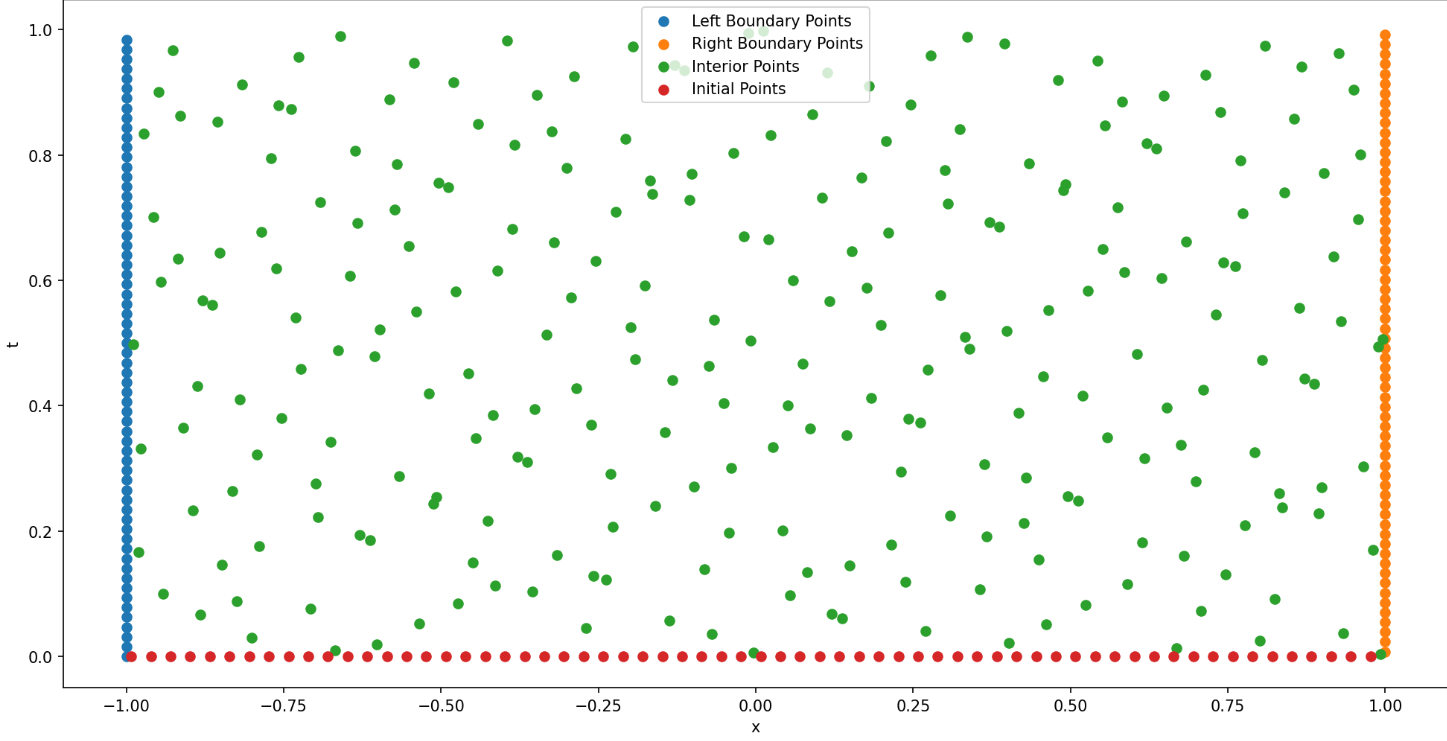}
\caption{Distribution of interior, boundary, and initial points sampled using a Sobol sequence.}
\label{fig:training_points}
\end{figure}

\begin{figure}[H]
\centering
\includegraphics[width=0.6\textwidth]{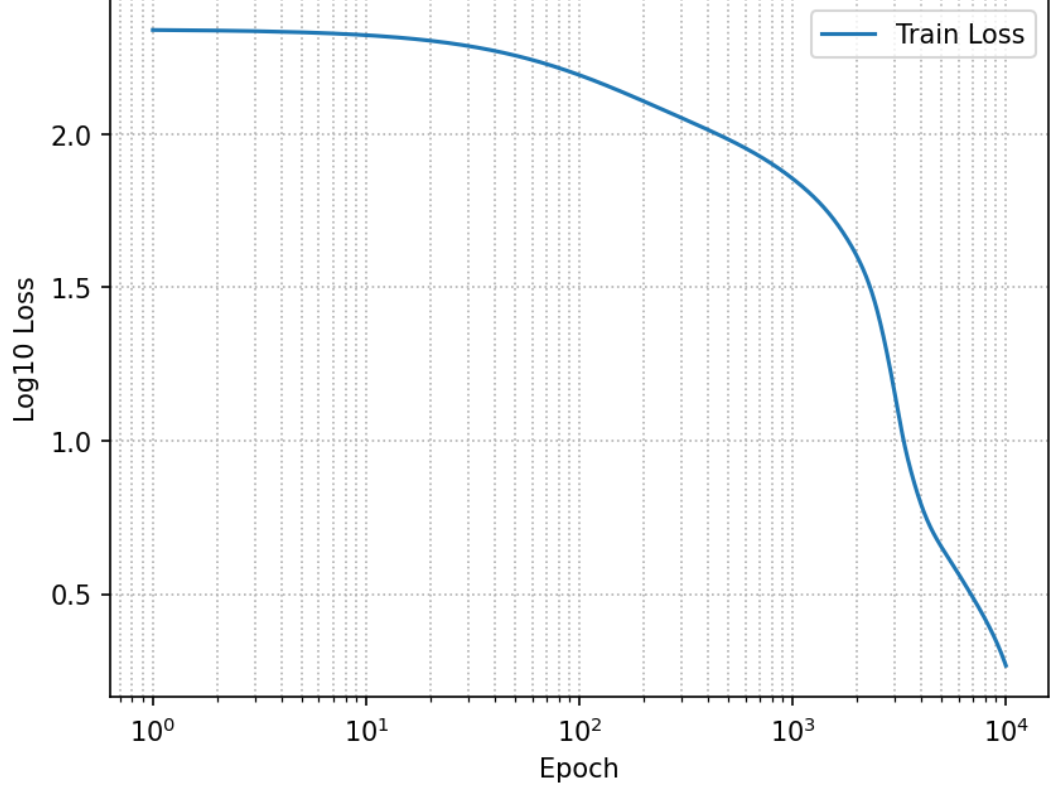}
\caption{Training loss (log scale) over epochs.}
\label{fig:train_loss}
\end{figure}

\begin{figure}[H]
\centering
\includegraphics[width=\textwidth]{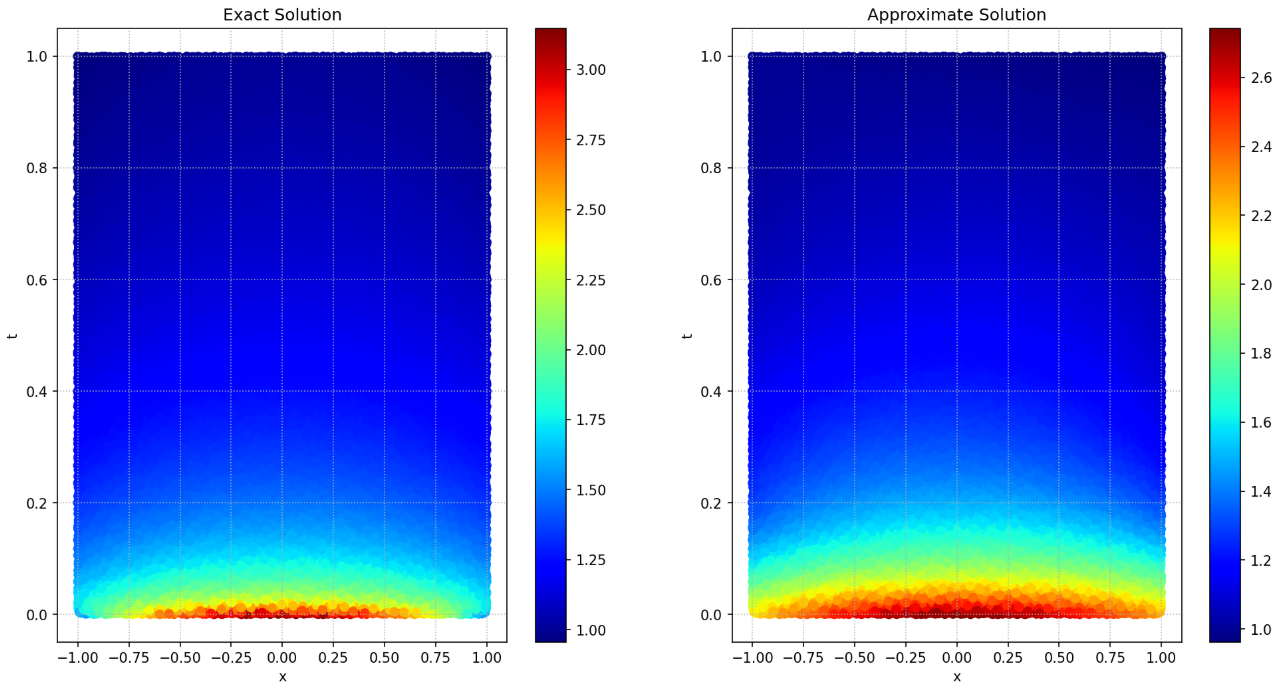}
\caption{Comparison of exact solution (left) and learned PINN solution (right) in the \( (t,x) \) domain.}
\label{fig:solution}
\end{figure}

\textbf{Case 2 : Adam  $+$ L-BFGS optimizers}

\begin{itemize}
    \item 37 epochs
    \item Final loss: -2.331878 
    \item PDE L2 Relative Error Norm: 1.095452 $\times 10^{-3}$
\end{itemize}

\begin{figure}[H]
\centering
\includegraphics[width=0.6\textwidth]{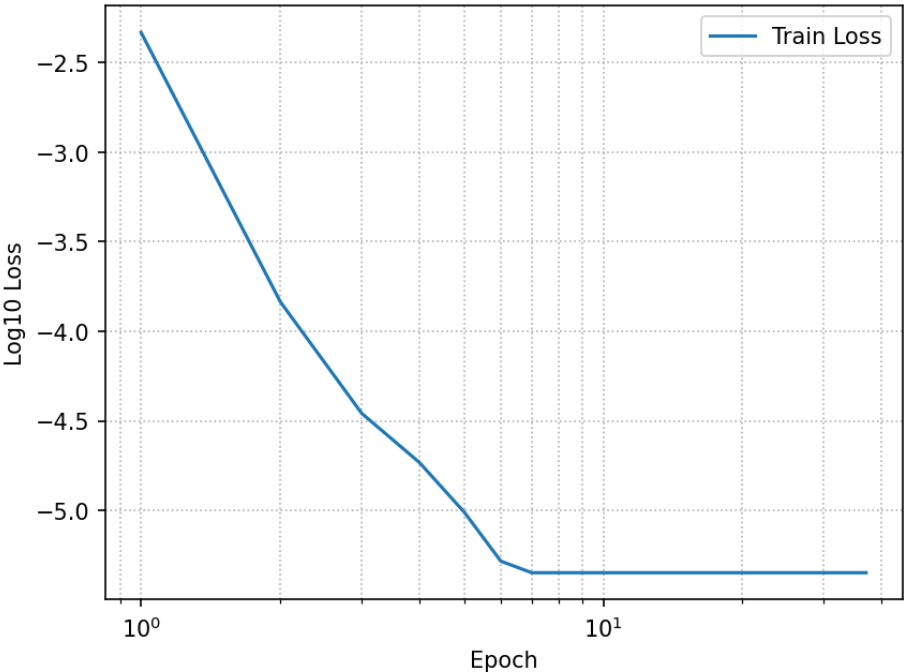}
\caption{Training loss (log scale) over epochs.}
\label{fig:train_loss}
\end{figure}

\begin{figure}[H]
\centering
\includegraphics[width=\textwidth]{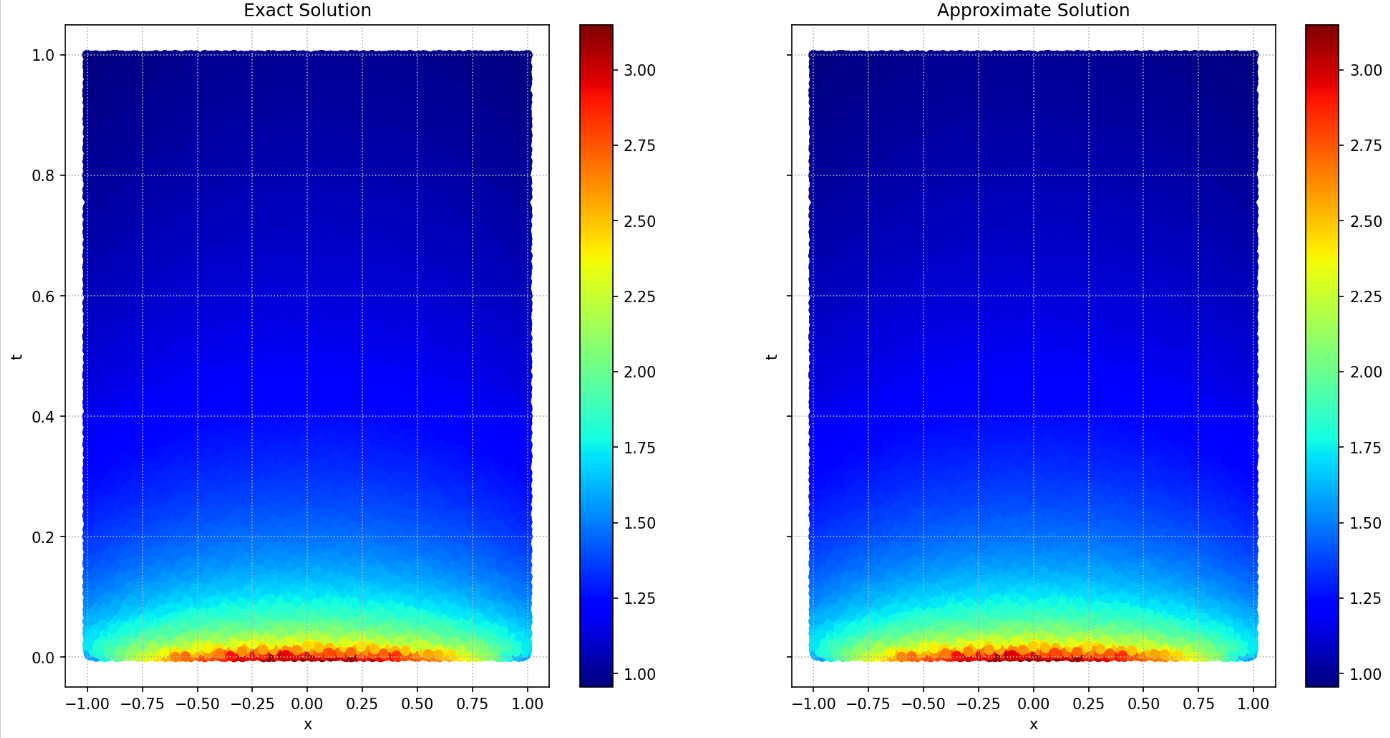}
\caption{Comparison of exact solution (left) and learned PINN solution (right) in the \( (t,x) \) domain.}
\label{fig:solution}
\end{figure}

\subsubsection{Analysis Comparison}
We started out with optimizing the PINN using the Adam optimizer, and gave us a decent start; we achieved a relative L2 error of around 4.43\%, and a final loss (in log scale) of around 0.267. Not bad, but clearly, there was still room for improvement. So we switched things around and brought in L-BFGS right after the Adam phase. That change had a profound effect: the loss fell to approximately -5.35 (log10), and relative L2 error fell to just 0.11\%, an enormous improvement- close to ten times better. Even more surprising, this jump in accuracy happened swiftly (within just 37 steps of training) and early stopping kicked in since the model had converged.

\noindent It is also worth noting that with the Adam optimizer only, the PINN performed worse than the classical method, which had a relative error of about 1.5\%. However, as soon as we combined L-BFGS with Adam, the PINN outperformed the classical solver accuracy.

\subsection{PME: Inverse problem for $\beta$ = 3}
We now consider the inverse formulation of the porous medium equation (PME), where the nonlinearity parameter \( m \) is treated as an unknown to be identified from data. The goal is to recover both the solution \( u(t,x) \) and the parameter \( m \) simultaneously, given a set of measurements \( \{ u_{\text{meas}}(t_i, x_i) \}_{i=1}^{N_{\text{meas}}} \) collected over the spatio-temporal domain.

\noindent The governing PDE is given by:
\begin{equation}
\label{eq:pme_inverse}
\frac{\partial u}{\partial t} = \frac{\partial}{\partial x} \left( \beta u^{\beta - 1} \frac{\partial u}{\partial x} \right), \quad x \in [-1, 1],\ t \in [0, 1],
\end{equation}
where \( \beta > 1 \) is an unknown constant.

\subsubsection{Neural Network Architecture}
We adopt a physics-informed neural network (PINN) to solve the inverse problem. The architecture consists of:

\begin{itemize}
    \item A fully connected feedforward neural network \( \mathcal{N}_\theta : \mathbb{R}^2 \rightarrow \mathbb{R} \), taking input \( (t, x) \) and outputting \( u_\theta(t,x) \).
    \item The parameter \( 
    \beta\in \mathbb{R} \) is treated as a trainable scalar variable, initialized to a value \( \beta_0 \neq 3 \).
    \item The network contains 4 hidden layers, each with 20 neurons and hyperbolic tangent \( \tanh \) activations.
\end{itemize}

\noindent The network weights are initialized using Xavier initialization, as described in the direct problem section.

\subsubsection{Loss Function and Optimization}
The total loss function for training includes physics constraints, initial/boundary data, and measurement data. It is defined as:
\begin{equation}
\mathcal{L}_\theta = \log_{10} \left( 
\lambda_u \cdot (\mathcal{L}_b + \mathcal{L}_t) +
\mathcal{L}_{\text{PDE}} +
\lambda_s \cdot \mathcal{L}_{\text{meas}}
\right),
\end{equation}
where:

\begin{itemize}
    \item \( \mathcal{L}_b \), \( \mathcal{L}_t \), and \( \mathcal{L}_{\text{PDE}} \) are defined identically as in the direct problem and enforce adherence to boundary, initial, and PDE constraints.
    \item \( \mathcal{L}_{\text{meas}} \) is the data mismatch term:
    \[
    \mathcal{L}_{\text{meas}} = \mathbb{E}_{(t,x)\in\mathcal{D}_{\text{meas}}} \left[ \left( u_\theta(t,x) - u_{\text{meas}}(t,x) \right)^2 \right],
    \]
    where \( \mathcal{D}_{\text{meas}} \) denotes the set of measurement points where we evaluate the exact solution. These measurement points are generated using the Barenblatt exact solution.
\end{itemize}

\subsubsection{Training Procedure}
The inverse problem was solved by training a physics-informed neural network to simultaneously learn the solution \( u(t, x) \) and estimate the unknown exponent \( m \) in the porous medium equation. The training used Sobol-generated points with the following setup:
\begin{itemize}
    \item Interior points: \( n_{\text{int}} = 256 \)
    \item Boundary points (left and right): \( n_{\text{sb}} = 64 \) on each side
    \item Temporal boundary (initial condition): \( n_{\text{tb}} = 64 \)
    \item Measurement data points: \( 40 \times 40 = 1600 \)
    \item Optimizer: ADAM for 10{,}000 epochs
\end{itemize}

\begin{figure}[H]
\centering
\includegraphics[width=0.85\textwidth]{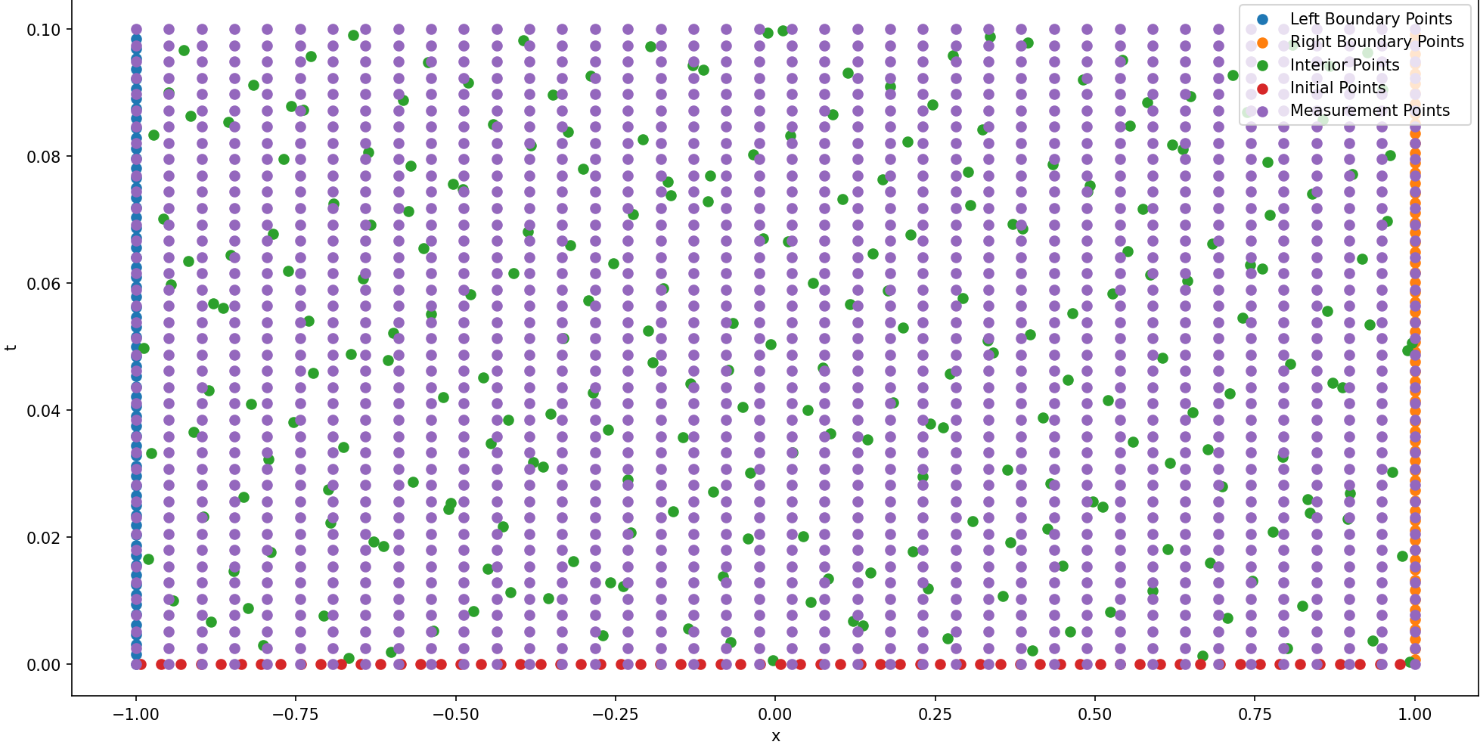}
\caption{Sobol-distributed training and measurement points across the domain.}
\label{fig:inv_training_points}
\end{figure}
For evaluation and plotting, the solution was tested on 50{,}000 Sobol points from the domain.

\subsubsection{Testing and Results}
\begin{figure}[H]
\centering
\includegraphics[width=0.65\textwidth]{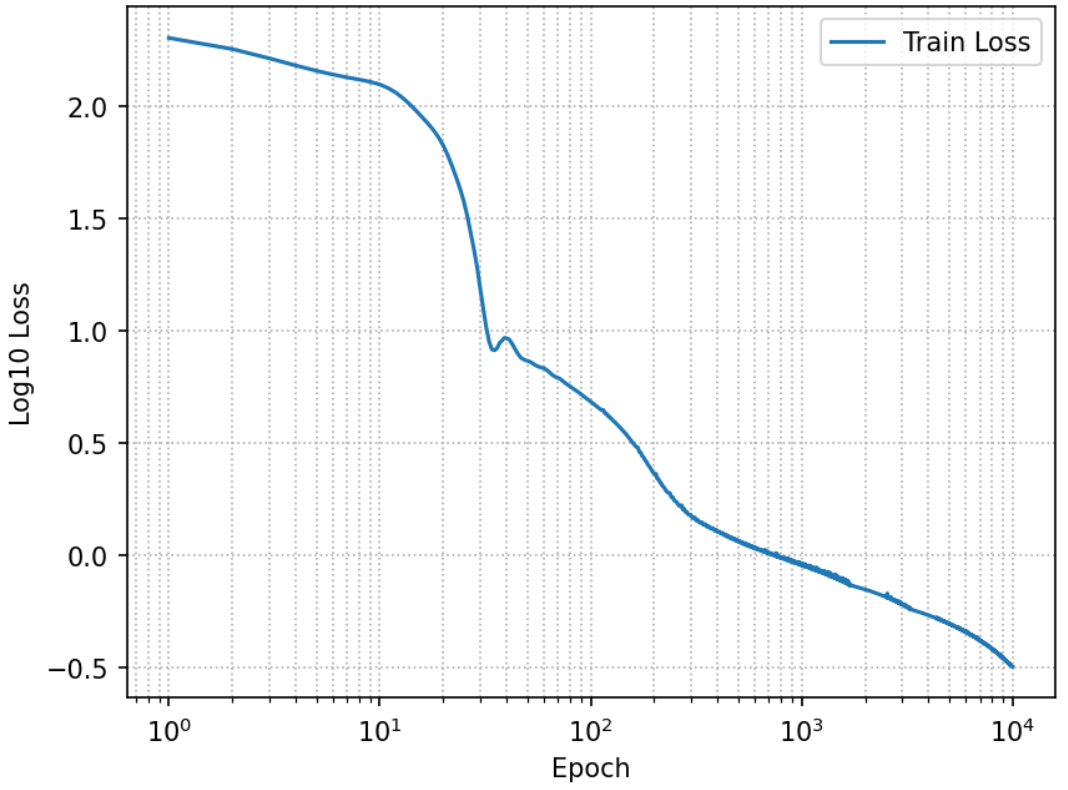}
\caption{Training loss (log scale) over epochs.}
\label{fig:inv_loss_curve}
\end{figure}
\begin{itemize}
    \item Final total log loss: \(-1.714040\)\vspace{-2mm}
    \item L2 Relative Error of the Solution: \(3.145434 \times 10^{-2}\)\vspace{-2mm}
    \item Recovered exponent \( m \): 2.3574\vspace{-2mm}
    \item True m = 3\vspace{-2mm}
    \item L2 Relative Error of m \(= 2.142 \times 10^{-1}\)
\end{itemize}
\begin{figure}[H]
\centering
\includegraphics[width=\textwidth]{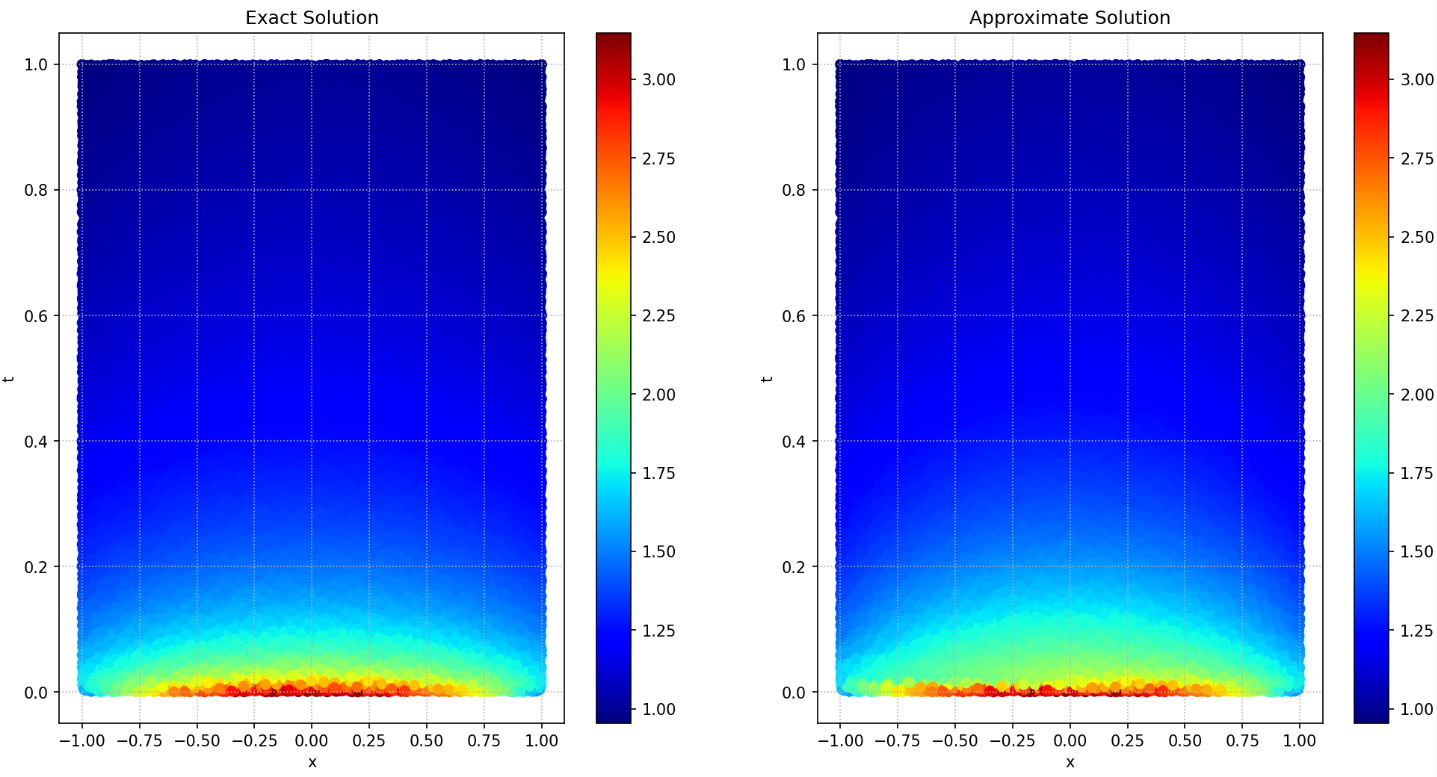}
\caption{Comparison of exact and learned solutions \( u(t,x) \) in the inverse problem.}
\label{fig:inv_solution}
\end{figure}

\subsubsection{Analysis and Comparison}
The relative error of the exponent \( m \) of the inverse problem was around 20\% when the initial guess was assumed as 2.0. This error was significantly greater than that of the classical solution.

\noindent To have a clearer picture of how much the initial guess affects accuracy, we tested the model using a range of different initial guesses for \( m \), ranging from 1.5 to 5.0 in steps of 0.5, excluding the true value.

\noindent For each initial guess, the model was trained independently and the relative error on the reconstructed \( m \) was calculated. We saw that the closer the initial guess was to the actual value \( m=3 \), the more precise the estimated exponent was. The highest precision achieved by all tests was when an initial guess of 2.5 was used, and a resulting relative error of about 9\% was achieved. But even that is significantly worse than the error achieved by the classical method.

\noindent This implies that the choice of initial guess for parameter \( m \) determines greatly the convergence and precision in finding the inverse problem of the PINN. Starting closer to the real value allows the model to learn a better approximation and increases reliability in general.

\section{Conclusion}
In this work, we successfully addressed the forward and inverse problems of the logistic equation and the porous medium equation with classical numerical approaches and Physics-Informed Neural Networks. Our results demonstrate that PINNs offer a promising platform at the intersection of applied mathematics and machine learning, providing rigorous and accurate solutions to these equations. It is also worth noting that despite being more complex tools, PINNs did not always outperform classical methods. Indeed, with proper parameterization and initial guesses, classical methods still output decent results.

\noindent Future research will focus on the PME, aiming to extend the analysis to more general cases of solutions broader than the classical Barenblatt profile , improve accuracy in inverse problems, and explore higher-dimensional geometries. In addition to PINNs, other machine learning algorithms will be investigated as potential solvers of the challenging inverse problem.

\section*{Acknowledgments}
This work originated as a second-year undergraduate project initiated by the authors and evolved over the course of a year through experimentation, refinement, and validation. We are deeply grateful to Professor Sophie Moufawad for her constant support, guidance, and encouragement throughout all the stages of this project. We also thank the Center for Advanced Mathematical Sciences (CAMS) and the Mathematics Department for providing the opportunity to initiate and develop this work. In addition, we acknowledge the contributions of Mr. Christoph Chabab and Mr. Ali Ibrahim during the 2024 Summer Research Camp. This article represents the culmination of a year of substained research and testing.

\newpage
\appendix
\addcontentsline{toc}{section}{Appendix}
\section{Testing of Classical Minimization Techniques}\label{sec:app}

This appendix is dedicated to the testing of classical minimization techniques. We begin by the logistic equation with an analytical minimizer for $r$ at each time point $t_i$. We then test numerical methods to find a constant value for the parameters when they are unknown. The cases where $r$, and $r$ and $K$ are unknown are considered, we also add white noise to our experiments in particular cases. For this purpose we use the MATLAB function \texttt{awgn} from the communication tool box (note that awgn stands for additive white noise) with a signal to noise ratio
of 100/3, and included a fallback option using additive Gaussian noise scaled by 3\% of the maximum
data magnitude. Finally, we turn our attention to the porous medium equation and apply minimization techniques to find its constant polytorpic parameter $\beta$ .\\[1mm]
We are interested in accurate predictions of the parameters. We define \textbf{feval} as the loss function evaluated using the estimated parameters. 
\subsection{Logistic equation}
\subsubsection{Analytical testing}
\label{app:analytical-logistic}

We assume that $r$ is a discrete function of time, with $r$=$r_i$ when $t$ = $t_i$. Hence, with these assumptions $r$ is a time series. We validated our approach using real-world data. Specifically, we tested on France's population data spanning the years 1950 to 2024. The parameters were set as follows:
\begin{itemize}
    \item Time interval: $t_0 = 1950$, $T = 2024$
    \item Initial population: $p_0 = p(1950)$
    \item Carrying capacity: $K = \max(p_i) + 1$, to ensure $K - p_i \neq 0$ in the denominator.
\end{itemize}

At each time point $t_i$, the growth rate $r_i$ is computed analytically using the formula:
\[
r_i = \frac{1}{t_i - t_0} \ln\left(\frac{p_i (K - p_0)}{p_0 (K - p_i)}\right)
\]

\noindent The method achieved a relative error of $1.70911 * 10^{-16}$ which is effectively zero within floating-point precision limits. The relative error (in populations) was computed as
\[
\frac{{p}(t_i; r_i) - p_{\text{i}}}{p_{\text{i}}},
\]
where ${p}(t_i; r_i)$ is the predicted population at time $t_i$ using the analytically determined $r_i$, and $p_{\text{i}}$ is the true population at time $t_i$. \\
We define the relative error in $r$ and in $k$ in the same way. 
We also denote by $feval$ the evalutation of the loss function at the given parameter(s) we found. The interpolation error and extrapolation error are defined as the evalution of the loss function at the given parameter(s) we found but on the interval of training and testing, respectfully.

\subsubsection{Classical Minimization Techniques for $r$  using simulated data}
Assume the case where $r = 0.13$, $K = 10^6$, $p_0 = 10^4$, and $t \in [0, 200]$.  
We simulate some data without noise and then with $3\%$ white noise, normalize the loss function, and then run our minimization algorithms on the first $50\%$ of the data. We normalized the loss function by dividing the mean squared error by the square of the maximum absolute value of the training data, resulting in the following expression:
\[
\text{Loss}(r) = \frac{1}{m \cdot \left( \max\limits_{1 \leq i \leq m} \lvert P_{\text{data}}(t_i) \rvert \right)^2} \sum_{i=1}^m \left( p_{\text{}}(t_i; r) - p_{\text{i}} \right)^2.
\]
Here we consider that our parameters ($r$, and $k$ later) are constants which we are after.\\[0.25mm]

\noindent We also consider multiple initial guesses, namely n*$r_0$ where n is: $25\%$, $50\%$, $75\%$, $90\%$, $110\%$, and $150\%$ of the true value of $r$, respectively. The timing were recorded by running each algorithm ten times and taking the average time to run. 
In table \ref{tab:nonoise} we show the obtained results for varying $n$ (shown in the )  consider the n from n*$r_0$ in the table. We get the following results.\vspace{-4mm}
\subsubsection{Without noise:}
\begin{table}[H]
\centering
\caption{Optimization results without noise}\label{tab:nonoise}
\renewcommand{\arraystretch}{1.15}
\begin{tabular}{ccccccc}
\toprule
\multicolumn{7}{c}{\textbf{fminunc}}\\
\midrule
$nr_0$ & iter & $r$ & $\lvert r-r_{\text{exact}}\rvert/r_{\text{exact}}$ & feval & Interp.\ error & Extrap.\ error \\
\midrule
$0.25$ & 3 & $1.200\,\mathrm{E}{-1}$ & $1.976\,\mathrm{E}{-7}$ & $1.111\,\mathrm{E}{-14}$ & $1.397\,\mathrm{E}{-7}$ & $1.101\,\mathrm{E}{-10}$ \\
0.5 & 7 & $1.300\mathrm{E}{-1}$ & $5.465\mathrm{E}{-7}$ & $8.500\mathrm{E}{-14}$ & $3.864\mathrm{E}{-7}$ & 3$.047\mathrm{E}{-10}$ \\
0.75 & 6 & $1.300\mathrm{E}{-1}$ & $6.360\mathrm{E}{-8}$ & $1.151\mathrm{E}{-15}$ & $4.496\mathrm{E}{-8}$ & $3.545\mathrm{E}{-11}$ \\
0.9 & 5 & $1.300\mathrm{E}{-1}$ & $5.630\mathrm{E}{-8}$ & $9.020\mathrm{E}{-16}$ & $3.981\mathrm{E}{-8}$ & $3.135\mathrm{E}{-11}$ \\
1.10 & 4 & $1.300\mathrm{E}{-1}$ & $5.918\mathrm{E}{-8}$ & $9.965\mathrm{E}{-16}$ & $4.184\mathrm{E}{-8}$ & $3.299\mathrm{E}{-11}$ \\
1.5 & 4 & $1.300\mathrm{E}{-1}$ & $5.717\mathrm{E}{-8}$ & $9.300\mathrm{E}{-16}$ & $4.042\mathrm{E}{-8}$ & $3.187\mathrm{E}{-11}$ \\
\midrule

\multicolumn{7}{c}{\textbf{fmincon}}\\
\midrule
$nr_0$ & iter & $r$ & $\lvert r-r_{\text{exact}}\rvert/r_{\text{exact}}$ & feval & Interp.\ error & Extrap.\ error \\\midrule
0.25 & 10 & $1.300\mathrm{E}{-1}$ & $5.958\mathrm{E}{-7}$ & $1.009\mathrm{E}{-13}$ & $4.211\mathrm{E}{-7}$ & $3.320\mathrm{E}{-10}$ \\
0.5 & 11 & $1.300\mathrm{E}{-1}$ & $2.217\mathrm{E}{-8}$ & $1.399\mathrm{E}{-16}$ & $1.567\mathrm{E}{-8}$ & $1.236\mathrm{E}{-11}$ \\
0.75 & 14 & $1.300\mathrm{E}{-1}$ & $6.444\mathrm{E}{-7}$ & $1.182\mathrm{E}{-13}$ & $4.556\mathrm{E}{-7}$ & $3.592\mathrm{E}{-10}$ \\
0.9 & 12 & $1.300\mathrm{E}{-1}$ & $6.416\mathrm{E}{-7} $&$ 1.171\mathrm{E}{-13}$ & $4.539\mathrm{E}{-7}$& $3.576\mathrm{E}{-10} $\\
1.1 & 7  & $1.300\mathrm{E}{-1}$ & $2.217\mathrm{E}{-8}$ & $1.399\mathrm{E}{-16}$ & $1.568\mathrm{E}{-8}$ & $1.236\mathrm{E}{-11}$ \\
1.5 & 9  & $1.300\mathrm{E}{-1}$ & $2.217\mathrm{E}{-8}$ &$ 1.398\mathrm{E}{-16} $&$ 1.567\mathrm{E}{-8} $& $1.236\mathrm{E}{-11}$ \\
\midrule

\multicolumn{7}{c}{\textbf{Secant}}\\
\midrule
$nr_0$ & iter & $r$ & $\lvert r-r_{\text{exact}}\rvert/r_{\text{exact}}$ & feval & Interp.\ error & Extrap.\ error \\ \midrule
0.25 & 50 & $-1.824\mathrm{E}{1}$ & $1.413\mathrm{E}{2}$ & $5.692\mathrm{E}{-1}$ & $1.000\mathrm{E}{0} $& $1.000\mathrm{E}{0}$ \\
0.5 & 8  & $1.300\mathrm{E}{-1}$ & $5.207\mathrm{E}{-11}$ & $7.713\mathrm{E}{-22}$ & $3.681\mathrm{E}{-11}$ & $2.902\mathrm{E}{-14}$ \\
0.75 & 6  & $1.300\mathrm{E}{-1}$ & $2.839\mathrm{E}{-10}$ & $2.929\mathrm{E}{-20}$ & $2.007\mathrm{E}{-10}$ & $1.582\mathrm{E}{-13}$ \\
0.9 & 4  & $1.300\mathrm{E}{-1} $& $2.021\mathrm{E}{-11}$ & $1.162\mathrm{E}{-22}$ & $1.429\mathrm{E}{-11}$ & $1.126\mathrm{E}{-14}$ \\
1.1 & 7  & $1.300\mathrm{E}{-1}$ & $8.711\mathrm{E}{-14}$ & $2.159\mathrm{E}{-27}$ & $6.159\mathrm{E}{-14}$ & $1.259\mathrm{E}{-16}$ \\
\bottomrule
\end{tabular}
\end{table}

\begin{table}[H]
\centering
\caption{Optimization results without noise (additional solvers)}
\renewcommand{\arraystretch}{1.15}
\begin{tabular}{ccccccc}
\toprule
\multicolumn{7}{c}{\textbf{Newton with symbolic diff (NewtonS)}}\\
\midrule
$nr_0$ & iter & $r$ & $\lvert r-r_{\text{exact}}\rvert/r_{\text{exact}}$ & feval & Interp.\ error & Extrap.\ error \\ \midrule
0.25 & 3  & $-3.000\mathrm{E}{-03}$ & $1.024\mathrm{E}{+00}$ & $5.590\mathrm{E}{-01}$ & $9.910\mathrm{E}{-01}$ & $9.940\mathrm{E}{-01}$ \\
0.5 & 7  & $1.300\mathrm{E}{-01}$  & $5.503\mathrm{E}{-12}$ & $8.616\mathrm{E}{-24}$ & $3.891\mathrm{E}{-12}$ & $3.085\mathrm{E}{-15}$ \\
0.75 & 5  & $1.300\mathrm{E}{-01}$  & $7.150\mathrm{E}{-07}$ & $1.455\mathrm{E}{-13}$ & $5.055\mathrm{E}{-07}$ & $3.986\mathrm{E}{-10}$ \\
0.9 & 5  & $1.300\mathrm{E}{-01}$  & $3.516\mathrm{E}{-12}$ & $3.518\mathrm{E}{-24}$ & $2.486\mathrm{E}{-12}$ & $1.956\mathrm{E}{-15}$ \\
1.1 & 5  & $1.300\mathrm{E}{-01}$  & $9.852\mathrm{E}{-11}$ & $2.762\mathrm{E}{-21}$ & $6.965\mathrm{E}{-11}$ & $5.491\mathrm{E}{-14}$ \\
1.5 & 2  & $-5.644\mathrm{E}{+00}$ & $4.442\mathrm{E}{+01}$ & $5.690\mathrm{E}{-01}$ & $1.000\mathrm{E}{+00}$ & $1.000\mathrm{E}{+00}$ \\
\midrule

\multicolumn{7}{c}{\textbf{Newton with exact diff (NewtonE)}}\\
\midrule
$nr_0$ & iter & $r$ & $\lvert r-r_{\text{exact}}\rvert/r_{\text{exact}}$ & feval & Interp.\ error & Extrap.\ error \\ \midrule
0.25 & 3  & $-3.000\mathrm{E}{-03}$ & $1.024\mathrm{E}{+00}$ & $5.590\mathrm{E}{-01}$ & $9.910\mathrm{E}{-01}$ & $9.940\mathrm{E}{-01}$ \\
0.5 & 7  & $1.300\mathrm{E}{-01}$  & $5.503\mathrm{E}{-12}$ & $8.616\mathrm{E}{-24}$ & $3.891\mathrm{E}{-12}$ & $3.085\mathrm{E}{-15}$ \\
0.7 & 5  & $1.300\mathrm{E}{-01}$  & $7.150\mathrm{E}{-07}$ & $1.455\mathrm{E}{-13}$ & $5.055\mathrm{E}{-07}$ & $3.986\mathrm{E}{-10}$ \\
0.9 & 5  & $1.300\mathrm{E}{-01}$  & $3.516\mathrm{E}{-12}$ & $3.518\mathrm{E}{-24}$ & $2.486\mathrm{E}{-12}$ & $1.956\mathrm{E}{-15}$ \\
1.1 & 5  & $1.300\mathrm{E}{-01}$  & $9.852\mathrm{E}{-11}$ & $2.762\mathrm{E}{-21}$ & $6.965\mathrm{E}{-11}$ & $5.491\mathrm{E}{-14}$ \\
1.5 & 2  & $-5.644\mathrm{E}{+00}$ & $4.442\mathrm{E}{+01}$ & $5.690\mathrm{E}{-01}$ & $1.000\mathrm{E}{+00}$ & $1.000\mathrm{E}{+00}$ \\
\midrule

\multicolumn{7}{c}{\textbf{Steepest Descent (SD)}}\\
\midrule
$nr_0$ & iter & $r$ & $\lvert r-r_{\text{exact}}\rvert/r_{\text{exact}}$ & feval & Interp.\ error & Extrap.\ error \\ \midrule
0.25 & 10 & $1.300\mathrm{E}{-01}$ & $4.045\mathrm{E}{-08}$ & $4.656\mathrm{E}{-16}$ & $2.860\mathrm{E}{-08}$ & $2.255\mathrm{E}{-11}$ \\
0.5 & 7  & $1.300\mathrm{E}{-01}$ & $4.101\mathrm{E}{-08}$ & $4.786\mathrm{E}{-16}$ & $2.900\mathrm{E}{-08}$ & $2.286\mathrm{E}{-11}$ \\
0.75 & 8  & $1.300\mathrm{E}{-01}$ & $3.850\mathrm{E}{-08}$ & $4.218\mathrm{E}{-16}$ & $2.722\mathrm{E}{-08}$ & $2.146\mathrm{E}{-11}$ \\
0.9 & 6  & $1.300\mathrm{E}{-01}$ & $9.134\mathrm{E}{-08}$ & $2.374\mathrm{E}{-15}$ & $6.458\mathrm{E}{-08}$ & $5.091\mathrm{E}{-11}$ \\
1.1 & 7  & $1.300\mathrm{E}{-01}$ & $3.599\mathrm{E}{-08}$ & $3.685\mathrm{E}{-16}$ & $2.544\mathrm{E}{-08}$ & $2.006\mathrm{E}{-11}$ \\
1.5 & 9  & $1.300\mathrm{E}{-01}$ & $2.428\mathrm{E}{-08}$ & $1.678\mathrm{E}{-16}$ & $1.717\mathrm{E}{-08}$ & $1.354\mathrm{E}{-11}$ \\
\bottomrule
\end{tabular}
\end{table}

\begin{table}[h!]
\centering
\begin{tabular}{|c|c|c|c|c|c|c|}
\hline
$nr_0$ & \textbf{fminunc} & \textbf{fmincon} & \textbf{Secant} & \textbf{NewtonS} & \textbf{NewtonE} & \textbf{SD} \\
\hline
0.25 & $1.238\mathrm{E{-}03}$ & $8.849\mathrm{E{-}03}$ & $1.588\mathrm{E{-}03}$ & $1.301\mathrm{E{-}02}$ & $1.833\mathrm{E{-}04}$ & $2.944$ \\
0.5 & $1.389\mathrm{E{-}03}$ & $1.010\mathrm{E{-}02}$ & $2.403\mathrm{E{-}04}$ & $4.552\mathrm{E{-}04}$ & $5.035\mathrm{E{-}04}$ & $1.801$ \\
0.75 & $7.652\mathrm{E{-}04}$ & $7.296\mathrm{E{-}03}$ & $1.315\mathrm{E{-}04}$ & $4.075\mathrm{E{-}04}$ & $4.215\mathrm{E{-}04}$ & $1.994$ \\
0.9 & $6.263\mathrm{E{-}04}$ & $7.064\mathrm{E{-}03}$ & $8.965\mathrm{E{-}05}$ & $3.989\mathrm{E{-}04}$ & $3.789\mathrm{E{-}04}$ & $1.480$ \\
1.1 & $6.322\mathrm{E{-}04}$ & $7.920\mathrm{E{-}03}$ & $3.872\mathrm{E{-}04}$ & $3.198\mathrm{E{-}04}$ & $2.813\mathrm{E{-}04}$ & $1.721$ \\
1.5 & $6.338\mathrm{E{-}04}$ & $5.625\mathrm{E{-}03}$ & $2.172\mathrm{E{-}04}$ & $1.666\mathrm{E{-}04}$ & $1.872\mathrm{E{-}04}$ & $2.153$ \\
\hline
\end{tabular}
\caption{Convergence time (in seconds)}
\end{table}

\subsubsection{With 3\% white noise added:}

\begin{table}[H]
\centering
\caption{Optimization results with $3\%$ noise}
\renewcommand{\arraystretch}{1.15}
\begin{tabular}{ccccccc}
\toprule
\multicolumn{7}{c}{\textbf{fminunc}}\\
\midrule
$nr_0$ iter & $r$ & $\lvert r-r_{\text{exact}}\rvert/r_{\text{exact}}$ & feval & Interp.\ error & Extrap.\ error \\ \midrule
0.25 & 3 & $1.300\mathrm{E{-}01}$ & $1.995\mathrm{E{-}07}$ & $1.151\mathrm{E{-}14}$ & $1.422\mathrm{E{-}07}$ & $2.276\mathrm{E{-}08}$ \\
0.5 & 7 & $1.300\mathrm{E{-}01}$ & $5.474\mathrm{E{-}07}$ & $8.545\mathrm{E{-}14}$ & $3.875\mathrm{E{-}07}$ & $2.346\mathrm{E{-}08}$ \\
0.75 & 6 & $1.300\mathrm{E{-}01}$ & $6.352\mathrm{E{-}08}$ & $1.548\mathrm{E{-}15}$ & $5.215\mathrm{E{-}08}$ & $2.108\mathrm{E{-}08}$ \\
0.9 & 5 & $1.300\mathrm{E{-}01}$ & $5.352\mathrm{E{-}08}$ & $1.367\mathrm{E{-}15}$ & $4.901\mathrm{E{-}08}$ & $2.135\mathrm{E{-}08}$ \\
1.1 & 4 & $1.300\mathrm{E{-}01}$ & $5.862\mathrm{E{-}08}$ & $1.443\mathrm{E{-}15}$ & $5.035\mathrm{E{-}08}$ & $1.992\mathrm{E{-}08}$ \\
1.5 & 4 & $1.300\mathrm{E{-}01}$ & $5.406\mathrm{E{-}08}$ & $1.340\mathrm{E{-}15}$ & $4.851\mathrm{E{-}08}$ & $2.129\mathrm{E{-}08}$ \\
\midrule

\multicolumn{7}{c}{\textbf{fmincon}}\\
\midrule
$nr_0$ & iter & $r$ & $\lvert r-r_{\text{exact}}\rvert/r_{\text{exact}}$ & feval & Interp.\ error & Extrap.\ error \\ \midrule
0.25 & 10 & $1.300\mathrm{E{-}01}$ & $5.937\mathrm{E{-}07}$ & $1.013\mathrm{E{-}13}$ & $4.219\mathrm{E{-}07}$ & $2.280\mathrm{E{-}08}$ \\
0.5 & 11 & $1.300\mathrm{E{-}01}$ & $2.136\mathrm{E{-}08}$ & $5.912\mathrm{E{-}16}$ & $3.223\mathrm{E{-}08}$ & $2.352\mathrm{E{-}08}$ \\
0.75 & 14 & $1.300\mathrm{E{-}01}$ & $6.445\mathrm{E{-}07}$ & $1.186\mathrm{E{-}13}$ & $4.564\mathrm{E{-}07}$ & $2.107\mathrm{E{-}08}$ \\
0.9 & 12 & $1.300\mathrm{E{-}01}$ & $6.443\mathrm{E{-}07}$ & $1.176\mathrm{E{-}13}$ & $4.545\mathrm{E{-}07}$ & $2.140\mathrm{E{-}08}$ \\
1.1 & 7  & $1.300\mathrm{E{-}01}$ & $2.161\mathrm{E{-}08}$ & $5.864\mathrm{E{-}16}$ & $3.210\mathrm{E{-}08}$ & $1.992\mathrm{E{-}08}$ \\
1.5 & 9  & $1.300\mathrm{E{-}01}$ & $1.905\mathrm{E{-}08}$ & $5.494\mathrm{E{-}16}$ & $3.107\mathrm{E{-}08}$ & $2.129\mathrm{E{-}08}$ \\
\midrule

\multicolumn{7}{c}{\textbf{Secant}}\\
\midrule
$nr_0$ & iter & $r$ & $\lvert r-r_{\text{exact}}\rvert/r_{\text{exact}}$ & feval & Interp.\ error & Extrap.\ error \\ \midrule
0.25 & 50 & $-1.824\mathrm{E{+}01}$ & $1.413\mathrm{E{+}02}$ & $5.692\mathrm{E{-}01}$ & $1.000\mathrm{E{+}00}$ & $1.000$ \\
0.5 & 8  & $1.300\mathrm{E{-}01}$ & $7.598\mathrm{E{-}10}$ & $4.513\mathrm{E{-}16}$ & $2.816\mathrm{E{-}08}$ & $2.351\mathrm{E{-}08}$ \\
0.75 & 6  & $1.300\mathrm{E{-}01}$ & $2.027\mathrm{E{-}10}$ & $3.975\mathrm{E{-}16}$ & $2.643\mathrm{E{-}08}$ & $2.107\mathrm{E{-}08}$ \\
0.9 & 4  & $1.300\mathrm{E{-}01}$ & $2.767\mathrm{E{-}09}$ & $4.654\mathrm{E{-}16}$ & $2.859\mathrm{E{-}08}$ & $2.135\mathrm{E{-}08}$ \\
1.1 & 7  & $1.300\mathrm{E{-}01}$ & $5.656\mathrm{E{-}10}$ & $4.466\mathrm{E{-}16}$ & $2.801\mathrm{E{-}08}$ & $1.991\mathrm{E{-}08}$ \\
\bottomrule
\end{tabular}
\end{table}

\begin{table}[h!]
\centering
\begin{tabular}{|c|c|c|c|c|c|c|}
\hline
$nro$ & 
\textbf{fminunc} & \textbf{fmincon} & \textbf{Secant} & \textbf{NewtonS} & \textbf{NewtonE} & \textbf{SD} \\
\hline
0.25& 
$8.254\mathrm{E{-}04}$ & $6.130\mathrm{E{-}03}$ & $1.118\mathrm{E{-}03}$ & $1.001\mathrm{E{-}02}$ & $1.352\mathrm{E{-}04}$ & $2.638\mathrm{E{+}00}$ \\
0.5&
$8.158\mathrm{E{-}04}$ & $1.087\mathrm{E{-}02}$ & $1.769\mathrm{E{-}04}$ & $3.429\mathrm{E{-}04}$ & $3.720\mathrm{E{-}04}$ & $1.816\mathrm{E{+}00}$ \\
0.75&
$6.854\mathrm{E{-}04}$ & $7.862\mathrm{E{-}03}$ & $1.354\mathrm{E{-}04}$ & $2.483\mathrm{E{-}04}$ & $2.596\mathrm{E{-}04}$ & $2.004\mathrm{E{+}00}$ \\
0.9&
$6.594\mathrm{E{-}04}$ & $7.635\mathrm{E{-}03}$ & $8.859\mathrm{E{-}05}$ & $2.477\mathrm{E{-}04}$ & $2.665\mathrm{E{-}04}$ & $1.572\mathrm{E{+}00}$ \\
1.1&
$5.858\mathrm{E{-}04}$ & $6.007\mathrm{E{-}03}$ & $2.270\mathrm{E{-}04}$ & $2.568\mathrm{E{-}04}$ & $2.506\mathrm{E{-}04}$ & $1.882\mathrm{E{+}00}$ \\
1.5&
$6.069\mathrm{E{-}04}$ & $5.653\mathrm{E{-}03}$ & $2.161\mathrm{E{-}04}$ & $1.574\mathrm{E{-}02}$ & $1.057\mathrm{E{-}04}$ & $2.504\mathrm{E{+}00}$ \\
\hline
\end{tabular}
\caption{Convergence time (in seconds) when adding 3\% white noise}
\end{table} 

\begin{table}[H]
\centering
\caption{Optimization results with $3\%$ noise}
\renewcommand{\arraystretch}{1.15}
\begin{tabular}{ccccccc}
\toprule
\multicolumn{7}{c}{\textbf{NewtonS}}\\
\midrule
$nr_0$ & iter & $r$ & $\lvert r - r_{\text{exact}} \rvert / r_{\text{exact}}$ & feval & Interp.\ error & Extrap.\ error \\ \midrule
0.25 & 3 & $-3.070\mathrm{E{-}03}$ & $1.024\mathrm{E{+}00}$ & $5.587\mathrm{E{-}01}$ & $9.908\mathrm{E{-}01}$ & $9.937\mathrm{E{-}01}$ \\
0.5 & 7 & $1.300\mathrm{E{-}01}$ & $8.064\mathrm{E{-}10}$ & $4.513\mathrm{E{-}16}$ & $2.816\mathrm{E{-}08}$ & $2.351\mathrm{E{-}08}$ \\
0.75 & 5 & $1.300\mathrm{E{-}01}$ & $7.149\mathrm{E{-}07}$ & $1.459\mathrm{E{-}13}$ & $5.062\mathrm{E{-}07}$ & $2.109\mathrm{E{-}08}$ \\
0.9 & 5 & $1.300\mathrm{E{-}01}$ & $2.784\mathrm{E{-}09}$ & $4.654\mathrm{E{-}16}$ & $2.859\mathrm{E{-}08}$ & $2.135\mathrm{E{-}08}$ \\
1.1 & 5 & $1.300\mathrm{E{-}01}$ & $4.672\mathrm{E{-}10}$ & $4.466\mathrm{E{-}16}$ & $2.801\mathrm{E{-}08}$ & $1.991\mathrm{E{-}08}$ \\
1.5 & 2 & $-5.644\mathrm{E{+}00}$ & $4.442\mathrm{E{+}01}$ & $5.692\mathrm{E{-}01}$ & $1.000\mathrm{E{+}00}$ & $1.000\mathrm{E{+}00}$ \\
\midrule

\multicolumn{7}{c}{\textbf{NewtonE}}\\
\midrule
$nr_0$ & iter & $r$ & $\lvert r - r_{\text{exact}} \rvert / r_{\text{exact}}$ & feval & Interp.\ error & Extrap.\ error \\ \midrule
0.25 & 3 & $-3.070\mathrm{E{-}03}$ & $1.024\mathrm{E{+}00}$ & $5.587\mathrm{E{-}01}$ & $9.908\mathrm{E{-}01}$ & $9.937\mathrm{E{-}01}$ \\
0.5 & 7 & $1.300\mathrm{E{-}01}$ & $8.064\mathrm{E{-}10}$ & $4.513\mathrm{E{-}16}$ & $2.816\mathrm{E{-}08}$ & $2.351\mathrm{E{-}08}$ \\
0.75 & 5 & $1.300\mathrm{E{-}01}$ & $7.149\mathrm{E{-}07}$ & $1.459\mathrm{E{-}13}$ & $5.062\mathrm{E{-}07}$ & $2.109\mathrm{E{-}08}$ \\
0.9 & 5 & $1.300\mathrm{E{-}01}$ & $2.784\mathrm{E{-}09}$ & $4.654\mathrm{E{-}16}$ & $2.859\mathrm{E{-}08}$ & $2.135\mathrm{E{-}08}$ \\
1.1 & 5 & $1.300\mathrm{E{-}01}$ & $4.672\mathrm{E{-}10}$ & $4.466\mathrm{E{-}16}$ & $2.801\mathrm{E{-}08}$ & $1.991\mathrm{E{-}08}$ \\
1.5 & 2 & $-5.644\mathrm{E{+}00}$ & $4.442\mathrm{E{+}01}$ & $5.692\mathrm{E{-}01}$ & $1.000\mathrm{E{+}00}$ & $1.000\mathrm{E{+}00}$ \\
\midrule

\multicolumn{7}{c}{\textbf{Steepest descent}}\\
\midrule
$nr_0$ & iter & $r$ & $\lvert r - r_{\text{exact}} \rvert / r_{\text{exact}}$ & feval & Interp.\ error & Extrap.\ error \\ \midrule
0.25 & 10 & $1.300\mathrm{E{-}01}$ & $3.861\mathrm{E{-}08}$ & $8.608\mathrm{E{-}16}$ & $3.889\mathrm{E{-}08}$ & $2.277\mathrm{E{-}08}$ \\
0.5 & 7  & $1.300\mathrm{E{-}01}$ & $4.183\mathrm{E{-}08}$ & $9.299\mathrm{E{-}16}$ & $4.042\mathrm{E{-}08}$ & $2.351\mathrm{E{-}08}$ \\
0.75 & 8  & $1.300\mathrm{E{-}01}$ & $3.859\mathrm{E{-}08}$ & $8.193\mathrm{E{-}16}$ & $3.794\mathrm{E{-}08}$ & $2.107\mathrm{E{-}08}$ \\
0.9 & 6  & $1.300\mathrm{E{-}01}$ & $8.855\mathrm{E{-}08}$ & $2.839\mathrm{E{-}15}$ & $7.062\mathrm{E{-}08}$ & $2.135\mathrm{E{-}08}$ \\
1.1 & 7  & $1.300\mathrm{E{-}01}$ & $3.655\mathrm{E{-}08}$ & $8.151\mathrm{E{-}16}$ & $3.784\mathrm{E{-}08}$ & $1.991\mathrm{E{-}08}$ \\
1.5 & 9  & $1.300\mathrm{E{-}01}$ & $2.740\mathrm{E{-}08}$ & $5.773\mathrm{E{-}16}$ & $3.185\mathrm{E{-}08}$ & $2.129\mathrm{E{-}08}$ \\
\bottomrule
\end{tabular}
\end{table}

We can see that all algorithms performed well overall. When the initial guess is far from the true value, the Newton methods and the secant do not converge to the true minimizer. However, when they do, it takes them the least amount of time. The steepest descent gave us good results but took (by far) the most time. \texttt{fminunc} and \texttt{fmincon} gave us good results with little time, making their performance faster than the steepest descent, but slower than the Newton and secant. They did, however, always converge to the true minimizer.

\subsubsection{Results for $r$ \& $k$:}\label{sec:rk}

The loss function is not convex in $k$, which makes it difficult for the steepest descent, \texttt{fminunc} \& \texttt{fmincon} to find the real $r$, as they are more prone to errors. Since the loss is not convex, the secant method \& the Newton's method are totally unable to find the solution. We do run some tests on data generated using $r$ = 0.13 and $k = 10^6$ and plot the loss. We see a ``plateau" near our true minimizer.\vspace{-5mm}

\begin{figure}[H]
    \centering
    \begin{minipage}{0.45\textwidth}
        \centering
\includegraphics[width=\linewidth]{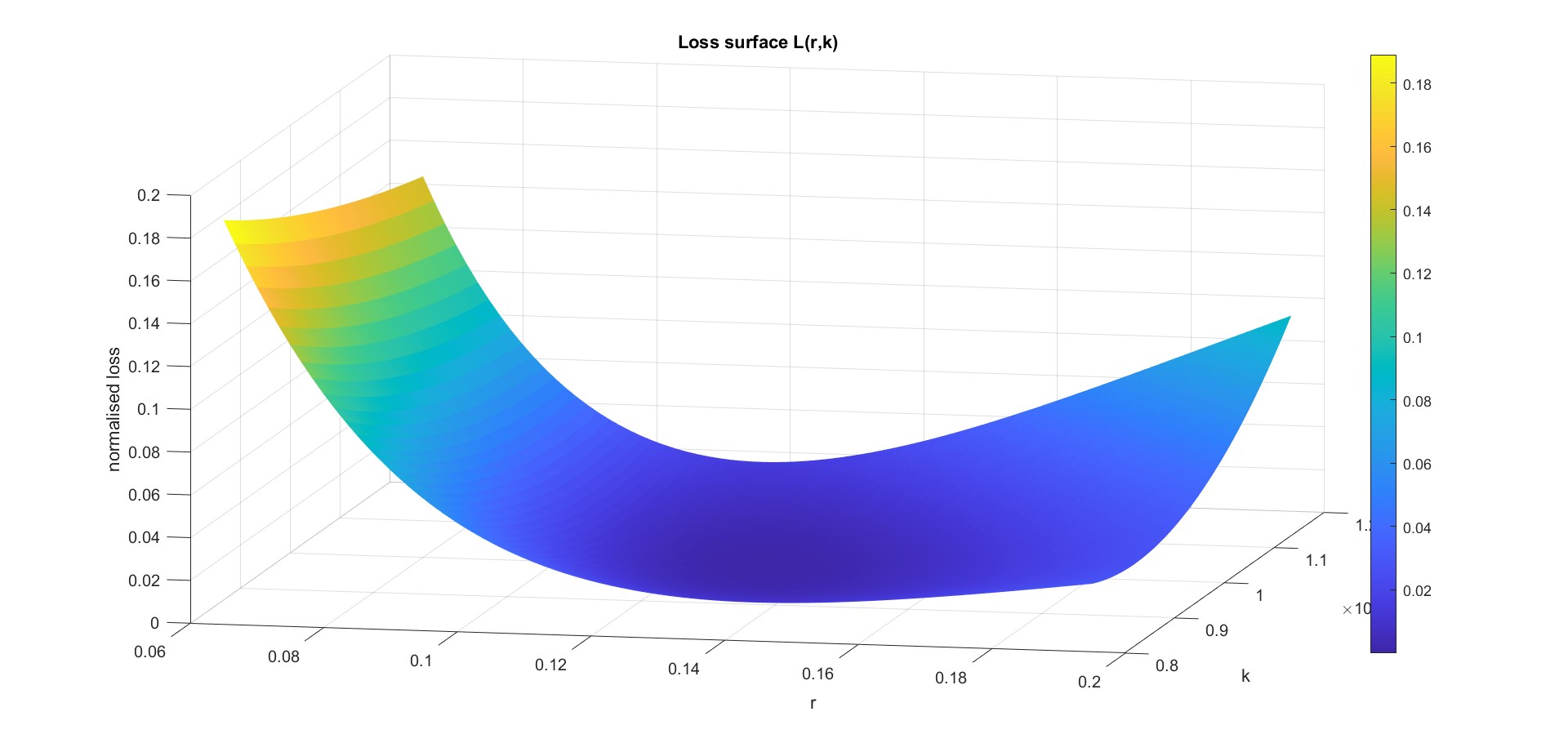}
    \end{minipage}%
    \hspace{0.01\textwidth}
    \begin{minipage}{0.45\textwidth}
        \centering
        \includegraphics[width=\linewidth]{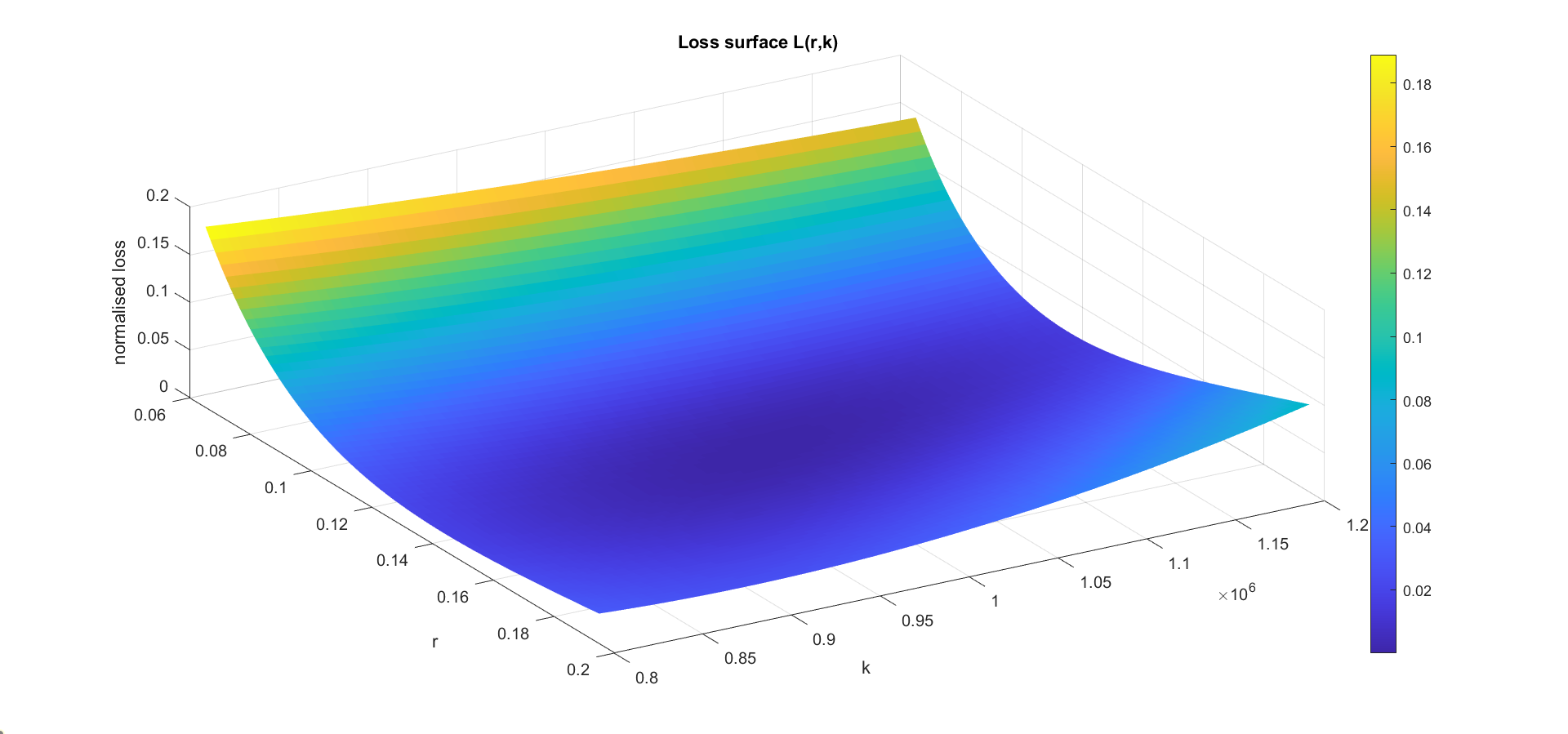}
    \end{minipage}
    \caption{Loss surface seen from different angles. The data has been generated using $r$ =0.13 and $k$ = $10^6$.}
    \label{fig:convexity_combined}
\end{figure}

\begin{table}[H]
    \centering
\caption{Optimization results for $r$ and $k$ without noise}
\begin{tabular}{cccccccc}
\toprule
\multicolumn{8}{c}{\textbf{fminunc}}\\
\midrule
$nr_0$ & iter & $r$ & $k$ & $\lvert r - r_{\text{exact}} \rvert / r_{\text{exact}}$ & $\lvert k - k_{\text{exact}} \rvert / k_{\text{exact}}$ & Interp.\ error & Extrap.\ error \\
\midrule
0.25 & 3 & $2.038\mathrm{E{-}01}$ & $2.500\mathrm{E{+}05}$ & $5.679\mathrm{E{-}01}$ & $7.500\mathrm{E{-}01}$ & $7.269\mathrm{E{+}05}$ & $7.498\mathrm{E{+}05}$ \\
0.5 & 3 & $1.693\mathrm{E{-}01}$ & $5.000\mathrm{E{+}05}$ & $3.020\mathrm{E{-}01}$ & $5.000\mathrm{E{-}01}$ & $4.691\mathrm{E{+}05}$ & $4.999\mathrm{E{+}05}$ \\
0.75 & 5 & $1.480\mathrm{E{-}01}$ & $7.500\mathrm{E{+}05}$ & $1.382\mathrm{E{-}01}$ & $2.500\mathrm{E{-}01}$ & $2.258\mathrm{E{+}05}$ & $2.499\mathrm{E{+}05}$ \\
0.9 & 6 & $1.370\mathrm{E{-}01}$ & $9.000\mathrm{E{+}05}$ & $5.421\mathrm{E{-}02}$ & $1.000\mathrm{E{-}01}$ & $8.790\mathrm{E{+}04}$ & $9.997\mathrm{E{+}04}$ \\
1.1 & 5 & $1.229\mathrm{E{-}01}$ & $1.100\mathrm{E{+}06}$ & $5.450\mathrm{E{-}02}$ & $1.000\mathrm{E{-}01}$ & $8.402\mathrm{E{+}04}$ & $9.995\mathrm{E{+}04}$ \\
1.5 & 6 & $9.103\mathrm{E{-}02}$ & $1.500\mathrm{E{+}06}$ & $2.998\mathrm{E{-}01}$ & $5.000\mathrm{E{-}01}$ & $3.573\mathrm{E{+}05}$ & $4.973\mathrm{E{+}05}$ \\
\midrule

\multicolumn{8}{c}{\textbf{fmincon}}\\
\midrule
0.25 & 15 & $2.038\mathrm{E{-}01}$ & $2.500\mathrm{E{+}05}$ & $5.679\mathrm{E{-}01}$ & $7.500\mathrm{E{-}01}$ & $7.269\mathrm{E{+}05}$ & $7.498\mathrm{E{+}05}$ \\
0.5 & 17 & $1.693\mathrm{E{-}01}$ & $5.000\mathrm{E{+}05}$ & $3.020\mathrm{E{-}01}$ & $5.000\mathrm{E{-}01}$ & $4.691\mathrm{E{+}05}$ & $4.999\mathrm{E{+}05}$ \\
0.75 & 14 & $1.480\mathrm{E{-}01}$ & $7.500\mathrm{E{+}05}$ & $1.382\mathrm{E{-}01}$ & $2.500\mathrm{E{-}01}$ & $2.258\mathrm{E{+}05}$ & $2.499\mathrm{E{+}05}$ \\
0.9 & 13 & $1.370\mathrm{E{-}01}$ & $9.000\mathrm{E{+}05}$ & $5.422\mathrm{E{-}02}$ & $1.000\mathrm{E{-}01}$ & $8.790\mathrm{E{+}04}$ & $9.997\mathrm{E{+}04}$ \\
1.1 & 36 & $1.301\mathrm{E{-}01}$ & $9.991\mathrm{E{+}05}$ & $5.060\mathrm{E{-}04}$ & $9.341\mathrm{E{-}04}$ & $8.042\mathrm{E{+}02}$ & $9.338\mathrm{E{+}02}$ \\
1.5 & 36 & $1.301\mathrm{E{-}01}$ & $9.991\mathrm{E{+}05}$ & $5.060\mathrm{E{-}04}$ & $9.341\mathrm{E{-}04}$ & $8.042\mathrm{E{+}02}$ & $9.338\mathrm{E{+}02}$ \\
\midrule

\multicolumn{8}{c}{\textbf{Steepest Descent}}\\
\midrule
0.25 & 200 & $1.446\mathrm{E{-}01}$ & $2.500\mathrm{E{+}05}$ & $1.121\mathrm{E{-}01}$ & $7.500\mathrm{E{-}01}$ & $7.301\mathrm{E{+}05}$ & $7.498\mathrm{E{+}05}$ \\
0.5 & 200 & $1.574\mathrm{E{-}01}$ & $5.000\mathrm{E{+}05}$ & $2.105\mathrm{E{-}01}$ & $5.000\mathrm{E{-}01}$ & $4.698\mathrm{E{+}05}$ & $4.999\mathrm{E{+}05}$ \\
0.75 & 200 & $1.470\mathrm{E{-}01}$ & $7.500\mathrm{E{+}05}$ & $1.306\mathrm{E{-}01}$ & $2.500\mathrm{E{-}01}$ & $2.258\mathrm{E{+}05}$ & $2.499\mathrm{E{+}05}$ \\
0.9 & 200 & $1.370\mathrm{E{-}01}$ & $9.000\mathrm{E{+}05}$ & $5.360\mathrm{E{-}02}$ & $1.000\mathrm{E{-}01}$ & $8.790\mathrm{E{+}04}$ & $9.997\mathrm{E{+}04}$ \\
1.1 & 200 & $1.229\mathrm{E{-}01}$ & $1.100\mathrm{E{+}06}$ & $5.447\mathrm{E{-}02}$ & $1.000\mathrm{E{-}01}$ & $8.402\mathrm{E{+}04}$ & $9.995\mathrm{E{+}04}$ \\
1.5 & 177 & $9.103\mathrm{E{-}02}$ & $1.500\mathrm{E{+}06}$ & $2.998\mathrm{E{-}01}$ & $5.000\mathrm{E{-}01}$ & $3.573\mathrm{E{+}05}$ & $4.973\mathrm{E{+}05}$ \\
\bottomrule
\end{tabular}
\end{table}

\begin{table}[h!]
\centering
\begin{tabular}{|c|c|c|c|}
\hline
$nr_0$&
\textbf{fminunc} & \textbf{fmincon} & \textbf{SD} \\
\hline
0.25&
$2.053\mathrm{E{-}03}$ & $8.941\mathrm{E{-}03}$ & $4.496\mathrm{E{-}03}$ \\
0.5&
$1.601\mathrm{E{-}03}$ & $1.050\mathrm{E{-}02}$ & $4.957\mathrm{E{-}03}$ \\
0.75&
$1.289\mathrm{E{-}03}$ & $7.894\mathrm{E{-}03}$ & $4.600\mathrm{E{-}03}$ \\
0.9&
$1.221\mathrm{E{-}03}$ & $8.352\mathrm{E{-}03}$ & $4.983\mathrm{E{-}03}$ \\
1.1&
$1.353\mathrm{E{-}03}$ & $1.550\mathrm{E{-}02}$ & $4.384\mathrm{E{-}03}$ \\
1.5&
$2.055\mathrm{E{-}03}$ & $1.421\mathrm{E{-}02}$ & $3.903\mathrm{E{-}03}$ \\
\hline
\end{tabular}
\caption{Convergence time for fminunc, fmincon, and steepest descent
\vspace{-2mm}}
\end{table}

Surprisingly enough, the Steepest descent converged very quickly, but to a wrong point. Also, there are two choices for which \texttt{fmincon} converged to the true minimizer, although they are not the closest initial guesses!\\

All in all, these results are bad enough for us to not consider the case with added white noise. We notice that the $k$ relative error is always bigger than the $r$ relative error. We decide to run our tests with an initial point ``at the true $k$".
In other words, our initial guess will have its $k$ coordinate as $10^6$. Here are the results:\\[0.25em]

\begin{table}[H]
    \caption{Optimization results for $r$ and $k$ without noise, with $k_0 = 10^6 = k_{exact}$}
\resizebox{\textwidth}{!}{%
\begin{tabular}{ccccccccc}
\toprule
\multicolumn{9}{c}{\textbf{fminunc}}\\
\midrule
$rn_0$ & iter & $r$ & $k$ & $\lvert r - r_{\text{exact}} \rvert / r_{\text{exact}}$ & $\lvert k - k_{\text{exact}} \rvert / k_{\text{exact}}$ & feval & Interp.\ error & Extrap.\ error \\
\midrule
0.25 & 3 & $1.300\mathrm{E{-}01}$ & $1.000\mathrm{E{+}06}$ & $1.976\mathrm{E{-}07}$ & $1.044\mathrm{E{-}11}$ & $1.111\mathrm{E{-}14}$ & $1.397\mathrm{E{-}01}$ & $1.148\mathrm{E{-}04}$ \\
0.5 & 7 & $1.300\mathrm{E{-}01}$ & $1.000\mathrm{E{+}06}$ & $5.466\mathrm{E{-}07}$ & $5.588\mathrm{E{-}15}$ & $8.500\mathrm{E{-}14}$ & $3.864\mathrm{E{-}01}$ & $3.046\mathrm{E{-}04}$ \\
0.75 & 6 & $1.300\mathrm{E{-}01}$ & $1.000\mathrm{E{+}06}$ & $6.360\mathrm{E{-}08}$ & $8.149\mathrm{E{-}16}$ & $1.151\mathrm{E{-}15}$ & $4.495\mathrm{E{-}02}$ & $3.544\mathrm{E{-}05}$ \\
0.9 & 5 & $1.300\mathrm{E{-}01}$ & $1.000\mathrm{E{+}06}$ & $5.630\mathrm{E{-}08}$ & $6.985\mathrm{E{-}16}$ & $9.020\mathrm{E{-}16}$ & $3.980\mathrm{E{-}02}$ & $3.138\mathrm{E{-}05}$ \\
1.1 & 4 & $1.300\mathrm{E{-}01}$ & $1.000\mathrm{E{+}06}$ & $5.918\mathrm{E{-}08}$ & $1.164\mathrm{E{-}15}$ & $9.965\mathrm{E{-}16}$ & $4.183\mathrm{E{-}02}$ & $3.298\mathrm{E{-}05}$ \\
1.5 & 4 & $1.300\mathrm{E{-}01}$ & $1.000\mathrm{E{+}06}$ & $5.717\mathrm{E{-}08}$ & $1.246\mathrm{E{-}14}$ & $9.300\mathrm{E{-}16}$ & $4.041\mathrm{E{-}02}$ & $3.187\mathrm{E{-}05}$ \\
\midrule

\multicolumn{9}{c}{\textbf{fmincon}}\\
\midrule
0.25 & 40 & $5.871\mathrm{E{-}01}$ & $9.999\mathrm{E{+}05}$ & $3.516\mathrm{E{+}00}$ & $2.631\mathrm{E{-}05}$ & $1.841\mathrm{E{-}01}$ & $5.685\mathrm{E{+}05}$ & $3.911\mathrm{E{+}01}$ \\
0.5 & 39 & $1.301\mathrm{E{-}01}$ & $9.991\mathrm{E{+}05}$ & $5.028\mathrm{E{-}04}$ & $9.291\mathrm{E{-}04}$ & $3.643\mathrm{E{-}07}$ & $7.999\mathrm{E{+}02}$ & $9.288\mathrm{E{+}02}$ \\
0.75 & 44 & $1.301\mathrm{E{-}01}$ & $9.993\mathrm{E{+}05}$ & $3.924\mathrm{E{-}04}$ & $7.251\mathrm{E{-}04}$ & $2.219\mathrm{E{-}07}$ & $6.242\mathrm{E{+}02}$ & $7.248\mathrm{E{+}02}$ \\
0.9 & 36 & $1.301\mathrm{E{-}01}$ & $9.991\mathrm{E{+}05}$ & $5.029\mathrm{E{-}04}$ & $9.292\mathrm{E{-}04}$ & $3.644\mathrm{E{-}07}$ & $7.999\mathrm{E{+}02}$ & $9.289\mathrm{E{+}02}$ \\
1.1 & 36 & $1.301\mathrm{E{-}01}$ & $9.991\mathrm{E{+}05}$ & $5.060\mathrm{E{-}04}$ & $9.341\mathrm{E{-}04}$ & $3.683\mathrm{E{-}07}$ & $8.042\mathrm{E{+}02}$ & $9.338\mathrm{E{+}02}$ \\
1.5 & 36 & $1.301\mathrm{E{-}01}$ & $9.991\mathrm{E{+}05}$ & $5.060\mathrm{E{-}04}$ & $9.341\mathrm{E{-}04}$ & $3.683\mathrm{E{-}07}$ & $8.042\mathrm{E{+}02}$ & $9.338\mathrm{E{+}02}$ \\
\midrule

\multicolumn{9}{c}{\textbf{Steepest Descent}}\\
\midrule
0.25 & 200 & $1.300\mathrm{E{-}01}$ & $1.000\mathrm{E{+}06}$ & $2.389\mathrm{E{-}04}$ & $1.746\mathrm{E{-}15}$ & $1.625\mathrm{E{-}08}$ & $1.689\mathrm{E{+}02}$ & $1.334\mathrm{E{-}01}$ \\
0.5 & 200 & $1.300\mathrm{E{-}01}$ & $1.000\mathrm{E{+}06}$ & $2.107\mathrm{E{-}04}$ & $1.513\mathrm{E{-}15}$ & $1.263\mathrm{E{-}08}$ & $1.489\mathrm{E{+}02}$ & $1.176\mathrm{E{-}01}$ \\
0.75 & 200 & $1.300\mathrm{E{-}01}$ & $1.000\mathrm{E{+}06}$ & $1.584\mathrm{E{-}04}$ & $5.821\mathrm{E{-}16}$ & $7.137\mathrm{E{-}09}$ & $1.120\mathrm{E{+}02}$ & $8.835\mathrm{E{-}02}$ \\
0.9 & 200 & $1.300\mathrm{E{-}01}$ & $1.000\mathrm{E{+}06}$ & $8.491\mathrm{E{-}05}$ & $0.000\mathrm{E{+}00}$ & $2.051\mathrm{E{-}09}$ & $6.002\mathrm{E{+}01}$ & $4.734\mathrm{E{-}02}$ \\
1.1 & 200 & $1.300\mathrm{E{-}01}$ & $1.000\mathrm{E{+}06}$ & $1.353\mathrm{E{-}04}$ & $0.000\mathrm{E{+}00}$ & $5.210\mathrm{E{-}09}$ & $9.565\mathrm{E{+}01}$ & $7.535\mathrm{E{-}02}$ \\
1.5 & 200 & $1.303\mathrm{E{-}01}$ & $1.000\mathrm{E{+}06}$ & $2.427\mathrm{E{-}03}$ & $6.985\mathrm{E{-}15}$ & $1.669\mathrm{E{-}06}$ & $1.712\mathrm{E{+}03}$ & $1.330\mathrm{E{+}00}$ \\
\midrule

\end{tabular}%
}
\end{table}

\begin{figure}[H]
    \centering
    \includegraphics[width=0.6\textwidth]{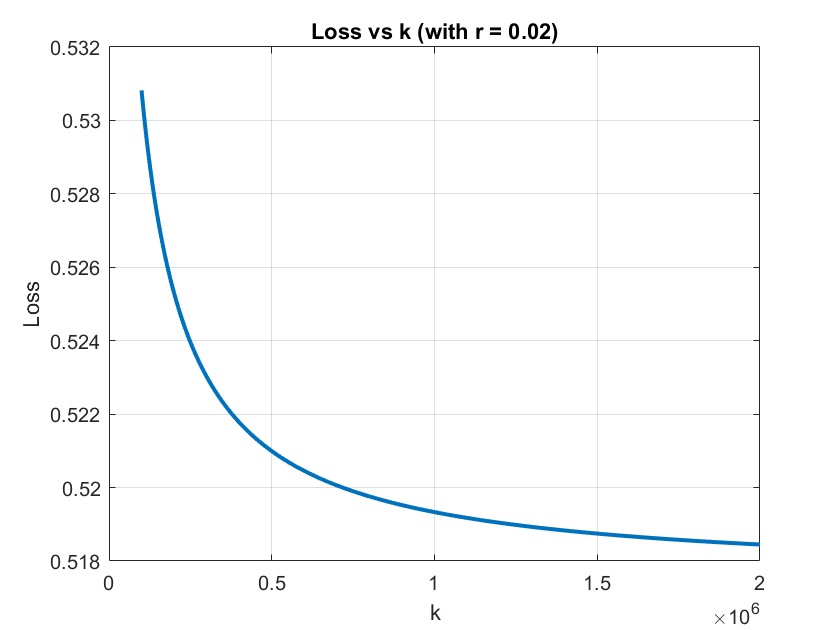}
    \caption{Loss function with $r = 0.02$ fixed and $k$ as the variable.}
    \label{fig:loss_r_fixed}
\end{figure}

\begin{table}[h!]
\centering
\begin{tabular}{|c|c|c|c|}
\hline
$nr_0$&
\textbf{fminunc} & \textbf{fmincon} & \textbf{SD} \\
\hline
0.25 &
$2.099\mathrm{E{-}03}$ & $1.921\mathrm{E{-}02}$ & $4.347\mathrm{E{-}03}$ \\
0.5 &
$1.214\mathrm{E{-}03}$ & $1.590\mathrm{E{-}02}$ & $4.314\mathrm{E{-}03}$ \\
0.75 & 
$1.501\mathrm{E{-}03}$ & $1.743\mathrm{E{-}02}$ & $4.554\mathrm{E{-}03}$ \\
0.9 & 
$1.589\mathrm{E{-}03}$ & $1.479\mathrm{E{-}02}$ & $4.728\mathrm{E{-}03}$ \\
1.1 & 
$1.089\mathrm{E{-}03}$ & $1.792\mathrm{E{-}02}$ & $4.403\mathrm{E{-}03}$ \\
1.5 & 
$1.304\mathrm{E{-}03}$ & $1.626\mathrm{E{-}02}$ & $5.578\mathrm{E{-}03}$ \\
\hline
\end{tabular}
\caption{Convergence time (in seconds) for fminunc, fmincon, and steepest descent (methods which did not converge omited)}
\end{table}
No change in runtime, however a drastic improvement in results. Every run was successful (with only the case where we started with $r_0 = 0.25 \times r^{true}$, which is a point far from our true minimizer. \\[0.5em]
To understand this behavior, we fix the value of \( r \) (chosen arbitrarily) and examine how the loss varies as a function of \( k \) (Figure \ref{fig:loss_r_fixed}).

We notice that the loss looks similar to the graph of $\exp(-x)$. Therefore, we decide to take the loss function with the logarithm of $k$ (i.e., $\Tilde{k} = \log(k)$), { where the logistic equation becomes $$\dfrac{dp(t)}{dt} = (1 - \dfrac{p(t)}{e^{\Tilde{k}}})\;r\;p(t)$$  with exact solution $p(t) =  \dfrac{e^{\Tilde{k}}\,p_0e^{r(t-t_0)}}{e^{\Tilde{k}} -p_0+p_0e^{r(t-t_0)}}$}.\\ 

Therefore, the loss became Loss($r$, log($k$)). Additionally, our initial guesses will not be chosen "at the true $k$" like before. On the contrary, the initial guesses will try multiple values for the initial $k$, exactly like with $r$.\\

We get the following results:
\begin{table}[H]
\centering
\begin{tabular}{|c|c|c|c|c|c|c|}
\hline
$nr_0$&
\textbf{fminunc} & \textbf{fmincon} & \textbf{Secant} & \textbf{NewtonS} & \textbf{NewtonE} & \textbf{Steepest Descent} \\
\hline
0.25 &
$2.274\mathrm{E}{-3}$ & $2.329\mathrm{E}{-2}$ & $5.347\mathrm{E}{-3}$ & $7.958\mathrm{E}{-3}$ & $6.390\mathrm{E}{-2}$ & $6.062\mathrm{E}{-3}$ \\
0.5 &
$3.760\mathrm{E}{-3}$ & $2.451\mathrm{E}{-2}$ & $8.170\mathrm{E}{-3}$ & $1.085\mathrm{E}{-2}$ & $1.192\mathrm{E}{-1}$ & $1.037\mathrm{E}{-2}$ \\
0.75 &
$4.996\mathrm{E}{-3}$ & $3.950\mathrm{E}{-2}$ & $1.249\mathrm{E}{-2}$ & $2.348\mathrm{E}{-2}$ & $8.248\mathrm{E}{-2}$ & $9.238\mathrm{E}{-3}$ \\
0.9 &
$2.941\mathrm{E}{-3}$ & $3.708\mathrm{E}{-2}$ & $1.215\mathrm{E}{-2}$ & $1.834\mathrm{E}{-2}$ & $8.456\mathrm{E}{-2}$ & $8.972\mathrm{E}{-3}$ \\
1.1 & 
$8.772\mathrm{E}{-4}$ & $2.645\mathrm{E}{-2}$ & $1.029\mathrm{E}{-2}$ & $1.295\mathrm{E}{-3}$ & $1.307\mathrm{E}{-2}$ & $8.784\mathrm{E}{-3}$ \\
1.5 & 
$9.005\mathrm{E}{-4}$ & $2.894\mathrm{E}{-2}$ & $9.021\mathrm{E}{-3}$ & $1.579\mathrm{E}{-2}$ & $8.211\mathrm{E}{-2}$ & $1.017\mathrm{E}{-2}$ \\
\hline
\end{tabular}
\caption{Convergence time (in seconds) for all optimization algorithms}
\end{table}
The secant, and in some cases the Newton using Symbolic differentiation failed to converge. We generally notice a fast convergence for most initial guesses. For a few points closest to the true minimizer, we have convergence to it. The steepest descent does a terrible job by (fastly) converging to the wrong point in every single run. \texttt{fmincon} surprisingly converged to the true minimizer every time, however it took more time than usual.
\begin{table}[H]
 \caption{Optimization results for $r$ and $log(k)$ without noise
 }
\resizebox{\textwidth}{!}{
\begin{tabular}{ccccccccc}
\toprule
\multicolumn{9}{c}{\textbf{fminunc}}\\
\midrule
$rn_0$ & iter & $r$ & $k$ & $\lvert r - r_{\text{exact}} \rvert / r_{\text{exact}}$ & $\lvert k - k_{\text{exact}} \rvert / k_{\text{exact}}$ & feval & Interp.\ error & Extrap.\ error \\
\midrule
0.25 & 1  & $-3.152\mathrm{E}{-5}$ & $3.172\mathrm{E}{1}$   & $1.000\mathrm{E}{0}$ & $1.000\mathrm{E}{0}$ & $4.901\mathrm{E}{-1}$ & $9.277\mathrm{E}{5}$ & $1.203\mathrm{E}{6}$ \\
0.5 & 1  & $-1.099\mathrm{E}{-3}$ & $1.048\mathrm{E}{3}$   & $1.008\mathrm{E}{0}$ & $9.990\mathrm{E}{-1}$ & $4.924\mathrm{E}{-1}$ & $9.299\mathrm{E}{5}$ & $1.156\mathrm{E}{6}$ \\
0.75 & 9  & $1.300\mathrm{E}{-1}$  & $1.000\mathrm{E}{6}$   & $1.569\mathrm{E}{-9}$ & $1.109\mathrm{E}{-7}$ & $6.269\mathrm{E}{-15}$ & $1.049\mathrm{E}{-1}$ & $1.109\mathrm{E}{-1}$ \\
0.9 & 11 & $1.300\mathrm{E}{-1}$  & $1.000\mathrm{E}{6}$   & $4.566\mathrm{E}{-9}$ & $1.064\mathrm{E}{-7}$ & $5.872\mathrm{E}{-15}$ & $1.015\mathrm{E}{-1}$ & $1.063\mathrm{E}{-1}$ \\
1.1 & 1  & $-8.570\mathrm{E}{-1}$ & $2.940\mathrm{E}{6}$   & $7.592\mathrm{E}{0}$ & $1.940\mathrm{E}{0}$ & $5.692\mathrm{E}{-1}$ & $9.998\mathrm{E}{5}$ & $9.998\mathrm{E}{5}$ \\
1.5 & 1  & $-8.050\mathrm{E}{-1}$ & $8.142\mathrm{E}{8}$   & $7.192\mathrm{E}{0}$ & $8.132\mathrm{E}{2}$ & $5.692\mathrm{E}{-1}$ & $9.998\mathrm{E}{5}$ & $9.998\mathrm{E}{5}$ \\
\midrule
\multicolumn{9}{c}{\textbf{fmincon}}\\
\midrule
0.25 & 30 & $1.301\mathrm{E}{-1}$ & $9.991\mathrm{E}{5}$  & $5.004\mathrm{E}{-4}$ & $9.238\mathrm{E}{-4}$ & $3.602\mathrm{E}{-7}$ & $7.953\mathrm{E}{2}$ & $9.234\mathrm{E}{2}$ \\
0.5 & 31 & $1.300\mathrm{E}{-1}$ & $9.993\mathrm{E}{5}$  & $3.748\mathrm{E}{-4}$ & $6.916\mathrm{E}{-4}$ & $2.018\mathrm{E}{-7}$ & $5.953\mathrm{E}{2}$ & $6.913\mathrm{E}{2}$ \\
0.75 & 25 & $1.301\mathrm{E}{-1}$ & $9.991\mathrm{E}{5}$  & $5.063\mathrm{E}{-4}$ & $9.341\mathrm{E}{-4}$ & $3.682\mathrm{E}{-7}$ & $8.041\mathrm{E}{2}$ & $9.337\mathrm{E}{2}$ \\
0.9 & 26 & $1.301\mathrm{E}{-1}$ & $9.993\mathrm{E}{5}$  & $3.924\mathrm{E}{-4}$ & $7.247\mathrm{E}{-4}$ & $2.217\mathrm{E}{-7}$ & $6.239\mathrm{E}{2}$ & $7.245\mathrm{E}{2}$ \\
1.1 & 17 & $1.301\mathrm{E}{-1}$ & $9.993\mathrm{E}{5}$  & $3.898\mathrm{E}{-4}$ & $7.191\mathrm{E}{-4}$ & $2.182\mathrm{E}{-7}$ & $6.191\mathrm{E}{2}$ & $7.188\mathrm{E}{2}$ \\
1.5 & 22 & $1.301\mathrm{E}{-1}$ & $9.993\mathrm{E}{5}$  & $3.898\mathrm{E}{-4}$ & $7.191\mathrm{E}{-4}$ & $2.182\mathrm{E}{-7}$ & $6.191\mathrm{E}{2}$ & $7.188\mathrm{E}{2}$ \\
\midrule
\multicolumn{9}{c}{\textbf{NewtonS}}\\
\midrule
0.75 & 50 & $-3.376\mathrm{E}{-1}$ & $9.889\mathrm{E}{3}$   & $3.597\mathrm{E}{0}$ & $9.901\mathrm{E}{-1}$ & $5.692\mathrm{E}{-1}$ & $9.998\mathrm{E}{5}$ & $9.998\mathrm{E}{5}$ \\
1.1 & 10 & $1.300\mathrm{E}{-1}$ & $1.000\mathrm{E}{6}$   & $2.135\mathrm{E}{-16}$ & $4.657\mathrm{E}{-16}$ & $1.116\mathrm{E}{-31}$ & $4.427\mathrm{E}{-10}$ & $4.007\mathrm{E}{-10}$ \\
1.5 & 50 & $-2.509\mathrm{E}{1}$ & $4.094\mathrm{E}{8}$   & $1.940\mathrm{E}{2}$ & $4.084\mathrm{E}{2}$ & $5.692\mathrm{E}{-1}$ & $9.998\mathrm{E}{5}$ & $9.998\mathrm{E}{5}$ \\
\midrule

\multicolumn{9}{c}{\textbf{NewtonE}}\\
\midrule
0.25 & 50 & $3.267\mathrm{E}{-2}$  & $6.078\mathrm{E}{-21}$ & $7.487\mathrm{E}{-1}$ & $1.000\mathrm{E}{0}$ & $6.554\mathrm{E}{4}$ & $6.554\mathrm{E}{4}$ & $9.998\mathrm{E}{5}$ \\
0.5 & 50 & $7.059\mathrm{E}{-2}$  & $1.898\mathrm{E}{-19}$ & $4.570\mathrm{E}{-1}$ & $1.000\mathrm{E}{0}$ & $6.554\mathrm{E}{4}$ & $6.554\mathrm{E}{4}$ & $9.998\mathrm{E}{5}$ \\
0.75 & 50 & $-2.548\mathrm{E}{0}$  & $1.027\mathrm{E}{4}$   & $2.060\mathrm{E}{1}$  & $9.897\mathrm{E}{-1}$ & $5.692\mathrm{E}{-1}$ & $9.998\mathrm{E}{5}$ & $9.998\mathrm{E}{5}$ \\
0.9 & 50 & $9.455\mathrm{E}{-2}$  & $4.096\mathrm{E}{-18}$ & $2.727\mathrm{E}{-1}$ & $1.000\mathrm{E}{0}$ & $6.554\mathrm{E}{4}$ & $6.554\mathrm{E}{4}$ & $9.998\mathrm{E}{5}$ \\
1.1 & 10 & $1.300\mathrm{E}{-1}$  & $1.000\mathrm{E}{6}$   & $2.135\mathrm{E}{-16}$ & $4.657\mathrm{E}{-16}$ & $1.116\mathrm{E}{-31}$ & $4.427\mathrm{E}{-10}$ & $4.007\mathrm{E}{-10}$ \\
1.5 & 50 & $-2.509\mathrm{E}{1}$  & $4.094\mathrm{E}{8}$   & $1.940\mathrm{E}{2}$  & $4.084\mathrm{E}{2}$ & $5.692\mathrm{E}{-1}$ & $9.998\mathrm{E}{5}$ & $9.998\mathrm{E}{5}$ \\
\midrule
\multicolumn{9}{c}{\textbf{Steepest Descent}}\\
\midrule
0.25 & 200 & $3.252\mathrm{E}{-2}$ & $1.002\mathrm{E}{4}$   & $7.499\mathrm{E}{-1}$ & $9.900\mathrm{E}{-1}$ & $5.564\mathrm{E}{-1}$ & $9.885\mathrm{E}{5}$ & $9.898\mathrm{E}{5}$ \\
0.5 & 200 & $6.500\mathrm{E}{-2}$ & $1.002\mathrm{E}{4}$   & $5.000\mathrm{E}{-1}$ & $9.900\mathrm{E}{-1}$ & $5.564\mathrm{E}{-1}$ & $9.885\mathrm{E}{5}$ & $9.898\mathrm{E}{5}$ \\
0.75 & 200 & $1.016\mathrm{E}{-1}$ & $3.186\mathrm{E}{4}$   & $2.185\mathrm{E}{-1}$ & $9.681\mathrm{E}{-1}$ & $5.297\mathrm{E}{-1}$ & $9.645\mathrm{E}{5}$ & $9.679\mathrm{E}{5}$ \\
0.9 & 200 & $1.531\mathrm{E}{-1}$ & $2.625\mathrm{E}{5}$   & $1.780\mathrm{E}{-1}$ & $7.375\mathrm{E}{-1}$ & $2.918\mathrm{E}{-1}$ & $7.158\mathrm{E}{5}$ & $7.374\mathrm{E}{5}$ \\
1.1 & 200 & $5.705\mathrm{E}{-2}$ & $3.883\mathrm{E}{6}$   & $5.612\mathrm{E}{-1}$ & $2.883\mathrm{E}{0}$ & $1.875\mathrm{E}{-1}$ & $5.739\mathrm{E}{5}$ & $2.411\mathrm{E}{6}$ \\
1.5 & 200 & $5.146\mathrm{E}{-2}$ & $8.494\mathrm{E}{8}$   & $6.041\mathrm{E}{-1}$ & $8.484\mathrm{E}{2}$ & $2.206\mathrm{E}{-1}$ & $6.224\mathrm{E}{5}$ & $7.616\mathrm{E}{7}$ \\
\midrule
\end{tabular}%
}
\end{table}

We finally choose all initial guesses to share the same \(\log(k)\) coordinate as the true minimizer. The results are comparable to the case where we were only optimizing over $r$. 
Other than the secant method and a couple of distant points in the Newton, everything converged to the true minimizer. The Newton method took more time than in the case with $r$ only. The steepest descent was $10^3$ times faster than in the case with $r$ only.\\
\begin{table}[H]
\caption{Optimization results for $r$ and $log(k)$ without noise, with $k_0 = log(10^6) = log(k_{exact})$\vspace{-1mm}}
\resizebox{\textwidth}{!}{
\begin{tabular}{ccccccccc}
\toprule
\multicolumn{9}{c}{\textbf{fminunc}}\\
\midrule
$nr_0$ & iter & $r$ & $k$ & rel.\ error $r$ & rel.\ error $k$ & feval & Interp.\ error & Extrap.\ error \\
\midrule
0.25 & 17 & $1.300\mathrm{E}{-1}$ & $1.000\mathrm{E}{6}$ & $5.090\mathrm{E}{-7}$ & $1.483\mathrm{E}{-7}$ & $6.159\mathrm{E}{-14}$ & $3.289\mathrm{E}{-1}$ & $1.482\mathrm{E}{-1}$ \\
0.5 & 11 & $1.300\mathrm{E}{-1}$ & $1.000\mathrm{E}{6}$ & $8.372\mathrm{E}{-9}$ & $9.694\mathrm{E}{-8}$ & $5.018\mathrm{E}{-15}$ & $9.387\mathrm{E}{-2}$ & $9.692\mathrm{E}{-2}$ \\
0.75 & 12 & $1.300\mathrm{E}{-1}$ & $1.000\mathrm{E}{6}$ & $6.098\mathrm{E}{-8}$ & $3.602\mathrm{E}{-7}$ & $5.986\mathrm{E}{-14}$ & $3.242\mathrm{E}{-1}$ & $3.601\mathrm{E}{-1}$ \\
0.9 & 11 & $1.300\mathrm{E}{-1}$ & $9.999\mathrm{E}{5}$ & $1.880\mathrm{E}{-7}$ & $8.244\mathrm{E}{-7}$ & $3.057\mathrm{E}{-13}$ & $7.327\mathrm{E}{-1}$ & $8.242\mathrm{E}{-1}$ \\
1.1 & 10 & $1.300\mathrm{E}{-1}$ & $1.000\mathrm{E}{6}$ & $1.905\mathrm{E}{-8}$ & $2.012\mathrm{E}{-7}$ & $2.173\mathrm{E}{-14}$ & $1.954\mathrm{E}{-1}$ & $2.011\mathrm{E}{-1}$ \\
1.5 & 6  & $1.300\mathrm{E}{-1}$ & $1.000\mathrm{E}{6}$ & $5.275\mathrm{E}{-8}$ & $1.118\mathrm{E}{-7}$ & $5.293\mathrm{E}{-15}$ & $9.641\mathrm{E}{-2}$ & $1.118\mathrm{E}{-1}$ \\
\midrule
\multicolumn{9}{c}{\textbf{fmincon}}\\
\midrule
0.25 & 27 & $1.301\mathrm{E}{-1}$ & $9.993\mathrm{E}{5}$ & $3.924\mathrm{E}{-4}$ & $7.247\mathrm{E}{-4}$ & $2.216\mathrm{E}{-7}$ & $6.239\mathrm{E}{2}$ & $7.244\mathrm{E}{2}$ \\
0.5 & 26 & $1.301\mathrm{E}{-1}$ & $9.993\mathrm{E}{5}$ & $3.924\mathrm{E}{-4}$ & $7.246\mathrm{E}{-4}$ & $2.216\mathrm{E}{-7}$ & $6.238\mathrm{E}{2}$ & $7.243\mathrm{E}{2}$ \\
0.75 & 26 & $1.301\mathrm{E}{-1}$ & $9.993\mathrm{E}{5}$ & $3.924\mathrm{E}{-4}$ & $7.247\mathrm{E}{-4}$ & $2.216\mathrm{E}{-7}$ & $6.239\mathrm{E}{2}$ & $7.244\mathrm{E}{2}$ \\
0.9 & 28 & $1.300\mathrm{E}{-1}$ & $9.994\mathrm{E}{5}$ & $3.483\mathrm{E}{-4}$ & $6.427\mathrm{E}{-4}$ & $1.743\mathrm{E}{-7}$ & $5.532\mathrm{E}{2}$ & $6.424\mathrm{E}{2}$ \\
1.1 & 17 & $1.301\mathrm{E}{-1}$ & $9.993\mathrm{E}{5}$ & $3.898\mathrm{E}{-4}$ & $7.191\mathrm{E}{-4}$ & $2.182\mathrm{E}{-7}$ & $6.191\mathrm{E}{2}$ & $7.188\mathrm{E}{2}$ \\
1.5 & 22 & $1.301\mathrm{E}{-1}$ & $9.993\mathrm{E}{5}$ & $3.898\mathrm{E}{-4}$ & $7.191\mathrm{E}{-4}$ & $2.182\mathrm{E}{-7}$ & $6.191\mathrm{E}{2}$ & $7.188\mathrm{E}{2}$ \\
\midrule
\multicolumn{9}{c}{\textbf{NewtonS}}\\
\midrule
0.25 & 8  & $-8.157\mathrm{E}{-19}$ & $1.000\mathrm{E}{4}$ & $1.000\mathrm{E}{0}$ & $9.900\mathrm{E}{-1}$ & $5.564\mathrm{E}{-1}$ & $9.885\mathrm{E}{5}$ & $9.898\mathrm{E}{5}$ \\
0.5 & 50 & $-3.041\mathrm{E}{1}$   & $1.069\mathrm{E}{10}$ & $2.349\mathrm{E}{2}$ & $1.069\mathrm{E}{4}$ & $5.692\mathrm{E}{-1}$ & $9.998\mathrm{E}{5}$ & $9.998\mathrm{E}{5}$ \\
0.75 & 5  & $1.300\mathrm{E}{-1}$   & $1.000\mathrm{E}{6}$ & $8.859\mathrm{E}{-12}$ & $2.158\mathrm{E}{-12}$ & $1.880\mathrm{E}{-23}$ & $5.745\mathrm{E}{-6}$ & $2.155\mathrm{E}{-6}$ \\
0.9 & 4  & $1.300\mathrm{E}{-1}$   & $1.000\mathrm{E}{6}$ & $2.039\mathrm{E}{-11}$ & $4.934\mathrm{E}{-12}$ & $9.961\mathrm{E}{-23}$ & $1.323\mathrm{E}{-5}$ & $4.929\mathrm{E}{-6}$ \\
1.1 & 4  & $1.300\mathrm{E}{-1}$   & $1.000\mathrm{E}{6}$ & $1.947\mathrm{E}{-9}$  & $4.721\mathrm{E}{-10}$ & $9.078\mathrm{E}{-19}$ & $1.263\mathrm{E}{-3}$ & $4.715\mathrm{E}{-4}$ \\
\midrule

\multicolumn{9}{c}{\textbf{ NewtonE}}\\
\midrule
0.25 & 8  & $-3.711\mathrm{E}{-18}$ & $1.000\mathrm{E}{4}$   & $1.000\mathrm{E}{0}$ & $9.900\mathrm{E}{-1}$ & $5.564\mathrm{E}{-1}$ & $9.885\mathrm{E}{5}$ & $9.898\mathrm{E}{5}$ \\
0.5 & 50 & $-3.041\mathrm{E}{1}$   & $1.069\mathrm{E}{10}$  & $2.349\mathrm{E}{2}$ & $1.069\mathrm{E}{4}$  & $5.692\mathrm{E}{-1}$ & $9.998\mathrm{E}{5}$ & $9.998\mathrm{E}{5}$ \\
0.75 & 5  & $1.300\mathrm{E}{-1}$   & $1.000\mathrm{E}{6}$   & $8.859\mathrm{E}{-12}$ & $2.156\mathrm{E}{-12}$ & $1.880\mathrm{E}{-23}$ & $5.745\mathrm{E}{-6}$ & $2.154\mathrm{E}{-6}$ \\
0.9 & 4  & $1.300\mathrm{E}{-1}$   & $1.000\mathrm{E}{6}$   & $2.039\mathrm{E}{-11}$ & $4.934\mathrm{E}{-12}$ & $9.961\mathrm{E}{-23}$ & $1.323\mathrm{E}{-5}$ & $4.929\mathrm{E}{-6}$ \\
1.1 & 4  & $1.300\mathrm{E}{-1}$   & $1.000\mathrm{E}{6}$   & $1.947\mathrm{E}{-9}$  & $4.721\mathrm{E}{-10}$ & $9.078\mathrm{E}{-19}$ & $1.263\mathrm{E}{-3}$ & $4.715\mathrm{E}{-4}$ \\
\midrule
\multicolumn{9}{c}{\textbf{Steepest Descent}}\\
\midrule
0.25 & 200 & $1.298\mathrm{E}{-1}$ & $1.002\mathrm{E}{6}$   & $1.499\mathrm{E}{-3}$ & $2.300\mathrm{E}{-3}$ & $2.248\mathrm{E}{-6}$ & $1.987\mathrm{E}{3}$ & $2.299\mathrm{E}{3}$ \\
0.5 & 200 & $1.298\mathrm{E}{-1}$ & $1.002\mathrm{E}{6}$   & $1.389\mathrm{E}{-3}$ & $2.148\mathrm{E}{-3}$ & $1.960\mathrm{E}{-6}$ & $1.855\mathrm{E}{3}$ & $2.147\mathrm{E}{3}$ \\
0.75 & 200 & $1.299\mathrm{E}{-1}$ & $1.001\mathrm{E}{6}$   & $9.596\mathrm{E}{-4}$ & $1.460\mathrm{E}{-3}$ & $9.067\mathrm{E}{-7}$ & $1.262\mathrm{E}{3}$ & $1.459\mathrm{E}{3}$ \\
0.9 & 200 & $1.299\mathrm{E}{-1}$ & $1.001\mathrm{E}{6}$   & $4.688\mathrm{E}{-4}$ & $6.986\mathrm{E}{-4}$ & $2.082\mathrm{E}{-7}$ & $6.046\mathrm{E}{2}$ & $6.983\mathrm{E}{2}$ \\
1.1 & 200 & $1.301\mathrm{E}{-1}$ & $9.991\mathrm{E}{5}$   & $6.259\mathrm{E}{-4}$ & $8.890\mathrm{E}{-4}$ & $3.395\mathrm{E}{-7}$ & $7.721\mathrm{E}{2}$ & $8.887\mathrm{E}{2}$ \\
1.5 & 200 & $1.308\mathrm{E}{-1}$ & $9.929\mathrm{E}{5}$   & $6.535\mathrm{E}{-3}$ & $7.123\mathrm{E}{-3}$ & $2.347\mathrm{E}{-5}$ & $6.420\mathrm{E}{3}$ & $7.120\mathrm{E}{3}$ \\
\midrule 
\end{tabular}\vspace{-5mm}
}
\end{table}

\begin{table}[H]
\centering
\begin{tabular}{|c|c|c|c|c|c|c|}
\hline
$nr_0$ & 
\textbf{fminunc} & \textbf{fmincon} & \textbf{NewtonS} & \textbf{NewtonE} & \textbf{SD} \\
\hline
0.25 & 
$3.577\mathrm{E}{-3}$ & $2.026\mathrm{E}{-2}$ & $6.091\mathrm{E}{-4}$ & $8.350\mathrm{E}{-3}$ & $6.559\mathrm{E}{-3}$ \\
0.5 & 
$1.180\mathrm{E}{-3}$ & $1.681\mathrm{E}{-2}$ & $4.829\mathrm{E}{-3}$ & $5.453\mathrm{E}{-2}$ & $5.658\mathrm{E}{-3}$ \\
0.75 & 
$1.951\mathrm{E}{-3}$ & $1.752\mathrm{E}{-2}$ & $3.064\mathrm{E}{-4}$ & $3.864\mathrm{E}{-3}$ & $5.285\mathrm{E}{-3}$ \\
0.9 & 
$1.814\mathrm{E}{-3}$ & $2.282\mathrm{E}{-2}$ & $3.232\mathrm{E}{-4}$ & $6.507\mathrm{E}{-2}$ & $5.754\mathrm{E}{-3}$ \\
1.1 & 
$1.556\mathrm{E}{-3}$ & $1.499\mathrm{E}{-2}$ & $3.094\mathrm{E}{-4}$ & $6.708\mathrm{E}{-2}$ & $5.042\mathrm{E}{-3}$ \\
1.5 & 
$1.117\mathrm{E}{-3}$ & $1.774\mathrm{E}{-2}$ & $1.873\mathrm{E}{-3}$ & $1.690\mathrm{E}{-2}$ & $6.408\mathrm{E}{-3}$ \\
\hline
\end{tabular}
\caption{Convergence time (in seconds) excluding Secant method(did not converge)}\vspace{-5mm}
\end{table}
We decide to add noise and see what happens.

\resizebox{\textwidth}{!}{
\begin{tabular}{ccccccccc}
\toprule
\multicolumn{9}{c}{\textbf{fminunc}}\\
\midrule
$nr_0$ & iter & $r$ & $k$ & rel.\ error $r$ & rel.\ error $k$ & feval & Interp.\ error & Extrap.\ error \\
\midrule
0.25 & 17 & $1.300\mathrm{E}{-1}$ & $1.000\mathrm{E}{6}$ & $5.107\mathrm{E}{-7}$ & $1.488\mathrm{E}{-7}$ & $6.253\mathrm{E}{-14}$ & $3.314\mathrm{E}{-1}$ & $1.494\mathrm{E}{-1}$ \\
0.5 & 11 & $1.300\mathrm{E}{-1}$ & $1.000\mathrm{E}{6}$ & $4.261\mathrm{E}{-9}$ & $9.808\mathrm{E}{-8}$ & $5.552\mathrm{E}{-15}$ & $9.874\mathrm{E}{-2}$ & $9.893\mathrm{E}{-2}$ \\
0.75 & 12 & $1.300\mathrm{E}{-1}$ & $1.000\mathrm{E}{6}$ & $6.128\mathrm{E}{-8}$ & $3.618\mathrm{E}{-7}$ & $6.024\mathrm{E}{-14}$ & $3.252\mathrm{E}{-1}$ & $3.630\mathrm{E}{-1}$ \\
0.9 & 11 & $1.300\mathrm{E}{-1}$ & $1.000\mathrm{E}{6}$ & $1.890\mathrm{E}{-7}$ & $8.308\mathrm{E}{-7}$ & $3.062\mathrm{E}{-13}$ & $7.333\mathrm{E}{-1}$ & $8.343\mathrm{E}{-1}$ \\
1.1 & 10 & $1.300\mathrm{E}{-1}$ & $1.000\mathrm{E}{6}$ & $2.395\mathrm{E}{-8}$ & $1.987\mathrm{E}{-7}$ & $2.218\mathrm{E}{-14}$ & $1.974\mathrm{E}{-1}$ & $1.978\mathrm{E}{-1}$ \\
1.5 & 6  & $1.300\mathrm{E}{-1}$ & $1.000\mathrm{E}{6}$ & $5.696\mathrm{E}{-8}$ & $1.120\mathrm{E}{-7}$ & $5.675\mathrm{E}{-15}$ & $9.983\mathrm{E}{-2}$ & $1.143\mathrm{E}{-1}$ \\
\midrule
\multicolumn{9}{c}{\textbf{fmincon}}\\
\midrule
0.25 & 27 & $1.301\mathrm{E}{-1}$ & $9.993\mathrm{E}{5}$ & $3.924\mathrm{E}{-4}$ & $7.247\mathrm{E}{-4}$ & $2.216\mathrm{E}{-7}$ & $6.239\mathrm{E}{2}$ & $7.244\mathrm{E}{2}$ \\
0.5 & 26 & $1.301\mathrm{E}{-1}$ & $9.993\mathrm{E}{5}$ & $3.924\mathrm{E}{-4}$ & $7.246\mathrm{E}{-4}$ & $2.216\mathrm{E}{-7}$ & $6.238\mathrm{E}{2}$ & $7.243\mathrm{E}{2}$ \\
0.75 & 26 & $1.301\mathrm{E}{-1}$ & $9.993\mathrm{E}{5}$ & $3.924\mathrm{E}{-4}$ & $7.247\mathrm{E}{-4}$ & $2.216\mathrm{E}{-7}$ & $6.239\mathrm{E}{2}$ & $7.244\mathrm{E}{2}$ \\
0.9 & 28 & $1.300\mathrm{E}{-1}$ & $9.994\mathrm{E}{5}$ & $3.483\mathrm{E}{-4}$ & $6.427\mathrm{E}{-4}$ & $1.743\mathrm{E}{-7}$ & $5.532\mathrm{E}{2}$ & $6.424\mathrm{E}{2}$ \\
1.1 & 17 & $1.301\mathrm{E}{-1}$ & $9.993\mathrm{E}{5}$ & $3.898\mathrm{E}{-4}$ & $7.191\mathrm{E}{-4}$ & $2.182\mathrm{E}{-7}$ & $6.191\mathrm{E}{2}$ & $7.188\mathrm{E}{2}$ \\
1.5 & 22 & $1.301\mathrm{E}{-1}$ & $9.993\mathrm{E}{5}$ & $3.898\mathrm{E}{-4}$ & $7.191\mathrm{E}{-4}$ & $2.182\mathrm{E}{-7}$ & $6.191\mathrm{E}{2}$ & $7.188\mathrm{E}{2}$ \\
\midrule
\multicolumn{9}{c}{\textbf{NewtonS}}\\
\midrule
0.25 & 8  & $-8.388\mathrm{E}{-19}$ & $1.000\mathrm{E}{4}$   & $1.000\mathrm{E}{0}$ & $9.900\mathrm{E}{-1}$ & $5.564\mathrm{E}{-1}$ & $9.885\mathrm{E}{5}$ & $9.898\mathrm{E}{5}$ \\
0.5 & 50 & $-3.041\mathrm{E}{1}$   & $1.069\mathrm{E}{10}$  & $2.349\mathrm{E}{2}$ & $1.069\mathrm{E}{4}$  & $5.692\mathrm{E}{-1}$ & $9.998\mathrm{E}{5}$ & $9.998\mathrm{E}{5}$ \\
0.75 & 5  & $1.300\mathrm{E}{-1}$   & $1.000\mathrm{E}{6}$   & $3.056\mathrm{E}{-10}$ & $1.523\mathrm{E}{-9}$ & $3.793\mathrm{E}{-16}$ & $2.581\mathrm{E}{-2}$ & $2.103\mathrm{E}{-2}$ \\
0.9 & 4  & $1.300\mathrm{E}{-1}$   & $1.000\mathrm{E}{6}$   & $9.634\mathrm{E}{-10}$ & $6.315\mathrm{E}{-9}$ & $4.844\mathrm{E}{-16}$ & $2.917\mathrm{E}{-2}$ & $2.359\mathrm{E}{-2}$ \\
1.1 & 4  & $1.300\mathrm{E}{-1}$   & $1.000\mathrm{E}{6}$   & $2.958\mathrm{E}{-9}$  & $2.007\mathrm{E}{-9}$ & $4.499\mathrm{E}{-16}$ & $2.811\mathrm{E}{-2}$ & $1.988\mathrm{E}{-2}$ \\
\midrule

\multicolumn{9}{c}{\textbf{NewtonE}}\\
\midrule
0.25 & 8  & $-1.901\mathrm{E}{-19}$ & $1.000\mathrm{E}{4}$   & $1.000\mathrm{E}{0}$   & $9.900\mathrm{E}{-1}$   & $5.564\mathrm{E}{-1}$   & $9.885\mathrm{E}{5}$   & $9.898\mathrm{E}{5}$ \\
0.5 & 50 & $-3.041\mathrm{E}{1}$   & $1.069\mathrm{E}{10}$  & $2.349\mathrm{E}{2}$   & $1.069\mathrm{E}{4}$    & $5.692\mathrm{E}{-1}$   & $9.998\mathrm{E}{5}$   & $9.998\mathrm{E}{5}$ \\
0.75 & 5  & $1.300\mathrm{E}{-1}$   & $1.000\mathrm{E}{6}$   & $3.056\mathrm{E}{-10}$ & $1.523\mathrm{E}{-9}$   & $3.793\mathrm{E}{-16}$  & $2.581\mathrm{E}{-2}$  & $2.103\mathrm{E}{-2}$ \\
0.9 & 4  & $1.300\mathrm{E}{-1}$   & $1.000\mathrm{E}{6}$   & $9.634\mathrm{E}{-10}$ & $6.315\mathrm{E}{-9}$   & $4.844\mathrm{E}{-16}$  & $2.917\mathrm{E}{-2}$  & $2.359\mathrm{E}{-2}$ \\
1.1 & 4  & $1.300\mathrm{E}{-1}$   & $1.000\mathrm{E}{6}$   & $2.958\mathrm{E}{-9}$  & $2.007\mathrm{E}{-9}$   & $4.499\mathrm{E}{-16}$  & $2.811\mathrm{E}{-2}$  & $1.988\mathrm{E}{-2}$ \\
\midrule
\multicolumn{9}{c}{\textbf{Steepest Descent}}\\
\midrule
0.25 & 200 & $1.298\mathrm{E}{-1}$ & $1.002\mathrm{E}{6}$   & $1.499\mathrm{E}{-3}$  & $2.300\mathrm{E}{-3}$  & $2.248\mathrm{E}{-6}$  & $1.987\mathrm{E}{3}$   & $2.299\mathrm{E}{3}$ \\
0.5 & 200 & $1.298\mathrm{E}{-1}$ & $1.002\mathrm{E}{6}$   & $1.389\mathrm{E}{-3}$  & $2.148\mathrm{E}{-3}$  & $1.960\mathrm{E}{-6}$  & $1.855\mathrm{E}{3}$   & $2.147\mathrm{E}{3}$ \\
0.75 & 200 & $1.299\mathrm{E}{-1}$ & $1.001\mathrm{E}{6}$   & $9.596\mathrm{E}{-4}$  & $1.460\mathrm{E}{-3}$  & $9.067\mathrm{E}{-7}$  & $1.262\mathrm{E}{3}$   & $1.459\mathrm{E}{3}$ \\
0.9 & 200 & $1.299\mathrm{E}{-1}$ & $1.001\mathrm{E}{6}$   & $4.688\mathrm{E}{-4}$  & $6.986\mathrm{E}{-4}$  & $2.082\mathrm{E}{-7}$  & $6.046\mathrm{E}{2}$   & $6.983\mathrm{E}{2}$ \\
1.1 & 200 & $1.301\mathrm{E}{-1}$ & $9.991\mathrm{E}{5}$   & $6.259\mathrm{E}{-4}$  & $8.890\mathrm{E}{-4}$  & $3.395\mathrm{E}{-7}$  & $7.721\mathrm{E}{2}$   & $8.887\mathrm{E}{2}$ \\
1.5 & 200 & $1.308\mathrm{E}{-1}$ & $9.929\mathrm{E}{5}$   & $6.535\mathrm{E}{-3}$  & $7.123\mathrm{E}{-3}$  & $2.347\mathrm{E}{-5}$  & $6.420\mathrm{E}{3}$   & $7.120\mathrm{E}{3}$ \\
\midrule

\end{tabular}
}
\begin{table}[h!]
\centering
\begin{tabular}{|c|c|c|c|c|c|c|}
\hline
$nr_0$ &
\textbf{fminunc} & \textbf{fmincon} & \textbf{NewtonS} & \textbf{NewtonE} & \textbf{SD} \\
\hline
0.25 &
$2.966\mathrm{E}{-3}$ & $1.828\mathrm{E}{-2}$ & $4.765\mathrm{E}{-4}$ & $7.216\mathrm{E}{-3}$ & $5.583\mathrm{E}{-3}$ \\
0.5 &
$1.145\mathrm{E}{-3}$ & $1.514\mathrm{E}{-2}$ & $5.036\mathrm{E}{-3}$ & $5.859\mathrm{E}{-2}$ & $5.025\mathrm{E}{-3}$ \\
0.75 &
$1.500\mathrm{E}{-3}$ & $1.537\mathrm{E}{-2}$ & $3.769\mathrm{E}{-4}$ & $4.394\mathrm{E}{-3}$ & $5.836\mathrm{E}{-3}$ \\
0.9 &
$1.299\mathrm{E}{-3}$ & $2.550\mathrm{E}{-2}$ & $3.912\mathrm{E}{-4}$ & $7.757\mathrm{E}{-2}$ & $5.374\mathrm{E}{-3}$ \\
1.1 &
$1.138\mathrm{E}{-3}$ & $1.298\mathrm{E}{-2}$ & $3.259\mathrm{E}{-4}$ & $8.939\mathrm{E}{-2}$ & $6.953\mathrm{E}{-3}$ \\
1.5 &
$1.216\mathrm{E}{-3}$ & $1.668\mathrm{E}{-2}$ & $4.714\mathrm{E}{-3}$ & $1.659\mathrm{E}{-2}$ & $5.421\mathrm{E}{-3}$ \\
\hline
\end{tabular}
\caption{Convergence time (in seconds) excluding Secant method (did not converge) for the case with noise}
\end{table}\\[0.5em]
The results are nearly the same as with the case with no added white noise.
\\[0.5em]
This series of tests highlighted the critical role of the initial guess. How it should be selected to ensure convergence, often near the true minimizer but not necessarily so, as observed. It also demonstrated the impact of re-parameterizing the space and how its topology influences the optimization process.

\subsection{PME testing algorithm:}\label{sec:appPME}
For the case of the PME, we simulate data using $\beta$ = 2, and Dirichlet boundary conditions and the explicit scheme. Our spatial domain is between 0 and 1, our final time is 0.2. dx = 0.02 , dt = $10^{-4}$. We then approximate the PDE solution using a forward in time and central in space (FTCS) numerical scheme. Notice that the solution here was approximated by a numerical scheme, unlike how we used an exact solution in ~\ref{sec: test3}. 
We will test fminunc, fmincon, and the steepest descent. First, consider the case with no additional white noise:\\[1em]
\begin{tabular}{ccccccc}
\toprule
$\beta_0$ & Method & $\hat{\beta}$ & feval & Interp. error & Extrap. error & Time (s) \\
\midrule
0.5 & fminunc         & $2.000\mathrm{E}{0}$ & $2.304\mathrm{E}{-3}$ & $9.404\mathrm{E}{-7}$ & $1.032\mathrm{E}{-6}$ & $6.364\mathrm{E}{-1}$ \\
0.5 & fmincon         & $2.000\mathrm{E}{0}$ & $2.304\mathrm{E}{-3}$ & $9.404\mathrm{E}{-7}$ & $1.032\mathrm{E}{-6}$ & $7.794\mathrm{E}{-1}$ \\
0.5 & SteepestDescent & $2.000\mathrm{E}{0}$ & $2.304\mathrm{E}{-3}$ & $9.404\mathrm{E}{-7}$ & $1.032\mathrm{E}{-6}$ & $5.094\mathrm{E}{0}$ \\
1   & fminunc         & $2.000\mathrm{E}{0}$ & $2.304\mathrm{E}{-3}$ & $9.404\mathrm{E}{-7}$ & $1.032\mathrm{E}{-6}$ & $3.682\mathrm{E}{-1}$ \\
1   & fmincon         & $2.000\mathrm{E}{0}$ & $2.304\mathrm{E}{-3}$ & $9.404\mathrm{E}{-7}$ & $1.032\mathrm{E}{-6}$ & $1.223\mathrm{E}{0}$ \\
1   & SteepestDescent & $2.000\mathrm{E}{0}$ & $2.304\mathrm{E}{-3}$ & $9.404\mathrm{E}{-7}$ & $1.032\mathrm{E}{-6}$ & $5.274\mathrm{E}{0}$ \\
1.5 & fminunc         & $2.000\mathrm{E}{0}$ & $2.304\mathrm{E}{-3}$ & $9.404\mathrm{E}{-7}$ & $1.032\mathrm{E}{-6}$ & $9.937\mathrm{E}{-1}$ \\
1.5 & fmincon         & $2.000\mathrm{E}{0}$ & $2.304\mathrm{E}{-3}$ & $9.404\mathrm{E}{-7}$ & $1.032\mathrm{E}{-6}$ & $1.366\mathrm{E}{0}$ \\
1.5 & SteepestDescent & $2.000\mathrm{E}{0}$ & $2.304\mathrm{E}{-3}$ & $9.404\mathrm{E}{-7}$ & $1.032\mathrm{E}{-6}$ & $5.011\mathrm{E}{0}$ \\
1.8 & fminunc         & $2.000\mathrm{E}{0}$ & $2.304\mathrm{E}{-3}$ & $9.404\mathrm{E}{-7}$ & $1.032\mathrm{E}{-6}$ & $5.896\mathrm{E}{-1}$ \\
1.8 & fmincon         & $2.000\mathrm{E}{0}$ & $2.304\mathrm{E}{-3}$ & $9.404\mathrm{E}{-7}$ & $1.032\mathrm{E}{-6}$ & $7.982\mathrm{E}{-1}$ \\
1.8 & SteepestDescent & $2.000\mathrm{E}{0}$ & $2.304\mathrm{E}{-3}$ & $9.404\mathrm{E}{-7}$ & $1.032\mathrm{E}{-6}$ & $3.598\mathrm{E}{0}$ \\
2.2 & fminunc         & $2.000\mathrm{E}{0}$ & $2.304\mathrm{E}{-3}$ & $9.404\mathrm{E}{-7}$ & $1.032\mathrm{E}{-6}$ & $8.649\mathrm{E}{-1}$ \\
2.2 & fmincon         & $2.000\mathrm{E}{0}$ & $2.304\mathrm{E}{-3}$ & $9.404\mathrm{E}{-7}$ & $1.032\mathrm{E}{-6}$ & $9.676\mathrm{E}{-1}$ \\
2.2 & SteepestDescent & $2.000\mathrm{E}{0}$ & $2.304\mathrm{E}{-3}$ & $9.404\mathrm{E}{-7}$ & $1.032\mathrm{E}{-6}$ & $2.661\mathrm{E}{0}$ \\
3   & fminunc         & $3.000\mathrm{E}{0}$ & $1.000\mathrm{E}{10}$ &                  &                  & $7.877\mathrm{E}{-3}$ \\
3   & fmincon         & $3.000\mathrm{E}{0}$ & $1.000\mathrm{E}{10}$ &                  &                  & $9.267\mathrm{E}{-3}$ \\
3   & SteepestDescent & $3.000\mathrm{E}{0}$ & $1.000\mathrm{E}{10}$ &                  &                  & $1.232\mathrm{E}{-2}$ \\
\bottomrule
\end{tabular}\\[3em]
Let us see what happens when we add a 3\% white noise.\\[1em]
\begin{tabular}{ccccccc}
\toprule
$\beta_0$ & Method & $\hat{\beta}$ & feval & Interp. error & Extrap. error & Time (s) \\
\midrule
0.5 & fminunc         & $1.996\mathrm{E}{0}$ & $2.259\mathrm{E}{0}$ & $9.219\mathrm{E}{-4}$ & $8.966\mathrm{E}{-4}$ & $6.398\mathrm{E}{-1}$ \\
0.5 & fmincon         & $1.996\mathrm{E}{0}$ & $2.259\mathrm{E}{0}$ & $9.219\mathrm{E}{-4}$ & $8.966\mathrm{E}{-4}$ & $8.360\mathrm{E}{-1}$ \\
0.5 & SteepestDescent & $1.996\mathrm{E}{0}$ & $2.259\mathrm{E}{0}$ & $9.219\mathrm{E}{-4}$ & $8.966\mathrm{E}{-4}$ & $3.538\mathrm{E}{0}$ \\
1   & fminunc         & $1.996\mathrm{E}{0}$ & $2.259\mathrm{E}{0}$ & $9.219\mathrm{E}{-4}$ & $8.966\mathrm{E}{-4}$ & $3.490\mathrm{E}{-1}$ \\
1   & fmincon         & $1.996\mathrm{E}{0}$ & $2.259\mathrm{E}{0}$ & $9.219\mathrm{E}{-4}$ & $8.966\mathrm{E}{-4}$ & $8.209\mathrm{E}{-1}$ \\
1   & SteepestDescent & $1.996\mathrm{E}{0}$ & $2.259\mathrm{E}{0}$ & $9.219\mathrm{E}{-4}$ & $8.966\mathrm{E}{-4}$ & $3.327\mathrm{E}{0}$ \\
1.5 & fminunc         & $1.996\mathrm{E}{0}$ & $2.259\mathrm{E}{0}$ & $9.219\mathrm{E}{-4}$ & $8.966\mathrm{E}{-4}$ & $6.476\mathrm{E}{-1}$ \\
1.5 & fmincon         & $1.996\mathrm{E}{0}$ & $2.259\mathrm{E}{0}$ & $9.219\mathrm{E}{-4}$ & $8.966\mathrm{E}{-4}$ & $8.466\mathrm{E}{-1}$ \\
1.5 & SteepestDescent & $1.996\mathrm{E}{0}$ & $2.259\mathrm{E}{0}$ & $9.219\mathrm{E}{-4}$ & $8.968\mathrm{E}{-4}$ & $3.382\mathrm{E}{0}$ \\
1.8 & fminunc         & $1.996\mathrm{E}{0}$ & $2.259\mathrm{E}{0}$ & $9.219\mathrm{E}{-4}$ & $8.966\mathrm{E}{-4}$ & $5.228\mathrm{E}{-1}$ \\
1.8 & fmincon         & $1.996\mathrm{E}{0}$ & $2.259\mathrm{E}{0}$ & $9.219\mathrm{E}{-4}$ & $8.966\mathrm{E}{-4}$ & $8.434\mathrm{E}{-1}$ \\
1.8 & SteepestDescent & $1.996\mathrm{E}{0}$ & $2.259\mathrm{E}{0}$ & $9.219\mathrm{E}{-4}$ & $8.966\mathrm{E}{-4}$ & $2.919\mathrm{E}{0}$ \\
2.2 & fminunc         & $1.996\mathrm{E}{0}$ & $2.259\mathrm{E}{0}$ & $9.219\mathrm{E}{-4}$ & $8.966\mathrm{E}{-4}$ & $5.558\mathrm{E}{-1}$ \\
2.2 & fmincon         & $1.996\mathrm{E}{0}$ & $2.259\mathrm{E}{0}$ & $9.219\mathrm{E}{-4}$ & $8.966\mathrm{E}{-4}$ & $7.232\mathrm{E}{-1}$ \\
2.2 & SteepestDescent & $1.996\mathrm{E}{0}$ & $2.259\mathrm{E}{0}$ & $9.219\mathrm{E}{-4}$ & $8.966\mathrm{E}{-4}$ & $2.375\mathrm{E}{0}$ \\
3   & fminunc         & $3.000\mathrm{E}{0}$ & $1.000\mathrm{E}{10}$ &                   &                   & $7.856\mathrm{E}{-3}$ \\
3   & fmincon         & $3.000\mathrm{E}{0}$ & $1.000\mathrm{E}{10}$ &                   &                   & $1.028\mathrm{E}{-2}$ \\
3   & SteepestDescent & $3.000\mathrm{E}{0}$ & $1.000\mathrm{E}{10}$ &                   &                   & $1.265\mathrm{E}{-2}$ \\
\bottomrule
\end{tabular}\\[1em]
These results reinforce the conclusions of the previous section. All three algorithms converged to the true minimizer when the initial guess was sufficiently close. However, steepest descent required significantly more time to do so. Moreover, when the initial guess was far from the true minimizer, all three methods rapidly converged to incorrect solutions and the loss functions exploded.

\section{FDM for the Heat Equation}\label{appendix:schemes}

To solve the heat equation numerically, we approximate the spatial and temporal derivatives at discrete grid points using finite-difference methods. Several discretization schemes are outlined below.

\subsection{Discretize Space}
We discretize the spatial domain $x \in [a,b]$ uniformly using a step size $h$. \\
\noindent we approximate $u_{xx}$ using central finite differences 

\[
u_{xx}(t, x_i) \approx \frac{u(t, x_{i-1}) - 2u(t,x_i) + u(t, x_{i+1})}{h^2}
\]\\
and we plug it back into the PME fo to get:
\[
u_t = \frac{u(t, x_{i-1}) - 2u(t,x_i) + u(t, x_{i+1})}{h^2}
\]
then for $i = 1,2,..., n-1$ and replacing the boundary conditions we get the following system:\\
\[ 
\frac{d}{dt} \begin{bmatrix}
    u(t,x_1)\\
    u(t,x_2)\\
    \vdots\\
    u(t,x_{n-2})\\
    u(t,x_{n-1})\\
\end{bmatrix} = \frac{1}{h^2} \begin{bmatrix}
    -2 & 1 & 0 & 0 &\cdots &0\\
    1 & -2 & 1 & 0 & \cdots &0\\
       0& 1 & -2 & 1 & 0 & \cdots \\
   & &\ddots & \ddots& \ddots \\
    0 & \cdots &  0 & 1 & -2 &1\\
    0 & 0 & \cdots & 0 & 1 & -2\\
\end{bmatrix} \begin{bmatrix}
    u(t,x_1)\\
    u(t,x_2)\\
    \vdots\\
    u(t,x_{n-2})\\
    u(t,x_{n-1})\\
\end{bmatrix}
\]
\\where we take \(
U(t,x_i) = 
\begin{bmatrix}
    u(t,x_1)\\
    u(t,x_2)\\
    \vdots\\
    u(t,x_{n-2})\\
    u(t,x_{n-1})\\
\end{bmatrix}
\)
so that  \(
\frac{d U(t,x_i)}{dt} = \frac{1}{h^2}\begin{bmatrix}
    -2 & 1 & 0 & 0 &\cdots &0\\
    1 & -2 & 1 & 0 & \cdots &0\\
       0& 1 & -2 & 1 & 0 & \cdots \\
   & &\ddots & \ddots& \ddots \\
    0 & \cdots &  0 & 1 & -2 &1\\
    0 & 0 & \cdots & 0 & 1 & -2\\
\end{bmatrix} U(t, x_i) 
\)\\
This transforms the PDE into a vectorized first order ODE that can be solved using RK schemes.

\subsection{Discritize Space and Time}

Recall that explicit schemes compute the next solution directly from known values at the current step. They are easy to implement but usually require very small time steps ($\tau$) to remain stable. In contrast, implicit schemes determine the solution by solving a system of equations at each step. While more computationally demanding, they offer much greater stability.
\\
\noindent We discretize the spatial domain $x \in [a,b]$ using a step size $h$ and the temporal domain $t \in [0,T]$ using a step size $\tau$.

\paragraph*{Explicit in Time - Forward Euler}

\noindent We first approximate $u_t$ using forward finite difference and $u_{xx}$ using central central finite difference.
\[ 
u_t(t_i,x_j)=\frac{u(t_{i+1},x_j)-u(t_i,x_j)}{\tau}\]

\[u_{xx}(t,x_j)=\frac{u(t,x_{j+1})-2u(t,x_j)+u(t,x_{j-1})}{h^2}
\]
We plug them into the PME with $m=1$ (heat equation) for $t = t_i$ to get :
\[ 
u(t_{i+1},x_j)=\tau \cdot 
\frac{u(t_i,x_{j+1})-2u(t_i,x_j)+u(t_i,x_{j-1})}{h^2}
+ u(t_i,x_j)
\]
which can be solved iteratively.

\paragraph*{Implicit in Time - Backward Euler}

\noindent We first approximate $u_t$ using backward finite difference and $u_{xx}$ using central finite difference:\\
\[ 
u_t(t_i,x_j) = \frac{u(t_{i},x_j ) - u (t_{i-1}, x_j)}{\tau}\\\]
\[ 
u_{xx}(t_i, x_j) = \frac{u(t_i, x_{j-1}) -2u(t_i,x_j) + u(t_i, x_{j+1})}{h^2}
\]\\
and we plug them into the PME with $m=1$ (heat equation)  for $t = t_i$ to get:
\[ 
u(t_{i-1},x_j) = -\frac{\tau}{h^2} u(t_i, x_{j-1}) + (1+\frac{2 \tau}{h^2}) u(t_i,x_j) - \frac{\tau}{h^2} u(t_i, x_{j+1})
\]
which leaves us with the following system to solve:\\
\[ 
\begin{bmatrix}
    (1+ \frac{2\tau}{h^2}) & -\frac{\tau}{h^2} & 0 & 0 &\cdots &0\\
    -\frac{\tau}{h^2} & (1+ \frac{2\tau}{h^2}) & -\frac{\tau}{h^2} & 0 & \cdots &0\\
   & &\ddots & \ddots& \ddots \\
    0 & \cdots &  0 & -\frac{\tau}{h^2} & (1+ \frac{2\tau}{h^2}) & -\frac{\tau}{h^2}\\
    0 & 0 & \cdots & 0 & -\frac{\tau}{h^2} & (1+ \frac{2\tau}{h^2})\\
\end{bmatrix} \begin{bmatrix}
    u(t_i, x_1)\\
    u(t_i, x_2)\\
    \vdots\\
    u(t_i, x_{n-1})\\
\end{bmatrix} = \begin{bmatrix}
    u(t_{i-1}, x_1)\\
    u(t_{i-1}, x_2)\\
    \vdots\\
    u(t_{i-1}, x_{n-1})\\
\end{bmatrix}
\]
this can be solved iteratively for \(
U_i = \begin{bmatrix}
    u(t_i, x_1)\\
    u(t_i, x_2)\\
    \vdots\\
    u(t_i, x_{n-1})\\
\end{bmatrix}
\)
given \(
U_{i-1} =  \begin{bmatrix}
    u(t_{i-1}, x_1)\\
    u(t_{i-1}, x_2)\\
    \vdots\\
    u(t_{i-1}, x_{n-1})\\
\end{bmatrix}
\) from the initial condition for i = 1.

\paragraph*{Implicit in Time - Crank-Nichelson} 
\noindent We use a central difference for the spatial discretization and a forward difference for time discretization, we get for $t = t_{i+ \frac{1}{2}}$:\\ 
\[\frac{u(t_{i+1}, x_j) - u(t_i, x_j)}{\tau} = \frac{u(t_{i+\frac{1}{2}}, x_{j-1}) - 2u(t_{i+\frac{1}{2}}, x_j) + u(t_{i+\frac{1}{2}}, x_{j+1})}{h^2} \]
We approximate 
\[ u(t_{i+\frac{1}{2}}, x_j) \approx \frac{1}{2} \left( u(t_i, x_j) + u(t_{i+1}, x_j) \right)\]
Plugging it into the previous equation, we get:

\[
    \frac{u(t_{i+1}, x_j) - u(t_i, x_j)}{\tau} = \frac{1}{2h^2}(
   [u(t_i, x_{j-1}) - 2u(t_i, x_j) + u(t_i, x_{j+1})] 
  +  [u(t_{i+1}, x_{j-1}) - 2u(t_{i+1}, x_j) + u(t_{i+1}, x_{j+1})])
\]
Rearranging, we obtain:
$$u(t_{i+1}, x_j) -  \frac{\tau}{2h^2} \left( u(t_{i+1}, x_{j-1}) - 2u(t_{i+1}, x_j) + u(t_{i+1}, x_{j+1}) \right) = u(t_i, x_j) + \frac{\tau}{2h^2} \left( u(t_i, x_{j-1}) - 2u(t_i, x_j) + u(t_i, x_{j+1}) \right) $$
Converting into matrix form, we get:

\begin{equation}
    \begin{bmatrix}
    1 + \frac{\tau}{h^2} & -\frac{\tau}{2h^2} & 0 & \dots & 0 \\
    -\frac{\tau}{2h^2} & 1 + \frac{\tau}{h^2} & -\frac{\tau}{2h^2} & \dots & 0 \\
    0 & -\frac{\tau}{2h^2} & 1 + \frac{\tau}{h^2} & \dots & 0 \\
    \vdots & \vdots & \vdots & \ddots & \vdots \\
    0 & 0 & 0 & \dots & 1 + \frac{\tau}{h^2}
    \end{bmatrix}
    \begin{bmatrix}
    u(t_{i+1}, x_1) \\
    u(t_{i+1}, x_2) \\
    \vdots \\
    u(t_{i+1}, x_{n-1})
    \end{bmatrix}
    =
    \begin{bmatrix}
    1 - \frac{\tau}{h^2} & \frac{\tau}{2h^2} & 0 & \dots & 0 \\
    \frac{\tau}{2h^2} & 1 - \frac{\tau}{h^2} & \frac{\tau}{2h^2} & \dots & 0 \\
    0 & \frac{\tau}{2h^2} & 1 - \frac{\tau}{h^2} & \dots & 0 \\
    \vdots & \vdots & \vdots & \ddots & \vdots \\
    0 & 0 & 0 & \dots & 1 - \frac{\tau}{h^2}
    \end{bmatrix}
    \begin{bmatrix}
   u(t_{i}, x_1) \\
    u(t_{i}, x_2) \\
    \vdots \\
    u(t_{i}, x_{n-1})
    \end{bmatrix}
\end{equation}
\end{document}